\input amstex
%

\def\next{AMS-SEKR}\ifx\styname\next \endinput\fi
\catcode`\@=11
\def\styname{AMS-SEKR}
\def\styversion{2.0}
{\W@{}\W@{\styname.STY - Version \styversion}\W@{}}
\hyphenation{acad-e-my acad-e-mies af-ter-thought anom-aly anom-alies
an-ti-deriv-a-tive an-tin-o-my an-tin-o-mies apoth-e-o-ses apoth-e-o-sis
ap-pen-dix ar-che-typ-al as-sign-a-ble as-sist-ant-ship as-ymp-tot-ic
asyn-chro-nous at-trib-uted at-trib-ut-able bank-rupt bank-rupt-cy
bi-dif-fer-en-tial blue-print busier busiest cat-a-stroph-ic
cat-a-stroph-i-cally con-gress cross-hatched data-base de-fin-i-tive
de-riv-a-tive dis-trib-ute dri-ver dri-vers eco-nom-ics econ-o-mist
elit-ist equi-vari-ant ex-quis-ite ex-tra-or-di-nary flow-chart
for-mi-da-ble forth-right friv-o-lous ge-o-des-ic ge-o-det-ic geo-met-ric
griev-ance griev-ous griev-ous-ly hexa-dec-i-mal ho-lo-no-my ho-mo-thetic
ideals idio-syn-crasy in-fin-ite-ly in-fin-i-tes-i-mal ir-rev-o-ca-ble
key-stroke lam-en-ta-ble light-weight mal-a-prop-ism man-u-script
mar-gin-al meta-bol-ic me-tab-o-lism meta-lan-guage me-trop-o-lis
met-ro-pol-i-tan mi-nut-est mol-e-cule mono-chrome mono-pole mo-nop-oly
mono-spline mo-not-o-nous mul-ti-fac-eted mul-ti-plic-able non-euclid-ean
non-iso-mor-phic non-smooth par-a-digm par-a-bol-ic pa-rab-o-loid
pa-ram-e-trize para-mount pen-ta-gon phe-nom-e-non post-script pre-am-ble
pro-ce-dur-al pro-hib-i-tive pro-hib-i-tive-ly pseu-do-dif-fer-en-tial
pseu-do-fi-nite pseu-do-nym qua-drat-ics quad-ra-ture qua-si-smooth
qua-si-sta-tion-ary qua-si-tri-an-gu-lar quin-tes-sence quin-tes-sen-tial
re-arrange-ment rec-tan-gle ret-ri-bu-tion retro-fit retro-fit-ted
right-eous right-eous-ness ro-bot ro-bot-ics sched-ul-ing se-mes-ter
semi-def-i-nite semi-ho-mo-thet-ic set-up se-vere-ly side-step sov-er-eign
spe-cious spher-oid spher-oid-al star-tling star-tling-ly
sta-tis-tics sto-chas-tic straight-est strange-ness strat-a-gem strong-hold
sum-ma-ble symp-to-matic syn-chro-nous topo-graph-i-cal tra-vers-a-ble
tra-ver-sal tra-ver-sals treach-ery turn-around un-at-tached un-err-ing-ly
white-space wide-spread wing-spread wretch-ed wretch-ed-ly Brown-ian
Eng-lish Euler-ian Feb-ru-ary Gauss-ian Grothen-dieck Hamil-ton-ian
Her-mit-ian Jan-u-ary Japan-ese Kor-te-weg Le-gendre Lip-schitz
Lip-schitz-ian Mar-kov-ian Noe-ther-ian No-vem-ber Rie-mann-ian
Schwarz-schild Sep-tem-ber
form per-iods Uni-ver-si-ty cri-ti-sism for-ma-lism}
\Invalid@\nofrills
\Invalid@\usualspace
\newif\ifnofrills@
\def\nofrills@#1#2{\relaxnext@
  \DN@{\ifx\next\nofrills
    \nofrills@true\let#2\relax\DN@\nofrills{\nextii@}%
  \else
    \nofrills@false\def#2{#1}\let\next@\nextii@\fi
\next@}}
\def\usualspace@#1{\ifnofrills@\def\usualspace{#1}\fi}
\def\addto#1#2{\csname \expandafter\eat@\string#1@\endcsname
  \expandafter{\the\csname \expandafter\eat@\string#1@\endcsname#2}}
\newdimen\bigsize@
\def\big@#1#2{{\hbox{$\left#2\vcenter to#1\bigsize@{}%
  \right.\nulldelimiterspace\z@\m@th$}}}
\def\big{\big@\@ne}
\def\Big{\big@{1.5}}
\def\bigg{\big@\tw@}
\def\Bigg{\big@{2.5}}
\def\raggedcenter@{\leftskip\z@ plus.4\hsize \rightskip\leftskip
 \parfillskip\z@ \parindent\z@ \spaceskip.3333em \xspaceskip.5em
 \pretolerance9999\tolerance9999 \exhyphenpenalty\@M
 \hyphenpenalty\@M \let\\\linebreak}
\def\upperspecialchars{\def\ss{SS}\let\i=I\let\j=J\let\ae\AE\let\oe\OE
  \let\o\O\let\aa\AA\let\l\L}
\def\uppercasetext@#1{%
  {\spaceskip1.2\fontdimen2\the\font plus1.2\fontdimen3\the\font
   \upperspecialchars\uctext@#1$\m@th\aftergroup\eat@$}}
\def\uctext@#1$#2${\endash@#1-\endash@$#2$\uctext@}
\def\endash@#1-#2\endash@{\uppercase{#1}\if\notempty{#2}--\endash@#2\endash@\fi}
\def\runaway@#1{\DN@{#1}\ifx\envir@\next@
  \Err@{You seem to have a missing or misspelled \string\end#1 ...}%
  \let\envir@\empty\fi}
\newif\iftemp@
\def\notempty#1{TT\fi\def\test@{#1}\ifx\test@\empty\temp@false
  \else\temp@true\fi \iftemp@}
\font@\tensmc=cmcsc10
\font@\sevenex=cmex7
\font@\sevenit=cmti7
\font@\eightrm=cmr8 
\font@\sixrm=cmr6 
\font@\eighti=cmmi8     \skewchar\eighti='177 
\font@\sixi=cmmi6       \skewchar\sixi='177   
\font@\eightsy=cmsy8    \skewchar\eightsy='60 
\font@\sixsy=cmsy6      \skewchar\sixsy='60   
\font@\eightex=cmex8
\font@\eightbf=cmbx8 
\font@\sixbf=cmbx6   
\font@\eightit=cmti8 
\font@\eightsl=cmsl8 
\font@\eightsmc=cmcsc8
\font@\eighttt=cmtt8 


\loadmsam
\loadmsbm
\loadeufm
\UseAMSsymbols
\newtoks\tenpoint@
\def\tenpoint{\normalbaselineskip12\p@
 \abovedisplayskip12\p@ plus3\p@ minus9\p@
 \belowdisplayskip\abovedisplayskip
 \abovedisplayshortskip\z@ plus3\p@
 \belowdisplayshortskip7\p@ plus3\p@ minus4\p@
 \textonlyfont@\rm\tenrm \textonlyfont@\it\tenit
 \textonlyfont@\sl\tensl \textonlyfont@\bf\tenbf
 \textonlyfont@\smc\tensmc \textonlyfont@\tt\tentt
 \textonlyfont@\bsmc\tenbsmc
 \ifsyntax@ \def\big##1{{\hbox{$\left##1\right.$}}}%
  \let\Big\big \let\bigg\big \let\Bigg\big
 \else
  \textfont\z@=\tenrm  \scriptfont\z@=\sevenrm  \scriptscriptfont\z@=\fiverm
  \textfont\@ne=\teni  \scriptfont\@ne=\seveni  \scriptscriptfont\@ne=\fivei
  \textfont\tw@=\tensy \scriptfont\tw@=\sevensy \scriptscriptfont\tw@=\fivesy
  \textfont\thr@@=\tenex \scriptfont\thr@@=\sevenex
        \scriptscriptfont\thr@@=\sevenex
  \textfont\itfam=\tenit \scriptfont\itfam=\sevenit
        \scriptscriptfont\itfam=\sevenit
  \textfont\bffam=\tenbf \scriptfont\bffam=\sevenbf
        \scriptscriptfont\bffam=\fivebf
  \setbox\strutbox\hbox{\vrule height8.5\p@ depth3.5\p@ width\z@}%
  \setbox\strutbox@\hbox{\lower.5\normallineskiplimit\vbox{%
        \kern-\normallineskiplimit\copy\strutbox}}%
 \setbox\z@\vbox{\hbox{$($}\kern\z@}\bigsize@=1.2\ht\z@
 \fi
 \normalbaselines\rm\ex@.2326ex\jot3\ex@\the\tenpoint@}
\newtoks\eightpoint@
\def\eightpoint{\normalbaselineskip10\p@
 \abovedisplayskip10\p@ plus2.4\p@ minus7.2\p@
 \belowdisplayskip\abovedisplayskip
 \abovedisplayshortskip\z@ plus2.4\p@
 \belowdisplayshortskip5.6\p@ plus2.4\p@ minus3.2\p@
 \textonlyfont@\rm\eightrm \textonlyfont@\it\eightit
 \textonlyfont@\sl\eightsl \textonlyfont@\bf\eightbf
 \textonlyfont@\smc\eightsmc \textonlyfont@\tt\eighttt
 \textonlyfont@\bsmc\eightbsmc
 \ifsyntax@\def\big##1{{\hbox{$\left##1\right.$}}}%
  \let\Big\big \let\bigg\big \let\Bigg\big
 \else
  \textfont\z@=\eightrm \scriptfont\z@=\sixrm \scriptscriptfont\z@=\fiverm
  \textfont\@ne=\eighti \scriptfont\@ne=\sixi \scriptscriptfont\@ne=\fivei
  \textfont\tw@=\eightsy \scriptfont\tw@=\sixsy \scriptscriptfont\tw@=\fivesy
  \textfont\thr@@=\eightex \scriptfont\thr@@=\sevenex
   \scriptscriptfont\thr@@=\sevenex
  \textfont\itfam=\eightit \scriptfont\itfam=\sevenit
   \scriptscriptfont\itfam=\sevenit
  \textfont\bffam=\eightbf \scriptfont\bffam=\sixbf
   \scriptscriptfont\bffam=\fivebf
 \setbox\strutbox\hbox{\vrule height7\p@ depth3\p@ width\z@}%
 \setbox\strutbox@\hbox{\raise.5\normallineskiplimit\vbox{%
   \kern-\normallineskiplimit\copy\strutbox}}%
 \setbox\z@\vbox{\hbox{$($}\kern\z@}\bigsize@=1.2\ht\z@
 \fi
 \normalbaselines\eightrm\ex@.2326ex\jot3\ex@\the\eightpoint@}

\font@\twelverm=cmr10 scaled\magstep1
\font@\twelveit=cmti10 scaled\magstep1
\font@\twelvesl=cmsl10 scaled\magstep1
\font@\twelvesmc=cmcsc10 scaled\magstep1
\font@\twelvett=cmtt10 scaled\magstep1
\font@\twelvebf=cmbx10 scaled\magstep1
\font@\twelvei=cmmi10 scaled\magstep1
\font@\twelvesy=cmsy10 scaled\magstep1
\font@\twelveex=cmex10 scaled\magstep1
\font@\twelvemsa=msam10 scaled\magstep1
\font@\twelveeufm=eufm10 scaled\magstep1
\font@\twelvemsb=msbm10 scaled\magstep1
\newtoks\twelvepoint@
\def\twelvepoint{\normalbaselineskip15\p@
 \abovedisplayskip15\p@ plus3.6\p@ minus10.8\p@
 \belowdisplayskip\abovedisplayskip
 \abovedisplayshortskip\z@ plus3.6\p@
 \belowdisplayshortskip8.4\p@ plus3.6\p@ minus4.8\p@
 \textonlyfont@\rm\twelverm \textonlyfont@\it\twelveit
 \textonlyfont@\sl\twelvesl \textonlyfont@\bf\twelvebf
 \textonlyfont@\smc\twelvesmc \textonlyfont@\tt\twelvett
 \textonlyfont@\bsmc\twelvebsmc
 \ifsyntax@ \def\big##1{{\hbox{$\left##1\right.$}}}%
  \let\Big\big \let\bigg\big \let\Bigg\big
 \else
  \textfont\z@=\twelverm  \scriptfont\z@=\tenrm  \scriptscriptfont\z@=\sevenrm
  \textfont\@ne=\twelvei  \scriptfont\@ne=\teni  \scriptscriptfont\@ne=\seveni
  \textfont\tw@=\twelvesy \scriptfont\tw@=\tensy \scriptscriptfont\tw@=\sevensy
  \textfont\thr@@=\twelveex \scriptfont\thr@@=\tenex
        \scriptscriptfont\thr@@=\tenex
  \textfont\itfam=\twelveit \scriptfont\itfam=\tenit
        \scriptscriptfont\itfam=\tenit
  \textfont\bffam=\twelvebf \scriptfont\bffam=\tenbf
        \scriptscriptfont\bffam=\sevenbf
  \setbox\strutbox\hbox{\vrule height10.2\p@ depth4.2\p@ width\z@}%
  \setbox\strutbox@\hbox{\lower.6\normallineskiplimit\vbox{%
        \kern-\normallineskiplimit\copy\strutbox}}%
 \setbox\z@\vbox{\hbox{$($}\kern\z@}\bigsize@=1.4\ht\z@
 \fi
 \normalbaselines\rm\ex@.2326ex\jot3.6\ex@\the\twelvepoint@}

\def\headfonts{\twelvepoint\bf}
\def\rrm{\twelvepoint\rm}

\font@\fourteenrm=cmr10 scaled\magstep2
\font@\fourteenit=cmti10 scaled\magstep2
\font@\fourteensl=cmsl10 scaled\magstep2
\font@\fourteensmc=cmcsc10 scaled\magstep2
\font@\fourteentt=cmtt10 scaled\magstep2
\font@\fourteenbf=cmbx10 scaled\magstep2
\font@\fourteeni=cmmi10 scaled\magstep2
\font@\fourteensy=cmsy10 scaled\magstep2
\font@\fourteenex=cmex10 scaled\magstep2
\font@\fourteenmsa=msam10 scaled\magstep2
\font@\fourteeneufm=eufm10 scaled\magstep2
\font@\fourteenmsb=msbm10 scaled\magstep2
\newtoks\fourteenpoint@
\def\fourteenpoint{\normalbaselineskip15\p@
 \abovedisplayskip18\p@ plus4.3\p@ minus12.9\p@
 \belowdisplayskip\abovedisplayskip
 \abovedisplayshortskip\z@ plus4.3\p@
 \belowdisplayshortskip10.1\p@ plus4.3\p@ minus5.8\p@
 \textonlyfont@\rm\fourteenrm \textonlyfont@\it\fourteenit
 \textonlyfont@\sl\fourteensl \textonlyfont@\bf\fourteenbf
 \textonlyfont@\smc\fourteensmc \textonlyfont@\tt\fourteentt
 \textonlyfont@\bsmc\fourteenbsmc
 \ifsyntax@ \def\big##1{{\hbox{$\left##1\right.$}}}%
  \let\Big\big \let\bigg\big \let\Bigg\big
 \else
  \textfont\z@=\fourteenrm  \scriptfont\z@=\twelverm  \scriptscriptfont\z@=\tenrm
  \textfont\@ne=\fourteeni  \scriptfont\@ne=\twelvei  \scriptscriptfont\@ne=\teni
  \textfont\tw@=\fourteensy \scriptfont\tw@=\twelvesy \scriptscriptfont\tw@=\tensy
  \textfont\thr@@=\fourteenex \scriptfont\thr@@=\twelveex
        \scriptscriptfont\thr@@=\twelveex
  \textfont\itfam=\fourteenit \scriptfont\itfam=\twelveit
        \scriptscriptfont\itfam=\twelveit
  \textfont\bffam=\fourteenbf \scriptfont\bffam=\twelvebf
        \scriptscriptfont\bffam=\tenbf
  \setbox\strutbox\hbox{\vrule height12.2\p@ depth5\p@ width\z@}%
  \setbox\strutbox@\hbox{\lower.72\normallineskiplimit\vbox{%
        \kern-\normallineskiplimit\copy\strutbox}}%
 \setbox\z@\vbox{\hbox{$($}\kern\z@}\bigsize@=1.7\ht\z@
 \fi
 \normalbaselines\rm\ex@.2326ex\jot4.3\ex@\the\fourteenpoint@}

\def\chapheadfonts{\fourteenpoint\bf}

\font@\seventeenrm=cmr10 scaled\magstep3
\font@\seventeenit=cmti10 scaled\magstep3
\font@\seventeensl=cmsl10 scaled\magstep3
\font@\seventeensmc=cmcsc10 scaled\magstep3
\font@\seventeentt=cmtt10 scaled\magstep3
\font@\seventeenbf=cmbx10 scaled\magstep3
\font@\seventeeni=cmmi10 scaled\magstep3
\font@\seventeensy=cmsy10 scaled\magstep3
\font@\seventeenex=cmex10 scaled\magstep3
\font@\seventeenmsa=msam10 scaled\magstep3
\font@\seventeeneufm=eufm10 scaled\magstep3
\font@\seventeenmsb=msbm10 scaled\magstep3
\newtoks\seventeenpoint@
\def\seventeenpoint{\normalbaselineskip18\p@
 \abovedisplayskip21.6\p@ plus5.2\p@ minus15.4\p@
 \belowdisplayskip\abovedisplayskip
 \abovedisplayshortskip\z@ plus5.2\p@
 \belowdisplayshortskip12.1\p@ plus5.2\p@ minus7\p@
 \textonlyfont@\rm\seventeenrm \textonlyfont@\it\seventeenit
 \textonlyfont@\sl\seventeensl \textonlyfont@\bf\seventeenbf
 \textonlyfont@\smc\seventeensmc \textonlyfont@\tt\seventeentt
 \textonlyfont@\bsmc\seventeenbsmc
 \ifsyntax@ \def\big##1{{\hbox{$\left##1\right.$}}}%
  \let\Big\big \let\bigg\big \let\Bigg\big
 \else
  \textfont\z@=\seventeenrm  \scriptfont\z@=\fourteenrm  \scriptscriptfont\z@=\twelverm
  \textfont\@ne=\seventeeni  \scriptfont\@ne=\fourteeni  \scriptscriptfont\@ne=\twelvei
  \textfont\tw@=\seventeensy \scriptfont\tw@=\fourteensy \scriptscriptfont\tw@=\twelvesy
  \textfont\thr@@=\seventeenex \scriptfont\thr@@=\fourteenex
        \scriptscriptfont\thr@@=\fourteenex
  \textfont\itfam=\seventeenit \scriptfont\itfam=\fourteenit
        \scriptscriptfont\itfam=\fourteenit
  \textfont\bffam=\seventeenbf \scriptfont\bffam=\fourteenbf
        \scriptscriptfont\bffam=\twelvebf
  \setbox\strutbox\hbox{\vrule height14.6\p@ depth6\p@ width\z@}%
  \setbox\strutbox@\hbox{\lower.86\normallineskiplimit\vbox{%
        \kern-\normallineskiplimit\copy\strutbox}}%
 \setbox\z@\vbox{\hbox{$($}\kern\z@}\bigsize@=2\ht\z@
 \fi
 \normalbaselines\rm\ex@.2326ex\jot5.2\ex@\the\seventeenpoint@}

\font@\rrrrrm=cmr10 scaled\magstep4
\font@\bigtitlefont=cmbx10 scaled\magstep4

\parindent1pc
\normallineskiplimit\p@
\newdimen\indenti \indenti=2pc
\def\pageheight#1{\vsize#1}
\def\pagewidth#1{\hsize#1%
   \captionwidth@\hsize \advance\captionwidth@-2\indenti}
\pagewidth{30pc} \pageheight{47pc}
\def\topmatter{%
 \ifx\undefined\msafam
 \else\font@\eightmsa=msam8 \font@\sixmsa=msam6
   \ifsyntax@\else \addto\tenpoint{\textfont\msafam=\tenmsa
              \scriptfont\msafam=\sevenmsa \scriptscriptfont\msafam=\fivemsa}%
     \addto\eightpoint{\textfont\msafam=\eightmsa \scriptfont\msafam=\sixmsa
              \scriptscriptfont\msafam=\fivemsa}%
   \fi
 \fi
 \ifx\undefined\msbfam
 \else\font@\eightmsb=msbm8 \font@\sixmsb=msbm6
   \ifsyntax@\else \addto\tenpoint{\textfont\msbfam=\tenmsb
         \scriptfont\msbfam=\sevenmsb \scriptscriptfont\msbfam=\fivemsb}%
     \addto\eightpoint{\textfont\msbfam=\eightmsb \scriptfont\msbfam=\sixmsb
         \scriptscriptfont\msbfam=\fivemsb}%
   \fi
 \fi
 \ifx\undefined\eufmfam
 \else \font@\eighteufm=eufm8 \font@\sixeufm=eufm6
   \ifsyntax@\else \addto\tenpoint{\textfont\eufmfam=\teneufm
       \scriptfont\eufmfam=\seveneufm \scriptscriptfont\eufmfam=\fiveeufm}%
     \addto\eightpoint{\textfont\eufmfam=\eighteufm
       \scriptfont\eufmfam=\sixeufm \scriptscriptfont\eufmfam=\fiveeufm}%
   \fi
 \fi
 \ifx\undefined\eufbfam
 \else \font@\eighteufb=eufb8 \font@\sixeufb=eufb6
   \ifsyntax@\else \addto\tenpoint{\textfont\eufbfam=\teneufb
      \scriptfont\eufbfam=\seveneufb \scriptscriptfont\eufbfam=\fiveeufb}%
    \addto\eightpoint{\textfont\eufbfam=\eighteufb
      \scriptfont\eufbfam=\sixeufb \scriptscriptfont\eufbfam=\fiveeufb}%
   \fi
 \fi
 \ifx\undefined\eusmfam
 \else \font@\eighteusm=eusm8 \font@\sixeusm=eusm6
   \ifsyntax@\else \addto\tenpoint{\textfont\eusmfam=\teneusm
       \scriptfont\eusmfam=\seveneusm \scriptscriptfont\eusmfam=\fiveeusm}%
     \addto\eightpoint{\textfont\eusmfam=\eighteusm
       \scriptfont\eusmfam=\sixeusm \scriptscriptfont\eusmfam=\fiveeusm}%
   \fi
 \fi
 \ifx\undefined\eusbfam
 \else \font@\eighteusb=eusb8 \font@\sixeusb=eusb6
   \ifsyntax@\else \addto\tenpoint{\textfont\eusbfam=\teneusb
       \scriptfont\eusbfam=\seveneusb \scriptscriptfont\eusbfam=\fiveeusb}%
     \addto\eightpoint{\textfont\eusbfam=\eighteusb
       \scriptfont\eusbfam=\sixeusb \scriptscriptfont\eusbfam=\fiveeusb}%
   \fi
 \fi
 \ifx\undefined\eurmfam
 \else \font@\eighteurm=eurm8 \font@\sixeurm=eurm6
   \ifsyntax@\else \addto\tenpoint{\textfont\eurmfam=\teneurm
       \scriptfont\eurmfam=\seveneurm \scriptscriptfont\eurmfam=\fiveeurm}%
     \addto\eightpoint{\textfont\eurmfam=\eighteurm
       \scriptfont\eurmfam=\sixeurm \scriptscriptfont\eurmfam=\fiveeurm}%
   \fi
 \fi
 \ifx\undefined\eurbfam
 \else \font@\eighteurb=eurb8 \font@\sixeurb=eurb6
   \ifsyntax@\else \addto\tenpoint{\textfont\eurbfam=\teneurb
       \scriptfont\eurbfam=\seveneurb \scriptscriptfont\eurbfam=\fiveeurb}%
    \addto\eightpoint{\textfont\eurbfam=\eighteurb
       \scriptfont\eurbfam=\sixeurb \scriptscriptfont\eurbfam=\fiveeurb}%
   \fi
 \fi
 \ifx\undefined\cmmibfam
 \else \font@\eightcmmib=cmmib8 \font@\sixcmmib=cmmib6
   \ifsyntax@\else \addto\tenpoint{\textfont\cmmibfam=\tencmmib
       \scriptfont\cmmibfam=\sevencmmib \scriptscriptfont\cmmibfam=\fivecmmib}%
    \addto\eightpoint{\textfont\cmmibfam=\eightcmmib
       \scriptfont\cmmibfam=\sixcmmib \scriptscriptfont\cmmibfam=\fivecmmib}%
   \fi
 \fi
 \ifx\undefined\cmbsyfam
 \else \font@\eightcmbsy=cmbsy8 \font@\sixcmbsy=cmbsy6
   \ifsyntax@\else \addto\tenpoint{\textfont\cmbsyfam=\tencmbsy
      \scriptfont\cmbsyfam=\sevencmbsy \scriptscriptfont\cmbsyfam=\fivecmbsy}%
    \addto\eightpoint{\textfont\cmbsyfam=\eightcmbsy
      \scriptfont\cmbsyfam=\sixcmbsy \scriptscriptfont\cmbsyfam=\fivecmbsy}%
   \fi
 \fi
 \let\topmatter\relax}
\def\chapterno@{\uppercase\expandafter{\romannumeral\chaptercount@}}
\newcount\chaptercount@
\def\chapter{\nofrills@{\afterassignment\chapterno@
                        CHAPTER \global\chaptercount@=}\chapter@
 \DNii@##1{\leavevmode\hskip-\leftskip
   \rlap{\vbox to\z@{\vss\centerline{\eightpoint
   \chapter@##1\unskip}\baselineskip2pc\null}}\hskip\leftskip
   \nofrills@false}%
 \FN@\next@}
\newbox\titlebox@

\def\title{\nofrills@{\relax}\title@%
 \DNii@##1\endtitle{\global\setbox\titlebox@\vtop{\tenpoint\bf
 \raggedcenter@\ignorespaces
 \baselineskip1.3\baselineskip\title@{##1}\endgraf}%
 \ifmonograph@ \edef\next{\the\leftheadtoks}\ifx\next\empty
    \leftheadtext{##1}\fi
 \fi
 \edef\next{\the\rightheadtoks}\ifx\next\empty \rightheadtext{##1}\fi
 }\FN@\next@}
\newbox\authorbox@
\def\author#1\endauthor{\global\setbox\authorbox@
 \vbox{\tenpoint\smc\raggedcenter@\ignorespaces
 #1\endgraf}\relaxnext@ \edef\next{\the\leftheadtoks}%
 \ifx\next\empty\leftheadtext{#1}\fi}
\newbox\affilbox@
\def\affil#1\endaffil{\global\setbox\affilbox@
 \vbox{\tenpoint\raggedcenter@\ignorespaces#1\endgraf}}
\newcount\addresscount@
\addresscount@\z@
\def\address#1\endaddress{\global\advance\addresscount@\@ne
  \expandafter\gdef\csname address\number\addresscount@\endcsname
  {\vskip12\p@ minus6\p@\noindent\eightpoint\smc\ignorespaces#1\par}}
\def\email{\nofrills@{\eightpoint{\it E-mail\/}:\enspace}\email@
  \DNii@##1\endemail{%
  \expandafter\gdef\csname email\number\addresscount@\endcsname
  {\def\usualspace{{\it\enspace}}\smallskip\noindent\eightpoint\email@
  \ignorespaces##1\par}}%
 \FN@\next@}
\def\thedate@{}
\def\date#1\enddate{\gdef\thedate@{\tenpoint\ignorespaces#1\unskip}}
\def\thethanks@{}
\def\thanks#1\endthanks{\gdef\thethanks@{\eightpoint\ignorespaces#1.\unskip}}
\def\thekeywords@{}
\def\keywords{\nofrills@{{\it Key words and phrases.\enspace}}\keywords@
 \DNii@##1\endkeywords{\def\thekeywords@{\def\usualspace{{\it\enspace}}%
 \eightpoint\keywords@\ignorespaces##1\unskip.}}%
 \FN@\next@}
\def\thesubjclass@{}
\def\subjclass{\nofrills@{{\rm2020 {\it Mathematics Subject
   Classification\/}.\enspace}}\subjclass@
 \DNii@##1\endsubjclass{\def\thesubjclass@{\def\usualspace
  {{\rm\enspace}}\eightpoint\subjclass@\ignorespaces##1\unskip.}}%
 \FN@\next@}
\newbox\abstractbox@
\def\abstract{\nofrills@{{\smc Abstract.\enspace}}\abstract@
 \DNii@{\setbox\abstractbox@\vbox\bgroup\noindent$$\vbox\bgroup
  \def\envir@{abstract}\advance\hsize-2\indenti
  \usualspace@{{\enspace}}\eightpoint \noindent\abstract@\ignorespaces}%
 \FN@\next@}
\def\endabstract{\par\unskip\egroup$$\egroup}
\def\widestnumber#1#2{\begingroup\let\head\null\let\subhead\empty
   \let\subsubhead\subhead
   \ifx#1\head\global\setbox\tocheadbox@\hbox{#2.\enspace}%
   \else\ifx#1\subhead\global\setbox\tocsubheadbox@\hbox{#2.\enspace}%
   \else\ifx#1\key\bgroup\let\endrefitem@\egroup
        \key#2\endrefitem@\global\refindentwd\wd\keybox@
   \else\ifx#1\no\bgroup\let\endrefitem@\egroup
        \no#2\endrefitem@\global\refindentwd\wd\nobox@
   \else\ifx#1\page\global\setbox\pagesbox@\hbox{\quad\bf#2}%
   \else\ifx#1\item\setboxz@h{#2}\global\rosteritemwd\wdz@
        \global\advance\rosteritemwd by.5\parindent
   \else\message{\string\widestnumber is not defined for this option
   (\string#1)}%
\fi\fi\fi\fi\fi\fi\endgroup}
\newif\ifmonograph@
\def\Monograph{\monograph@true \let\headmark\rightheadtext
  \let\varindent@\indent \def\headfont@{\bf}\def\proclaimheadfont@{\smc}%
  \def\demofont@{\smc}}
\let\varindent@\indent

\newbox\tocheadbox@    \newbox\tocsubheadbox@
\newbox\tocbox@
\def\toc{\toc@{Contents}}
\def\newtocdefs{%
   \def \title##1\endtitle
       {\penaltyandskip@\z@\smallskipamount
        \hangindent\wd\tocheadbox@\noindent{\bf##1}}%
   \def \chapter##1{%
        Chapter \uppercase\expandafter{\romannumeral##1.\unskip}\enspace}%
   \def \specialhead##1\endspecialhead
       {\par\hangindent\wd\tocheadbox@ \noindent##1\par}%
   \def \head##1 ##2\endhead
       {\par\hangindent\wd\tocheadbox@ \noindent
        \if\notempty{##1}\hbox to\wd\tocheadbox@{\hfil##1\enspace}\fi
        ##2\par}%
   \def \subhead##1 ##2\endsubhead
       {\par\vskip-\parskip {\normalbaselines
        \advance\leftskip\wd\tocheadbox@
        \hangindent\wd\tocsubheadbox@ \noindent
        \if\notempty{##1}\hbox to\wd\tocsubheadbox@{##1\unskip\hfil}\fi
         ##2\par}}%
   \def \subsubhead##1 ##2\endsubsubhead
       {\par\vskip-\parskip {\normalbaselines
        \advance\leftskip\wd\tocheadbox@
        \hangindent\wd\tocsubheadbox@ \noindent
        \if\notempty{##1}\hbox to\wd\tocsubheadbox@{##1\unskip\hfil}\fi
        ##2\par}}}
\def\toc@#1{\relaxnext@
   \def\page##1%
       {\unskip\penalty0\null\hfil
        \rlap{\hbox to\wd\pagesbox@{\quad\hfil##1}}\hfilneg\penalty\@M}%
 \DN@{\ifx\next\nofrills\DN@\nofrills{\nextii@}%
      \else\DN@{\nextii@{{#1}}}\fi
      \next@}%
 \DNii@##1{%
\ifmonograph@\bgroup\else\setbox\tocbox@\vbox\bgroup
   \centerline{\headfont@\ignorespaces##1\unskip}\nobreak
   \vskip\belowheadskip \fi
   \setbox\tocheadbox@\hbox{0.\enspace}%
   \setbox\tocsubheadbox@\hbox{0.0.\enspace}%
   \leftskip\indenti \rightskip\leftskip
   \setbox\pagesbox@\hbox{\bf\quad000}\advance\rightskip\wd\pagesbox@
   \newtocdefs
 }%
 \FN@\next@}
\def\endtoc{\par\egroup}
\let\pretitle\relax
\let\preauthor\relax
\let\preaffil\relax
\let\predate\relax
\let\preabstract\relax
\let\prepaper\relax
\def\dedicatory #1\enddedicatory{\def\preabstract{{\medskip
  \eightpoint\it \raggedcenter@#1\endgraf}}}
\def\thetranslator@{}
\def\translator#1\endtranslator{\def\thetranslator@{\nobreak\medskip
 \line{\eightpoint\hfil Translated by \uppercase{#1}\qquad\qquad}\nobreak}}
\outer\def\endtopmatter{\runaway@{abstract}%
 \edef\next{\the\leftheadtoks}\ifx\next\empty
  \expandafter\leftheadtext\expandafter{\the\rightheadtoks}\fi
 \ifmonograph@\else
   \ifx\thesubjclass@\empty\else \makefootnote@{}{\thesubjclass@}\fi
   \ifx\thekeywords@\empty\else \makefootnote@{}{\thekeywords@}\fi
   \ifx\thethanks@\empty\else \makefootnote@{}{\thethanks@}\fi
 \fi
  \pretitle
  \ifmonograph@ \topskip7pc \else \topskip4pc \fi
  \box\titlebox@
  \topskip10pt
  \preauthor
  \ifvoid\authorbox@\else \vskip2.5pc plus1pc \unvbox\authorbox@\fi
  \preaffil
  \ifvoid\affilbox@\else \vskip1pc plus.5pc \unvbox\affilbox@\fi
  \predate
  \ifx\thedate@\empty\else \vskip1pc plus.5pc \line{\hfil\thedate@\hfil}\fi
  \preabstract
  \ifvoid\abstractbox@\else \vskip1.5pc plus.5pc \unvbox\abstractbox@ \fi
  \ifvoid\tocbox@\else\vskip1.5pc plus.5pc \unvbox\tocbox@\fi
  \prepaper
  \vskip2pc plus1pc
}
\def\document{\let\fontlist@\relax\let\alloclist@\relax
  \tenpoint}

\newskip\aboveheadskip       \aboveheadskip1.8\bigskipamount
\newdimen\belowheadskip      \belowheadskip1.8\medskipamount

\def\headfont@{\smc}
\def\penaltyandskip@#1#2{\relax\ifdim\lastskip<#2\relax\removelastskip
      \ifnum#1=\z@\else\penalty@#1\relax\fi\vskip#2%
  \else\ifnum#1=\z@\else\penalty@#1\relax\fi\fi}
\def\nobreak{\penalty\@M
  \ifvmode\def\penalty@{\let\penalty@\penalty\count@@@}%
  \everypar{\let\penalty@\penalty\everypar{}}\fi}
\let\penalty@\penalty
\def\heading#1\endheading{\head#1\endhead}

\def\specialheadfont@{\bf}
\outer\def\specialhead{\par\penaltyandskip@{-200}\aboveheadskip
  \begingroup\interlinepenalty\@M\rightskip\z@ plus\hsize \let\\\linebreak
  \specialheadfont@\noindent\ignorespaces}
\def\endspecialhead{\par\endgroup\nobreak\vskip\belowheadskip}
\let\headmark\eat@
\newskip\subheadskip       \subheadskip\medskipamount
\def\subheadfont@{\bf}
\outer\def\subhead{\nofrills@{.\enspace}\subhead@
 \DNii@##1\endsubhead{\par\penaltyandskip@{-100}\subheadskip
  \varindent@{\usualspace@{{\subheadfont@\enspace}}%
 \subheadfont@\ignorespaces##1\unskip\subhead@}\ignorespaces}%
 \FN@\next@}
\outer\def\subsubhead{\nofrills@{.\enspace}\subsubhead@
 \DNii@##1\endsubsubhead{\par\penaltyandskip@{-50}\medskipamount
      {\usualspace@{{\it\enspace}}%
  \it\ignorespaces##1\unskip\subsubhead@}\ignorespaces}%
 \FN@\next@}
\def\proclaimheadfont@{\bf}
\outer\def\proclaim{\runaway@{proclaim}\def\envir@{proclaim}%
  \nofrills@{.\enspace}\proclaim@
 \DNii@##1{\penaltyandskip@{-100}\medskipamount\varindent@
   \usualspace@{{\proclaimheadfont@\enspace}}\proclaimheadfont@
   \ignorespaces##1\unskip\proclaim@
  \sl\ignorespaces}%
 \FN@\next@}
\outer\def\endproclaim{\let\envir@\relax\par\rm
  \penaltyandskip@{55}\medskipamount}
\def\demoheadfont@{\it}
\def\demo{\runaway@{proclaim}\nofrills@{.\enspace}\demo@
     \DNii@##1{\par\penaltyandskip@\z@\medskipamount
  {\usualspace@{{\demoheadfont@\enspace}}%
  \varindent@\demoheadfont@\ignorespaces##1\unskip\demo@}\rm
  \ignorespaces}\FN@\next@}
\def\enddemo{\par\medskip}
\def\qed{\ifhmode\unskip\nobreak\fi\quad\ifmmode\square\else$\m@th\square$\fi}
\let\remark\demo
\let\endremark\enddemo
\def\definition{\runaway@{proclaim}%
  \nofrills@{.\demoheadfont@\enspace}\definition@
        \DNii@##1{\penaltyandskip@{-100}\medskipamount
        {\usualspace@{{\demoheadfont@\enspace}}%
        \varindent@\demoheadfont@\ignorespaces##1\unskip\definition@}%
        \rm \ignorespaces}\FN@\next@}


\newdimen\rosteritemwd
\newcount\rostercount@
\newif\iffirstitem@
\let\plainitem@\item
\newtoks\everypartoks@
\def\par@{\everypartoks@\expandafter{\the\everypar}\everypar{}}
\def\roster{\edef\leftskip@{\leftskip\the\leftskip}%
 \relaxnext@
 \rostercount@\z@  
 \def\item{\FN@\rosteritem@}%
 \DN@{\ifx\next\runinitem\let\next@\nextii@\else
  \let\next@\nextiii@\fi\next@}%
 \DNii@\runinitem  
  {\unskip  
   \DN@{\ifx\next[\let\next@\nextii@\else
    \ifx\next"\let\next@\nextiii@\else\let\next@\nextiv@\fi\fi\next@}%
   \DNii@[####1]{\rostercount@####1\relax
    \enspace{\rm(\number\rostercount@)}~\ignorespaces}%
   \def\nextiii@"####1"{\enspace{\rm####1}~\ignorespaces}%
   \def\nextiv@{\enspace{\rm(1)}\rostercount@\@ne~}%
   \par@\firstitem@false  
   \FN@\next@}%
 \def\nextiii@{\par\par@  
  \penalty\@m\smallskip\vskip-\parskip
  \firstitem@true}%
 \FN@\next@}
\def\rosteritem@{\iffirstitem@\firstitem@false\else\par\vskip-\parskip\fi
 \leftskip3\parindent\noindent  
 \DNii@[##1]{\rostercount@##1\relax
  \llap{\hbox to2.5\parindent{\hss\rm(\number\rostercount@)}%
   \hskip.5\parindent}\ignorespaces}%
 \def\nextiii@"##1"{%
  \llap{\hbox to2.5\parindent{\hss\rm##1}\hskip.5\parindent}\ignorespaces}%
 \def\nextiv@{\advance\rostercount@\@ne
  \llap{\hbox to2.5\parindent{\hss\rm(\number\rostercount@)}%
   \hskip.5\parindent}}%
 \ifx\next[\let\next@\nextii@\else\ifx\next"\let\next@\nextiii@\else
  \let\next@\nextiv@\fi\fi\next@}

\newif\ifnextRunin@
\def\endroster{\relaxnext@
 \par\leftskip@  
 \penalty-50 \vskip-\parskip\smallskip  
 \DN@{\ifx\next\Runinitem\let\next@\relax
  \else\nextRunin@false\let\item\plainitem@  
   \ifx\next\par 
    \DN@\par{\everypar\expandafter{\the\everypartoks@}}%
   \else  
    \DN@{\noindent\everypar\expandafter{\the\everypartoks@}}%
  \fi\fi\next@}%
 \FN@\next@}
\newcount\rosterhangafter@
\def\Runinitem#1\roster\runinitem{\relaxnext@
 \rostercount@\z@ 
 \def\item{\FN@\rosteritem@}%
 \def\runinitem@{#1}%
 \DN@{\ifx\next[\let\next\nextii@\else\ifx\next"\let\next\nextiii@
  \else\let\next\nextiv@\fi\fi\next}%
 \DNii@[##1]{\rostercount@##1\relax
  \def\item@{{\rm(\number\rostercount@)}}\nextv@}%
 \def\nextiii@"##1"{\def\item@{{\rm##1}}\nextv@}%
 \def\nextiv@{\advance\rostercount@\@ne
  \def\item@{{\rm(\number\rostercount@)}}\nextv@}%
 \def\nextv@{\setbox\z@\vbox  
  {\ifnextRunin@\noindent\fi  
  \runinitem@\unskip\enspace\item@~\par  
  \global\rosterhangafter@\prevgraf}%
  \firstitem@false  
  \ifnextRunin@\else\par\fi  
  \hangafter\rosterhangafter@\hangindent3\parindent
  \ifnextRunin@\noindent\fi  
  \runinitem@\unskip\enspace 
  \item@~\ifnextRunin@\else\par@\fi  
  \nextRunin@true\ignorespaces}%
 \FN@\next@}
\def\footmarkform@#1{$\m@th^{#1}$}
\let\thefootnotemark\footmarkform@
\def\makefootnote@#1#2{\insert\footins
 {\interlinepenalty\interfootnotelinepenalty
 \eightpoint\splittopskip\ht\strutbox\splitmaxdepth\dp\strutbox
 \floatingpenalty\@MM\leftskip\z@\rightskip\z@\spaceskip\z@\xspaceskip\z@
 \leavevmode{#1}\footstrut\ignorespaces#2\unskip\lower\dp\strutbox
 \vbox to\dp\strutbox{}}}
\newcount\footmarkcount@
\footmarkcount@\z@
\def\footnotemark{\let\@sf\empty\relaxnext@
 \ifhmode\edef\@sf{\spacefactor\the\spacefactor}\/\fi
 \DN@{\ifx[\next\let\next@\nextii@\else
  \ifx"\next\let\next@\nextiii@\else
  \let\next@\nextiv@\fi\fi\next@}%
 \DNii@[##1]{\footmarkform@{##1}\@sf}%
 \def\nextiii@"##1"{{##1}\@sf}%
 \def\nextiv@{\iffirstchoice@\global\advance\footmarkcount@\@ne\fi
  \footmarkform@{\number\footmarkcount@}\@sf}%
 \FN@\next@}
\def\footnotetext{\relaxnext@
 \DN@{\ifx[\next\let\next@\nextii@\else
  \ifx"\next\let\next@\nextiii@\else
  \let\next@\nextiv@\fi\fi\next@}%
 \DNii@[##1]##2{\makefootnote@{\footmarkform@{##1}}{##2}}%
 \def\nextiii@"##1"##2{\makefootnote@{##1}{##2}}%
 \def\nextiv@##1{\makefootnote@{\footmarkform@{\number\footmarkcount@}}{##1}}%
 \FN@\next@}
\def\footnote{\let\@sf\empty\relaxnext@
 \ifhmode\edef\@sf{\spacefactor\the\spacefactor}\/\fi
 \DN@{\ifx[\next\let\next@\nextii@\else
  \ifx"\next\let\next@\nextiii@\else
  \let\next@\nextiv@\fi\fi\next@}%
 \DNii@[##1]##2{\footnotemark[##1]\footnotetext[##1]{##2}}%
 \def\nextiii@"##1"##2{\footnotemark"##1"\footnotetext"##1"{##2}}%
 \def\nextiv@##1{\footnotemark\footnotetext{##1}}%
 \FN@\next@}
\def\adjustfootnotemark#1{\advance\footmarkcount@#1\relax}
\def\footnoterule{\kern-3\p@
  \hrule width 5pc\kern 2.6\p@} 
\def\captionfont@{\smc}
\def\topcaption#1#2\endcaption{%
  {\dimen@\hsize \advance\dimen@-\captionwidth@
   \rm\raggedcenter@ \advance\leftskip.5\dimen@ \rightskip\leftskip
  {\captionfont@#1}%
  \if\notempty{#2}.\enspace\ignorespaces#2\fi
  \endgraf}\nobreak\bigskip}
\def\botcaption#1#2\endcaption{%
  \nobreak\bigskip
  \setboxz@h{\captionfont@#1\if\notempty{#2}.\enspace\rm#2\fi}%
  {\dimen@\hsize \advance\dimen@-\captionwidth@
   \leftskip.5\dimen@ \rightskip\leftskip
   \noindent \ifdim\wdz@>\captionwidth@ 
   \else\hfil\fi 
  {\captionfont@#1}\if\notempty{#2}.\enspace\rm#2\fi\endgraf}}
\def\@ins{\par\begingroup\def\vspace##1{\vskip##1\relax}%
  \def\captionwidth##1{\captionwidth@##1\relax}%
  \setbox\z@\vbox\bgroup} 
\def\block{\RIfMIfI@\nondmatherr@\block\fi
       \else\ifvmode\vskip\abovedisplayskip\noindent\fi
        $$\def\endblock{\par\egroup$$}\fi
  \vbox\bgroup\advance\hsize-2\indenti\noindent}
\def\endblock{\par\egroup}
\def\cite#1{{\rm[{\citefont@\m@th#1}]}}
\def\citefont@{\rm}
\def\refsfont@{\eightpoint}
\outer\def\Refs{\runaway@{proclaim}%
 \relaxnext@ \DN@{\ifx\next\nofrills\DN@\nofrills{\nextii@}\else
  \DN@{\nextii@{References}}\fi\next@}%
 \DNii@##1{\penaltyandskip@{-200}\aboveheadskip
  \line{\hfil\headfont@\ignorespaces##1\unskip\hfil}\nobreak
  \vskip\belowheadskip
  \begingroup\refsfont@\sfcode`.=\@m}%
 \FN@\next@}
\def\endRefs{\par\endgroup}
\newbox\nobox@            \newbox\keybox@           \newbox\bybox@
\newbox\paperbox@         \newbox\paperinfobox@     \newbox\jourbox@
\newbox\volbox@           \newbox\issuebox@         \newbox\yrbox@
\newbox\pagesbox@         \newbox\bookbox@          \newbox\bookinfobox@
\newbox\publbox@          \newbox\publaddrbox@      \newbox\finalinfobox@
\newbox\edsbox@           \newbox\langbox@
\newif\iffirstref@        \newif\iflastref@
\newif\ifprevjour@        \newif\ifbook@            \newif\ifprevinbook@
\newif\ifquotes@          \newif\ifbookquotes@      \newif\ifpaperquotes@
\newdimen\bysamerulewd@
\setboxz@h{\refsfont@\kern3em}
\bysamerulewd@\wdz@
\newdimen\refindentwd
\setboxz@h{\refsfont@ 00. }
\refindentwd\wdz@
\outer\def\ref{\begingroup \noindent\hangindent\refindentwd
 \firstref@true \def\nofrills{\def\refkern@{\kern3sp}}%
 \ref@}
\def\ref@{\book@false \bgroup\let\endrefitem@\egroup \ignorespaces}
\def\moreref{\endrefitem@\endref@\firstref@false\ref@}%
\def\transl{\endrefitem@\endref@\firstref@false
  \book@false
  \prepunct@
  \setboxz@h\bgroup \aftergroup\unhbox\aftergroup\z@
    \def\endrefitem@{\unskip\refkern@\egroup}\ignorespaces}%
\def\emptyifempty@{\dimen@\wd\currbox@
  \advance\dimen@-\wd\z@ \advance\dimen@-.1\p@
  \ifdim\dimen@<\z@ \setbox\currbox@\copy\voidb@x \fi}
\let\refkern@\relax
\def\endrefitem@{\unskip\refkern@\egroup
  \setboxz@h{\refkern@}\emptyifempty@}\ignorespaces
\def\refdef@#1#2#3{\edef\next@{\noexpand\endrefitem@
  \let\noexpand\currbox@\csname\expandafter\eat@\string#1box@\endcsname
    \noexpand\setbox\noexpand\currbox@\hbox\bgroup}%
  \toks@\expandafter{\next@}%
  \if\notempty{#2#3}\toks@\expandafter{\the\toks@
  \def\endrefitem@{\unskip#3\refkern@\egroup
  \setboxz@h{#2#3\refkern@}\emptyifempty@}#2}\fi
  \toks@\expandafter{\the\toks@\ignorespaces}%
  \edef#1{\the\toks@}}
\refdef@\no{}{. }
\refdef@\key{[\m@th}{] }
\refdef@\by{}{}
\def\bysame{\by\hbox to\bysamerulewd@{\hrulefill}\thinspace
   \kern0sp}
\def\manyby{\message{\string\manyby is no longer necessary; \string\by
  can be used instead, starting with version 2.0 of \styname.STY}\by}
\refdef@\paper{\ifpaperquotes@``\fi\it}{}
\refdef@\paperinfo{}{}
\def\jour{\endrefitem@\let\currbox@\jourbox@
  \setbox\currbox@\hbox\bgroup
  \def\endrefitem@{\unskip\refkern@\egroup
    \setboxz@h{\refkern@}\emptyifempty@
    \ifvoid\jourbox@\else\prevjour@true\fi}%
\ignorespaces}
\refdef@\vol{\ifbook@\else\bf\fi}{}
\refdef@\issue{no. }{}
\refdef@\yr{}{}
\refdef@\pages{}{}
\def\page{\endrefitem@\def\pp@{\def\pp@{pp.~}p.~}\let\currbox@\pagesbox@
  \setbox\currbox@\hbox\bgroup\ignorespaces}
\def\pp@{pp.~}
\def\book{\endrefitem@ \let\currbox@\bookbox@
 \setbox\currbox@\hbox\bgroup\def\endrefitem@{\unskip\refkern@\egroup
  \setboxz@h{\ifbookquotes@``\fi}\emptyifempty@
  \ifvoid\bookbox@\else\book@true\fi}%
  \ifbookquotes@``\fi\it\ignorespaces}
\def\inbook{\endrefitem@
  \let\currbox@\bookbox@\setbox\currbox@\hbox\bgroup
  \def\endrefitem@{\unskip\refkern@\egroup
  \setboxz@h{\ifbookquotes@``\fi}\emptyifempty@
  \ifvoid\bookbox@\else\book@true\previnbook@true\fi}%
  \ifbookquotes@``\fi\ignorespaces}
\refdef@\eds{(}{, eds.)}
\def\ed{\endrefitem@\let\currbox@\edsbox@
 \setbox\currbox@\hbox\bgroup
 \def\endrefitem@{\unskip, ed.)\refkern@\egroup
  \setboxz@h{(, ed.)}\emptyifempty@}(\ignorespaces}
\refdef@\bookinfo{}{}
\refdef@\publ{}{}
\refdef@\publaddr{}{}
\refdef@\finalinfo{}{}
\refdef@\lang{(}{)}

\let\refdef@\relax 
\def\ppunbox@#1{\ifvoid#1\else\prepunct@\unhbox#1\fi}
\def\nocomma@#1{\ifvoid#1\else\changepunct@3\prepunct@\unhbox#1\fi}
\def\changepunct@#1{\ifnum\lastkern<3 \unkern\kern#1sp\fi}
\def\prepunct@{\count@\lastkern\unkern
  \ifnum\lastpenalty=0
    \let\penalty@\relax
  \else
    \edef\penalty@{\penalty\the\lastpenalty\relax}%
  \fi
  \unpenalty
  \let\refspace@\ \ifcase\count@,
\or;\or.\or 
  \or\let\refspace@\relax
  \else,\fi
  \ifquotes@''\quotes@false\fi \penalty@ \refspace@
}
\def\transferpenalty@#1{\dimen@\lastkern\unkern
  \ifnum\lastpenalty=0\unpenalty\let\penalty@\relax
  \else\edef\penalty@{\penalty\the\lastpenalty\relax}\unpenalty\fi
  #1\penalty@\kern\dimen@}
\def\endref{\endrefitem@\lastref@true\endref@
  \par\endgroup \prevjour@false \previnbook@false }
\def\endref@{%
\iffirstref@
  \ifvoid\nobox@\ifvoid\keybox@\indent\fi
  \else\hbox to\refindentwd{\hss\unhbox\nobox@}\fi
  \ifvoid\keybox@
  \else\ifdim\wd\keybox@>\refindentwd
         \box\keybox@
       \else\hbox to\refindentwd{\unhbox\keybox@\hfil}\fi\fi
  \kern4sp\ppunbox@\bybox@
\fi 
  \ifvoid\paperbox@
  \else\prepunct@\unhbox\paperbox@
    \ifpaperquotes@\quotes@true\fi\fi
  \ppunbox@\paperinfobox@
  \ifvoid\jourbox@
    \ifprevjour@ \nocomma@\volbox@
      \nocomma@\issuebox@
      \ifvoid\yrbox@\else\changepunct@3\prepunct@(\unhbox\yrbox@
        \transferpenalty@)\fi
      \ppunbox@\pagesbox@
    \fi 
  \else \prepunct@\unhbox\jourbox@
    \nocomma@\volbox@
    \nocomma@\issuebox@
    \ifvoid\yrbox@\else\changepunct@3\prepunct@(\unhbox\yrbox@
      \transferpenalty@)\fi
    \ppunbox@\pagesbox@
  \fi 
  \ifbook@\prepunct@\unhbox\bookbox@ \ifbookquotes@\quotes@true\fi \fi
  \nocomma@\edsbox@
  \ppunbox@\bookinfobox@
  \ifbook@\ifvoid\volbox@\else\prepunct@ vol.~\unhbox\volbox@
  \fi\fi
  \ppunbox@\publbox@ \ppunbox@\publaddrbox@
  \ifbook@ \ppunbox@\yrbox@
    \ifvoid\pagesbox@
    \else\prepunct@\pp@\unhbox\pagesbox@\fi
  \else
    \ifprevinbook@ \ppunbox@\yrbox@
      \ifvoid\pagesbox@\else\prepunct@\pp@\unhbox\pagesbox@\fi
    \fi \fi
  \ppunbox@\finalinfobox@
  \iflastref@
    \ifvoid\langbox@.\ifquotes@''\fi
    \else\changepunct@2\prepunct@\unhbox\langbox@\fi
  \else
    \ifvoid\langbox@\changepunct@1%
    \else\changepunct@3\prepunct@\unhbox\langbox@
      \changepunct@1\fi
  \fi
}
\outer\def\enddocument{%
 \runaway@{proclaim}%
\ifmonograph@ 
\else
 \nobreak
 \thetranslator@
 \count@\z@ \loop\ifnum\count@<\addresscount@\advance\count@\@ne
 \csname address\number\count@\endcsname
 \csname email\number\count@\endcsname
 \repeat
\fi
 \vfill\supereject\end}

\def\headfont@{\headfonts}
\def\proclaimheadfont@{\bf}
\def\specialheadfont@{\bf}
\def\subheadfont@{\bf}
\def\demoheadfont@{\smc}

\newif\ifThisToToc \ThisToTocfalse
\newif\iftocloaded \tocloadedfalse

\def\C@L{\noexpand\Cal}\def\B@B{\noexpand\Bbb}\def\fR@K{\noexpand\frak}
\def\S@{\noexpand\S}\def\P@P{\noexpand\"}
\def\xpar{\\}

\def\writetoc#1{\iftocloaded\ifThisToToc\begingroup\def\totoc{}
  \def\Cal{\noexpand\C@L}\def\Bbb{\noexpand\B@B}
  \def\frak{\noexpand\fR@K}\def\goth{\frak}\def\S{\noexpand\S@}
  \def\"{\noexpand\P@P}
  \def\xpar{\par\penalty100000 }\def\idx##1{##1}\def\\{\xpar}
  \edef\next@{\write\toc{\noindent#1\leaderfill\noexpand\folio\par}}%
  \next@\endgroup\global\ThisToTocfalse\fi\fi}
\def\leaderfill{\leaders\hbox to 1em{\hss.\hss}\hfill}

\newif\ifindexloaded \indexloadedfalse
\def\idx#1{\ifindexloaded\begingroup\def\ign{}\def\it{}\def\/{}%
 \def\smc{}\def\bf{}\def\tt{}%
 \def\Cal{\noexpand\C@L}\def\Bbb{\noexpand\B@B}%
 \def\frak{\noexpand\fR@K}\def\goth{\frak}\def\S{\noexpand\S@}%
  \def\"{\noexpand\P@P}%
 {\edef\next@{\write\index{#1, \noexpand\folio}}\next@}%
 \endgroup\fi{#1}}
\def\ign#1{}

\def\totoc{\global\ThisToToctrue}

\outer\def\head#1\endhead{\par\penaltyandskip@{-200}\aboveheadskip
 {\headfont@\raggedcenter@\interlinepenalty\@M
 \ignorespaces#1\endgraf}\nobreak
 \vskip\belowheadskip
 \headmark{#1}\writetoc{#1}}

\outer\def\chaphead#1\endchaphead{\par\penaltyandskip@{-200}\aboveheadskip
 {\chapheadfonts\raggedcenter@\interlinepenalty\@M
 \ignorespaces#1\endgraf}\nobreak
 \vskip3\belowheadskip
 \headmark{#1}\writetoc{#1}}

\def\folio{{\foliofont@\ifnum\pageno<\z@ \romannumeral-\pageno
 \else\number\pageno \fi}}
\newtoks\leftheadtoks
\newtoks\rightheadtoks

\def\leftheadtext{\nofrills@{\relax}\lht@
  \DNii@##1{\leftheadtoks\expandafter{\lht@{##1}}%
    \mark{\the\leftheadtoks\noexpand\else\the\rightheadtoks}
    \ifsyntax@\setboxz@h{\def\\{\unskip\space\ignorespaces}%
        \headlinefont@##1}\fi}%
  \FN@\next@}
\def\rightheadtext{\nofrills@{\relax}\rht@
  \DNii@##1{\rightheadtoks\expandafter{\rht@{##1}}%
    \mark{\the\leftheadtoks\noexpand\else\the\rightheadtoks}%
    \ifsyntax@\setboxz@h{\def\\{\unskip\space\ignorespaces}%
        \headlinefont@##1}\fi}%
  \FN@\next@}
\def\NoRunningHeads{\global\runheads@false\global\let\headmark\eat@}

\newif\iffirstpage@     \firstpage@true
\newif\ifrunheads@      \runheads@true

\newdimen\fullhsize \fullhsize=\hsize
\newdimen\fullvsize \fullvsize=\vsize
\def\fullline{\hbox to\fullhsize}

\def\pagenumbers{\gdef\folio{\folio@}}

\let\norunningheads\NoRunningHeads
\def\userunningheads{\global\runheads@true}
\norunningheads

\headline={\def\chapter#1{\chapterno@. }%
  \def\\{\unskip\space\ignorespaces}\ifrunheads@\headlinefont@
    \ifodd\pageno\rightheadline \else\leftheadline\fi
   \else\hfil\fi\ifNoRunHeadline\global\NoRunHeadlinefalse\fi}
\let\folio@\folio
\def\foliofont@{\foliofont}
\def\foliofont{\eightrm}
\def\headlinefont@{\headlinefont}
\def\headlinefont{\eightpoint\smc}
\def\leftheadline{\rlap{\folio}\hfill
   \ifNoRunHeadline\else\iftrue\topmark\fi\fi \hfill}
\def\rightheadline{\hfill\ifNoRunHeadline
   \else \expandafter\fi
  \hfill \llap{\folio}}
\footline={{\eightpoint\bottremark}%
   \ifrunheads@\else\hfil{\let\foliofont\tenrm\folio}\fi\hfil}
\def\bottremark{}
 
\newif\ifNoRunHeadline      
\def\norunninghead{\global\NoRunHeadlinetrue}
\norunninghead

\output={\output@}
%
\newif\ifoffset\offsetfalse
\output={\output@}
\def\output@{%
 \ifoffset 
  \ifodd\count0\advance\hoffset by0.5truecm
   \else\advance\hoffset by-0.5truecm\fi\fi
 \shipout\vbox{%
  \makeheadline \pagebody \makefootline }%
 \advancepageno \ifnum\outputpenalty>-\@MM\else\dosupereject\fi}

\def\indexoutput#1{%
  \ifoffset 
   \ifodd\count0\advance\hoffset by0.5truecm
    \else\advance\hoffset by-0.5truecm\fi\fi
  \shipout\vbox{\makeheadline
  \vbox to\fullvsize{\boxmaxdepth\maxdepth%
  \ifvoid\topins\else\unvbox\topins\fi%
  #1 %
  \ifvoid\footins\else 
    \vskip\skip\footins
    \footnoterule
    \unvbox\footins\fi
  \ifr@ggedbottom \kern-\dimen@ \vfil \fi}%
  \baselineskip2pc
  \makefootline}%
 \global\advance\pageno\@ne
 \ifnum\outputpenalty>-\@MM\else\dosupereject\fi}
 
 \newbox\partialpage \newdimen\halfsize \halfsize=0.5\fullhsize
 \advance\halfsize by-0.5em

 \def\begindoublecolumns{\output={\indexoutput{\unvbox255}}%
   \begingroup \def\line{\fullline}
   \output={\global\setbox\partialpage=\vbox{\unvbox255\bigskip}}\eject
   \output={\doublecolumnout}\hsize=\halfsize \vsize=2\fullvsize}
 \def\enddoublecolumns{\output={\balancecolumns}\eject
  \endgroup \pagegoal=\fullvsize%
  \output={\output@}}
\def\doublecolumnout{\splittopskip=\topskip \splitmaxdepth=\maxdepth
  \dimen@=\fullvsize \advance\dimen@ by-\ht\partialpage
  \setbox0=\vsplit255 to \dimen@ \setbox2=\vsplit255 to \dimen@
  \indexoutput{\pagesofar} \unvbox255 \penalty\outputpenalty}
\def\pagesofar{\unvbox\partialpage
  \wd0=\hsize \wd2=\hsize \hbox to\fullhsize{\box0\hfil\box2}}
\def\balancecolumns{\setbox0=\vbox{\unvbox255} \dimen@=\ht0
  \advance\dimen@ by\topskip \advance\dimen@ by-\baselineskip
  \divide\dimen@ by2 \splittopskip=\topskip
  {\vbadness=10000 \loop \global\setbox3=\copy0
    \global\setbox1=\vsplit3 to\dimen@
    \ifdim\ht3>\dimen@ \global\advance\dimen@ by1pt \repeat}
  \setbox0=\vbox to\dimen@{\unvbox1} \setbox2=\vbox to\dimen@{\unvbox3}
  \pagesofar}

\tenpoint
\catcode`\@=\active

\def\smallheadings{\let\chapheadfonts\tenpoint\let\headfonts\tenpoint}

\tenpoint
\catcode`\@=\active

\def\LL{\leavevmode\setbox0=\hbox{L}\hbox to\wd0{\hss\char'40L}}

\def\la{\lambda}
\def\rh{\rho}

\def\De{\Delta}


\def\today{\ifcase\month\or
 January\or February\or March\or April\or May\or June\or
 July\or August\or September\or October\or November\or December\fi
 \space\number\day, \number\year}

\def\({\left(}
\def\){\right)}
\def\[{\left[}
\def\]{\right]}

\def\3{\ss}
\catcode`\@=11
\def\dddot#1{\vbox{\ialign{##\crcr
      .\hskip-.5pt.\hskip-.5pt.\crcr\noalign{\kern1.5\p@\nointerlineskip}
      $\hfil\displaystyle{#1}\hfil$\crcr}}}

\newif\iftab@\tab@false
\newif\ifvtab@\vtab@false
\def\tab{\bgroup\tab@true\vtab@false\vst@bfalse\Strich@false%
   \def\\{\global\hline@@false%
     \ifhline@\global\hline@false\global\hline@@true\fi\cr}
   \edef\l@{\the\leftskip}\ialign\bgroup\hskip\l@##\hfil&&##\hfil\cr}
\def\endtab{\cr\egroup\egroup}
\def\vtab{\vtop\bgroup\vst@bfalse\vtab@true\tab@true\Strich@false%
   \bgroup\def\\{\cr}\ialign\bgroup&##\hfil\cr}
\def\endvtab{\cr\egroup\egroup\egroup}
\def\stab{\D@cke0.5pt\null 
 \bgroup\tab@true\vtab@false\vst@bfalse\Strich@true\Let@@\vspace@
 \normalbaselines\offinterlineskip
  \openup\spreadmlines@
 \edef\l@{\the\leftskip}\ialign
 \bgroup\hskip\l@##\hfil&&##\hfil\crcr}
\def\endstab{\crcr\egroup
 \egroup}
\newif\ifvst@b\vst@bfalse
\def\vstab{\D@cke0.5pt\null
 \vtop\bgroup\tab@true\vtab@false\vst@btrue\Strich@true\bgroup\Let@@\vspace@
 \normalbaselines\offinterlineskip
  \openup\spreadmlines@\bgroup}
\def\endvstab{\crcr\egroup\egroup
 \egroup\tab@false\Strich@false}

\newdimen\htstrut@
\htstrut@8.5\p@
\newdimen\htStrut@
\htStrut@12\p@
\newdimen\dpstrut@
\dpstrut@3.5\p@
\newdimen\dpStrut@
\dpStrut@3.5\p@
\def\openup{\afterassignment\@penup\dimen@=}
\def\@penup{\advance\lineskip\dimen@
  \advance\baselineskip\dimen@
  \advance\lineskiplimit\dimen@
  \divide\dimen@ by2
  \advance\htstrut@\dimen@
  \advance\htStrut@\dimen@
  \advance\dpstrut@\dimen@
  \advance\dpStrut@\dimen@}
\def\Let@@{\relax%
    \def\\{\global\hline@@false%
     \ifhline@\global\hline@false\global\hline@@true\fi\cr}%
    \iffalse}\fi}
\def\matrix{\null\,\vcenter\bgroup
 \tab@false\vtab@false\vst@bfalse\Strich@false\Let@@\vspace@
 \normalbaselines\openup\spreadmlines@\ialign
 \bgroup\hfil$\m@th##$\hfil&&\quad\hfil$\m@th##$\hfil\crcr
 \Mathstrut@\crcr\noalign{\kern-\baselineskip}}
\def\endmatrix{\crcr\Mathstrut@\crcr\noalign{\kern-\baselineskip}\egroup
 \egroup\,}
\def\smatrix{\D@cke0.5pt\null\,
 \vcenter\bgroup\tab@false\vtab@false\vst@bfalse\Strich@true\Let@@\vspace@
 \normalbaselines\offinterlineskip
  \openup\spreadmlines@\ialign
 \bgroup\hfil$\m@th##$\hfil&&\quad\hfil$\m@th##$\hfil\crcr}
\def\endsmatrix{\crcr\egroup
 \egroup\,\Strich@false}
\newdimen\D@cke
\def\Dicke#1{\global\D@cke#1}
\newtoks\tabs@\tabs@{&}
\newif\ifStrich@\Strich@false
\newif\iff@rst

\def\Stricherr@{\iftab@\ifvtab@\errmessage{\noexpand\s not allowed
     here. Use \noexpand\vstab!}%
  \else\errmessage{\noexpand\s not allowed here. Use \noexpand\stab!}%
  \fi\else\errmessage{\noexpand\s not allowed
     here. Use \noexpand\smatrix!}\fi}
\def\format{\ifvst@b\else\crcr\fi\egroup\iffalse{\fi\ifnum`}=0 \fi\format@}
\def\format@#1\\{\def\preamble@{#1}%
 \def\Str@chfehlt##1{\ifx##1\s\Stricherr@\fi\ifx##1\\\let\Next\relax%
   \else\let\Next\Str@chfehlt\fi\Next}%
 \def\c{\hfil\noexpand\ifhline@@\hbox{\vrule height\htStrut@%
   depth\dpstrut@ width\z@}\noexpand\fi%
   \ifStrich@\hbox{\vrule height\htstrut@ depth\dpstrut@ width\z@}%
   \fi\iftab@\else$\m@th\fi\the\hashtoks@\iftab@\else$\fi\hfil}%
 \def\r{\hfil\noexpand\ifhline@@\hbox{\vrule height\htStrut@%
   depth\dpstrut@ width\z@}\noexpand\fi%
   \ifStrich@\hbox{\vrule height\htstrut@ depth\dpstrut@ width\z@}%
   \fi\iftab@\else$\m@th\fi\the\hashtoks@\iftab@\else$\fi}%
 \def\l{\noexpand\ifhline@@\hbox{\vrule height\htStrut@%
   depth\dpstrut@ width\z@}\noexpand\fi%
   \ifStrich@\hbox{\vrule height\htstrut@ depth\dpstrut@ width\z@}%
   \fi\iftab@\else$\m@th\fi\the\hashtoks@\iftab@\else$\fi\hfil}%
 \def\s{\ifStrich@\ \the\tabs@\vrule width\D@cke\the\hashtoks@%
          \fi\the\tabs@\ }%
 \def\sa{\ifStrich@\vrule width\D@cke\the\hashtoks@%
            \the\tabs@\ %
            \fi}%
 \def\se{\ifStrich@\ \the\tabs@\vrule width\D@cke\the\hashtoks@\fi}%
 \def\cd{\hfil\noexpand\ifhline@@\hbox{\vrule height\htStrut@%
   depth\dpstrut@ width\z@}\noexpand\fi%
   \ifStrich@\hbox{\vrule height\htstrut@ depth\dpstrut@ width\z@}%
   \fi$\dsize\m@th\the\hashtoks@$\hfil}%
 \def\rd{\hfil\noexpand\ifhline@@\hbox{\vrule height\htStrut@%
   depth\dpstrut@ width\z@}\noexpand\fi%
   \ifStrich@\hbox{\vrule height\htstrut@ depth\dpstrut@ width\z@}%
   \fi$\dsize\m@th\the\hashtoks@$}%
 \def\ld{\noexpand\ifhline@@\hbox{\vrule height\htStrut@%
   depth\dpstrut@ width\z@}\noexpand\fi%
   \ifStrich@\hbox{\vrule height\htstrut@ depth\dpstrut@ width\z@}%
   \fi$\dsize\m@th\the\hashtoks@$\hfil}%
 \ifStrich@\else\Str@chfehlt#1\\\fi%
 \setbox\z@\hbox{\xdef\Preamble@{\preamble@}}\ifnum`{=0 \fi\iffalse}\fi
 \ialign\bgroup\span\Preamble@\crcr}
\newif\ifhline@\hline@false
\newif\ifhline@@\hline@@false
\def\hlinefor#1{\multispan@{\strip@#1 }\leaders\hrule height\D@cke\hfill%
    \global\hline@true\ignorespaces}
\def\Item "#1"{\par\noindent\hangindent2\parindent%
  \hangafter1\setbox0\hbox{\rm#1\enspace}\ifdim\wd0>2\parindent%
  \box0\else\hbox to 2\parindent{\rm#1\hfil}\fi\ignorespaces}
\def\ITEM #1"#2"{\par\noindent\hangafter1\hangindent#1%
  \setbox0\hbox{\rm#2\enspace}\ifdim\wd0>#1%
  \box0\else\hbox to 0pt{\rm#2\hss}\hskip#1\fi\ignorespaces}
\def\item"#1"{\par\noindent\hang%
  \setbox0=\hbox{\rm#1\enspace}\ifdim\wd0>\the\parindent%
  \box0\else\hbox to \parindent{\rm#1\hfil}\enspace\fi\ignorespaces}
\let\plainitem@\item
\catcode`\@=13

\hsize13cm
\vsize19cm
\newdimen\fullhsize
\newdimen\fullvsize
\newdimen\halfsize
\fullhsize13cm
\fullvsize19cm
\halfsize=0.5\fullhsize
\advance\halfsize by-0.5em

\magnification1200

\catcode`\@=11
\font\tenln    = line10
\font\tenlnw   = linew10

\newskip\Einheit \Einheit=0.5cm
\newcount\xcoord \newcount\ycoord
\newdimen\xdim \newdimen\ydim \newdimen\PfadD@cke \newdimen\Pfadd@cke

\newcount\@tempcnta
\newcount\@tempcntb

\newdimen\@tempdima
\newdimen\@tempdimb

\newdimen\@wholewidth
\newdimen\@halfwidth

\newcount\@xarg
\newcount\@yarg
\newcount\@yyarg
\newbox\@linechar
\newbox\@tempboxa
\newdimen\@linelen
\newdimen\@clnwd
\newdimen\@clnht

\newif\if@negarg

\def\@whilenoop#1{}
\def\@whiledim#1\do #2{\ifdim #1\relax#2\@iwhiledim{#1\relax#2}\fi}
\def\@iwhiledim#1{\ifdim #1\let\@nextwhile=\@iwhiledim
        \else\let\@nextwhile=\@whilenoop\fi\@nextwhile{#1}}

\def\@whileswnoop#1\fi{}
\def\@whilesw#1\fi#2{#1#2\@iwhilesw{#1#2}\fi\fi}
\def\@iwhilesw#1\fi{#1\let\@nextwhile=\@iwhilesw
         \else\let\@nextwhile=\@whileswnoop\fi\@nextwhile{#1}\fi}

\def\thinlines{\let\@linefnt\tenln \let\@circlefnt\tencirc
  \@wholewidth\fontdimen8\tenln \@halfwidth .5\@wholewidth}
\def\thicklines{\let\@linefnt\tenlnw \let\@circlefnt\tencircw
  \@wholewidth\fontdimen8\tenlnw \@halfwidth .5\@wholewidth}
\thinlines

\PfadD@cke1pt \Pfadd@cke0.5pt
\def\PfadDicke#1{\PfadD@cke#1 \divide\PfadD@cke by2 \Pfadd@cke\PfadD@cke \multiply\PfadD@cke by2}
\long\def\LOOP#1\REPEAT{\def\BODY{#1}\ITERATE}
\def\ITERATE{\BODY \let\next\ITERATE \else\let\next\relax\fi \next}
\let\REPEAT=\fi
\def\Punkt{\hbox{\raise-2pt\hbox to0pt{\hss$\ssize\bullet$\hss}}}
\def\DuennPunkt(#1,#2){\unskip
  \raise#2 \Einheit\hbox to0pt{\hskip#1 \Einheit
          \raise-2.5pt\hbox to0pt{\hss$\bullet$\hss}\hss}}
\def\NormalPunkt(#1,#2){\unskip
  \raise#2 \Einheit\hbox to0pt{\hskip#1 \Einheit
          \raise-3pt\hbox to0pt{\hss\twelvepoint$\bullet$\hss}\hss}}
\def\DickPunkt(#1,#2){\unskip
  \raise#2 \Einheit\hbox to0pt{\hskip#1 \Einheit
          \raise-4pt\hbox to0pt{\hss\fourteenpoint$\bullet$\hss}\hss}}
\def\Kreis(#1,#2){\unskip
  \raise#2 \Einheit\hbox to0pt{\hskip#1 \Einheit
          \raise-4pt\hbox to0pt{\hss\fourteenpoint$\circ$\hss}\hss}}

\def\Line@(#1,#2)#3{\@xarg #1\relax \@yarg #2\relax
\@linelen=#3\Einheit
\ifnum\@xarg =0 \@vline
  \else \ifnum\@yarg =0 \@hline \else \@sline\fi
\fi}

\def\@sline{\ifnum\@xarg< 0 \@negargtrue \@xarg -\@xarg \@yyarg -\@yarg
  \else \@negargfalse \@yyarg \@yarg \fi
\ifnum \@yyarg >0 \@tempcnta\@yyarg \else \@tempcnta -\@yyarg \fi
\ifnum\@tempcnta>6 \@badlinearg\@tempcnta0 \fi
\ifnum\@xarg>6 \@badlinearg\@xarg 1 \fi
\setbox\@linechar\hbox{\@linefnt\@getlinechar(\@xarg,\@yyarg)}%
\ifnum \@yarg >0 \let\@upordown\raise \@clnht\z@
   \else\let\@upordown\lower \@clnht \ht\@linechar\fi
\@clnwd=\wd\@linechar
\if@negarg \hskip -\wd\@linechar \def\@tempa{\hskip -2\wd\@linechar}\else
     \let\@tempa\relax \fi
\@whiledim \@clnwd <\@linelen \do
  {\@upordown\@clnht\copy\@linechar
   \@tempa
   \advance\@clnht \ht\@linechar
   \advance\@clnwd \wd\@linechar}%
\advance\@clnht -\ht\@linechar
\advance\@clnwd -\wd\@linechar
\@tempdima\@linelen\advance\@tempdima -\@clnwd
\@tempdimb\@tempdima\advance\@tempdimb -\wd\@linechar
\if@negarg \hskip -\@tempdimb \else \hskip \@tempdimb \fi
\multiply\@tempdima \@m
\@tempcnta \@tempdima \@tempdima \wd\@linechar \divide\@tempcnta \@tempdima
\@tempdima \ht\@linechar \multiply\@tempdima \@tempcnta
\divide\@tempdima \@m
\advance\@clnht \@tempdima
\ifdim \@linelen <\wd\@linechar
   \hskip \wd\@linechar
  \else\@upordown\@clnht\copy\@linechar\fi}

\def\@hline{\ifnum \@xarg <0 \hskip -\@linelen \fi
\vrule height\Pfadd@cke width \@linelen depth\Pfadd@cke
\ifnum \@xarg <0 \hskip -\@linelen \fi}

\def\@getlinechar(#1,#2){\@tempcnta#1\relax\multiply\@tempcnta 8
\advance\@tempcnta -9 \ifnum #2>0 \advance\@tempcnta #2\relax\else
\advance\@tempcnta -#2\relax\advance\@tempcnta 64 \fi
\char\@tempcnta}

\def\Vektor(#1,#2)#3(#4,#5){\unskip\leavevmode
  \xcoord#4\relax \ycoord#5\relax
      \raise\ycoord \Einheit\hbox to0pt{\hskip\xcoord \Einheit
         \Vector@(#1,#2){#3}\hss}}

\def\Vector@(#1,#2)#3{\@xarg #1\relax \@yarg #2\relax
\@tempcnta \ifnum\@xarg<0 -\@xarg\else\@xarg\fi
\ifnum\@tempcnta<5\relax
\@linelen=#3\Einheit
\ifnum\@xarg =0 \@vvector
  \else \ifnum\@yarg =0 \@hvector \else \@svector\fi
\fi
\else\@badlinearg\fi}

\def\@hvector{\@hline\hbox to 0pt{\@linefnt
\ifnum \@xarg <0 \@getlarrow(1,0)\hss\else
    \hss\@getrarrow(1,0)\fi}}

\def\@vvector{\ifnum \@yarg <0 \@downvector \else \@upvector \fi}

\def\@svector{\@sline
\@tempcnta\@yarg \ifnum\@tempcnta <0 \@tempcnta=-\@tempcnta\fi
\ifnum\@tempcnta <5
  \hskip -\wd\@linechar
  \@upordown\@clnht \hbox{\@linefnt  \if@negarg
  \@getlarrow(\@xarg,\@yyarg) \else \@getrarrow(\@xarg,\@yyarg) \fi}%
\else\@badlinearg\fi}

\def\@upline{\hbox to \z@{\hskip -.5\Pfadd@cke \vrule width \Pfadd@cke
   height \@linelen depth \z@\hss}}

\def\@downline{\hbox to \z@{\hskip -.5\Pfadd@cke \vrule width \Pfadd@cke
   height \z@ depth \@linelen \hss}}

\def\@upvector{\@upline\setbox\@tempboxa\hbox{\@linefnt\char'66}\raise
     \@linelen \hbox to\z@{\lower \ht\@tempboxa\box\@tempboxa\hss}}

\def\@downvector{\@downline\lower \@linelen
      \hbox to \z@{\@linefnt\char'77\hss}}

\def\@getlarrow(#1,#2){\ifnum #2 =\z@ \@tempcnta='33\else
\@tempcnta=#1\relax\multiply\@tempcnta \sixt@@n \advance\@tempcnta
-9 \@tempcntb=#2\relax\multiply\@tempcntb \tw@
\ifnum \@tempcntb >0 \advance\@tempcnta \@tempcntb\relax
\else\advance\@tempcnta -\@tempcntb\advance\@tempcnta 64
\fi\fi\char\@tempcnta}

\def\@getrarrow(#1,#2){\@tempcntb=#2\relax
\ifnum\@tempcntb < 0 \@tempcntb=-\@tempcntb\relax\fi
\ifcase \@tempcntb\relax \@tempcnta='55 \or
\ifnum #1<3 \@tempcnta=#1\relax\multiply\@tempcnta
24 \advance\@tempcnta -6 \else \ifnum #1=3 \@tempcnta=49
\else\@tempcnta=58 \fi\fi\or
\ifnum #1<3 \@tempcnta=#1\relax\multiply\@tempcnta
24 \advance\@tempcnta -3 \else \@tempcnta=51\fi\or
\@tempcnta=#1\relax\multiply\@tempcnta
\sixt@@n \advance\@tempcnta -\tw@ \else
\@tempcnta=#1\relax\multiply\@tempcnta
\sixt@@n \advance\@tempcnta 7 \fi\ifnum #2<0 \advance\@tempcnta 64 \fi
\char\@tempcnta}

\def\Diagonale(#1,#2)#3{\unskip\leavevmode
  \xcoord#1\relax \ycoord#2\relax
      \raise\ycoord \Einheit\hbox to0pt{\hskip\xcoord \Einheit
         \Line@(1,1){#3}\hss}}
\def\AntiDiagonale(#1,#2)#3{\unskip\leavevmode
  \xcoord#1\relax \ycoord#2\relax 
      \raise\ycoord \Einheit\hbox to0pt{\hskip\xcoord \Einheit
         \Line@(1,-1){#3}\hss}}
\def\Pfad(#1,#2),#3\endPfad{\unskip\leavevmode
  \xcoord#1 \ycoord#2 \thicklines\ZeichnePfad#3\endPfad\thinlines}
\def\ZeichnePfad#1{\ifx#1\endPfad\let\next\relax
  \else\let\next\ZeichnePfad
    \ifnum#1=1
      \raise\ycoord \Einheit\hbox to0pt{\hskip\xcoord \Einheit
         \vrule height\Pfadd@cke width1 \Einheit depth\Pfadd@cke\hss}%
      \advance\xcoord by 1
    \else\ifnum#1=2
      \raise\ycoord \Einheit\hbox to0pt{\hskip\xcoord \Einheit
        \hbox{\hskip-\PfadD@cke\vrule height1 \Einheit width\PfadD@cke depth0pt}\hss}%
      \advance\ycoord by 1
    \else\ifnum#1=3
      \raise\ycoord \Einheit\hbox to0pt{\hskip\xcoord \Einheit
         \Line@(1,1){1}\hss}
      \advance\xcoord by 1
      \advance\ycoord by 1
    \else\ifnum#1=4
      \raise\ycoord \Einheit\hbox to0pt{\hskip\xcoord \Einheit
         \Line@(1,-1){1}\hss}
      \advance\xcoord by 1
      \advance\ycoord by -1
    \else\ifnum#1=5
      \advance\xcoord by -1
      \raise\ycoord \Einheit\hbox to0pt{\hskip\xcoord \Einheit
         \vrule height\Pfadd@cke width1 \Einheit depth\Pfadd@cke\hss}%
    \else\ifnum#1=6
      \advance\ycoord by -1
      \raise\ycoord \Einheit\hbox to0pt{\hskip\xcoord \Einheit
        \hbox{\hskip-\PfadD@cke\vrule height1 \Einheit width\PfadD@cke depth0pt}\hss}%
    \else\ifnum#1=7
      \advance\xcoord by -1
      \advance\ycoord by -1
      \raise\ycoord \Einheit\hbox to0pt{\hskip\xcoord \Einheit
         \Line@(1,1){1}\hss}
    \else\ifnum#1=8
      \advance\xcoord by -1
      \advance\ycoord by +1
      \raise\ycoord \Einheit\hbox to0pt{\hskip\xcoord \Einheit
         \Line@(1,-1){1}\hss}
    \fi\fi\fi\fi
    \fi\fi\fi\fi
  \fi\next}
\def\hSSchritt{\leavevmode\raise-.4pt\hbox to0pt{\hss.\hss}\hskip.2\Einheit
  \raise-.4pt\hbox to0pt{\hss.\hss}\hskip.2\Einheit
  \raise-.4pt\hbox to0pt{\hss.\hss}\hskip.2\Einheit
  \raise-.4pt\hbox to0pt{\hss.\hss}\hskip.2\Einheit
  \raise-.4pt\hbox to0pt{\hss.\hss}\hskip.2\Einheit}
\def\vSSchritt{\vbox{\baselineskip.2\Einheit\lineskiplimit0pt
\hbox{.}\hbox{.}\hbox{.}\hbox{.}\hbox{.}}}
\def\DSSchritt{\leavevmode\raise-.4pt\hbox to0pt{%
  \hbox to0pt{\hss.\hss}\hskip.2\Einheit
  \raise.2\Einheit\hbox to0pt{\hss.\hss}\hskip.2\Einheit
  \raise.4\Einheit\hbox to0pt{\hss.\hss}\hskip.2\Einheit
  \raise.6\Einheit\hbox to0pt{\hss.\hss}\hskip.2\Einheit
  \raise.8\Einheit\hbox to0pt{\hss.\hss}\hss}}
\def\dSSchritt{\leavevmode\raise-.4pt\hbox to0pt{%
  \hbox to0pt{\hss.\hss}\hskip.2\Einheit
  \raise-.2\Einheit\hbox to0pt{\hss.\hss}\hskip.2\Einheit
  \raise-.4\Einheit\hbox to0pt{\hss.\hss}\hskip.2\Einheit
  \raise-.6\Einheit\hbox to0pt{\hss.\hss}\hskip.2\Einheit
  \raise-.8\Einheit\hbox to0pt{\hss.\hss}\hss}}
\def\SPfad(#1,#2),#3\endSPfad{\unskip\leavevmode
  \xcoord#1 \ycoord#2 \ZeichneSPfad#3\endSPfad}
\def\ZeichneSPfad#1{\ifx#1\endSPfad\let\next\relax
  \else\let\next\ZeichneSPfad
    \ifnum#1=1
      \raise\ycoord \Einheit\hbox to0pt{\hskip\xcoord \Einheit
         \hSSchritt\hss}%
      \advance\xcoord by 1
    \else\ifnum#1=2
      \raise\ycoord \Einheit\hbox to0pt{\hskip\xcoord \Einheit
        \hbox{\hskip-2pt \vSSchritt}\hss}%
      \advance\ycoord by 1
    \else\ifnum#1=3
      \raise\ycoord \Einheit\hbox to0pt{\hskip\xcoord \Einheit
         \DSSchritt\hss}
      \advance\xcoord by 1
      \advance\ycoord by 1
    \else\ifnum#1=4
      \raise\ycoord \Einheit\hbox to0pt{\hskip\xcoord \Einheit
         \dSSchritt\hss}
      \advance\xcoord by 1
      \advance\ycoord by -1
    \else\ifnum#1=5
      \advance\xcoord by -1
      \raise\ycoord \Einheit\hbox to0pt{\hskip\xcoord \Einheit
         \hSSchritt\hss}%
    \else\ifnum#1=6
      \advance\ycoord by -1
      \raise\ycoord \Einheit\hbox to0pt{\hskip\xcoord \Einheit
        \hbox{\hskip-2pt \vSSchritt}\hss}%
    \else\ifnum#1=7
      \advance\xcoord by -1
      \advance\ycoord by -1
      \raise\ycoord \Einheit\hbox to0pt{\hskip\xcoord \Einheit
         \DSSchritt\hss}
    \else\ifnum#1=8
      \advance\xcoord by -1
      \advance\ycoord by 1
      \raise\ycoord \Einheit\hbox to0pt{\hskip\xcoord \Einheit
         \dSSchritt\hss}
    \fi\fi\fi\fi
    \fi\fi\fi\fi
  \fi\next}
\def\Koordinatenachsen(#1,#2){\unskip
 \hbox to0pt{\hskip-.5pt\vrule height#2 \Einheit width.5pt depth1 \Einheit}%
 \hbox to0pt{\hskip-1 \Einheit \xcoord#1 \advance\xcoord by1
    \vrule height0.25pt width\xcoord \Einheit depth0.25pt\hss}}
\def\Koordinatenachsen(#1,#2)(#3,#4){\unskip
 \hbox to0pt{\hskip-.5pt \ycoord-#4 \advance\ycoord by1
    \vrule height#2 \Einheit width.5pt depth\ycoord \Einheit}%
 \hbox to0pt{\hskip-1 \Einheit \hskip#3\Einheit 
    \xcoord#1 \advance\xcoord by1 \advance\xcoord by-#3 
    \vrule height0.25pt width\xcoord \Einheit depth0.25pt\hss}}
\def\Gitter(#1,#2){\unskip \xcoord0 \ycoord0 \leavevmode
  \LOOP\ifnum\ycoord<#2
    \loop\ifnum\xcoord<#1
      \raise\ycoord \Einheit\hbox to0pt{\hskip\xcoord \Einheit\Punkt\hss}%
      \advance\xcoord by1
    \repeat
    \xcoord0
    \advance\ycoord by1
  \REPEAT}
\def\Gitter(#1,#2)(#3,#4){\unskip \xcoord#3 \ycoord#4 \leavevmode
  \LOOP\ifnum\ycoord<#2
    \loop\ifnum\xcoord<#1
      \raise\ycoord \Einheit\hbox to0pt{\hskip\xcoord \Einheit\Punkt\hss}%
      \advance\xcoord by1
    \repeat
    \xcoord#3
    \advance\ycoord by1
  \REPEAT}
\def\Label#1#2(#3,#4){\unskip \xdim#3 \Einheit \ydim#4 \Einheit
  \def\lo{\advance\xdim by-.5 \Einheit \advance\ydim by.5 \Einheit}%
  \def\llo{\advance\xdim by-.25cm \advance\ydim by.5 \Einheit}%
  \def\loo{\advance\xdim by-.5 \Einheit \advance\ydim by.25cm}%
  \def\o{\advance\ydim by.25cm}%
  \def\ro{\advance\xdim by.5 \Einheit \advance\ydim by.5 \Einheit}%
  \def\rro{\advance\xdim by.25cm \advance\ydim by.5 \Einheit}%
  \def\roo{\advance\xdim by.5 \Einheit \advance\ydim by.25cm}%
  \def\l{\advance\xdim by-.30cm}%
  \def\r{\advance\xdim by.30cm}%
  \def\lu{\advance\xdim by-.5 \Einheit \advance\ydim by-.6 \Einheit}%
  \def\llu{\advance\xdim by-.25cm \advance\ydim by-.6 \Einheit}%
  \def\luu{\advance\xdim by-.5 \Einheit \advance\ydim by-.30cm}%
  \def\u{\advance\ydim by-.30cm}%
  \def\ru{\advance\xdim by.5 \Einheit \advance\ydim by-.6 \Einheit}%
  \def\rru{\advance\xdim by.25cm \advance\ydim by-.6 \Einheit}%
  \def\ruu{\advance\xdim by.5 \Einheit \advance\ydim by-.30cm}%
  #1\raise\ydim\hbox to0pt{\hskip\xdim
     \vbox to0pt{\vss\hbox to0pt{\hss$#2$\hss}\vss}\hss}%
}
\catcode`\@=13

\TagsOnRight

\def\BeriAA{1}
\def\BrFoAA{2}
\def\ChDDAB{3}
\def\FomiAZ{4}
\def\FomiAB{5}
\def\FomiAF{6}
\def\GreCAA{7}
\def\HaLeAA{8}
\def\KratCE{9}
\def\MaRoAA{10}
\def\RobyAA{11}
\def\RobyAD{12}
\def\SagaAQ{13}
\def\StanBI{14}

\def\AA{1.1}
\def\AB{1.2}
\def\AC{1.3}
\def\AD{1.4}
\def\AE{1.5}
\def\AF{1.6}

\def\BAa{3.1a}
\def\BAb{3.1b}
\def\BAc{3.1c}
\def\BB{3.2}

\def\CACA{4.1}

\def\CAa{4.2a}
\def\CAb{4.2b}
\def\CAc{4.2c}

\def\DAa{5.1a}
\def\DAb{5.1b}
\def\DAc{5.1c}
\def\DB{5.2}

\def\EAa{6.1a}
\def\EAb{6.1b}
\def\EAc{6.1c}

\def\FA{7.1}
\def\FAa{7.1a}
\def\FAb{7.1b}
\def\FAc{7.1c}
\def\FE{7.2}
\def\FE{7.3}
\def\FE{7.4}

\def\TA{1}
\def\TB{2}
\def\TC{3}

%

\def\FA{1}
\def\FB{2}
\def\FC{3}
\def\FD{4}
\def\FE{5}
\def\FEa{6}
\def\FF{7}
\def\FFa{8}
\def\FG{9}
\def\FH{10}
\def\FI{11}
\def\FJ{12}
\def\FK{13}

\topmatter 
\title Identities for vacillating tableaux via growth diagrams
\endtitle 
\author C.~Krattenthaler
\endauthor 
\affil 
Fakult\"at f\"ur Mathematik, Universit\"at Wien,\\
Oskar-Morgenstern-Platz~1, A-1090 Vienna, Austria.\\
WWW: \tt http://www.mat.univie.ac.at/\~{}kratt
\endaffil
\address Fakult\"at f\"ur Mathematik, Universit\"at Wien,
Oskar-Morgenstern-Platz~1, A-1090 Vienna, Austria.\newline
http://www.mat.univie.ac.at/\~{}kratt
\endaddress

\subjclass Primary 05A15;
 Secondary 05A17 05A19 05E10
\endsubjclass
\keywords Growth diagrams, Robinson--Schensted correspondence, standard Young ta\-bleaux,
oscillating tableaux, vacillating tableaux,
set partitions, Stirling numbers of the second kind, integer partitions
\endkeywords
\abstract 
We give bijective proofs using Fomin's growth diagrams
for identities involving numbers of
vacillating tableaux that arose in the representation theory of
partition algebras or are inspired by such identities.
\endabstract
\endtopmatter
\document

\subhead 1. Introduction\endsubhead
Recently, there has been a resurgence of interest in algorithmic
bijections involving vacillating tableaux, due to Halverson and
Lewandowski~\cite{\HaLeAA}, which produce combinatorial proofs of
identities that arise in the representation theory of partition
algebras. Here, given a positive integer~$K$,
a {\it vacillating tableau} of length~$K$ from the (integer) partition~$\la$ to the
partition~$\mu$ (see Section~2 for all definitions) is a sequence of partitions
$\la=\la^0\supset
\la^1\subset\la^2\supset\cdots\subset \la^{2K}=\mu$, for some positive integer~$K$,
where $\la^i$ and $\la^{i+1}$ differ by exactly one cell for all~$i$.
We write $m^\la_\mu(K)$ for the number of vacillating tableaux
from~$\la$ to~$\mu$ and, as usual,
$f^\la$ for the number of {\it standard Young tableaux} of shape~$\la$.
Using this notation, the identities considered in~\cite{\HaLeAA} are
$$
n^k=\sum_{\la\vdash n}f^\la m_{(n)}^\la(k),\quad \text{for }n\ge1,
\tag\AA$$
and
$$
B_{2k}=\sum_{\la\vdash n}\left(m_{(n)}^\la(k)\right)^2
=m_{(n)}^{(n)}(2k),
\quad \text{for }n\ge 2k,
\tag\AB$$
where $B_m$ is the $m$-th {\it Bell number} (the number of all set
partitions of~$\{1,2,\dots,m\}$).
(The expression in the centre results trivially from the right-hand side by
cutting the vacillating tableau from~$(n)$ to~$(n)$ halfway into two
vacillating tableaux.)
As a matter of fact, Martin and Rollet~\cite{\MaRoAA}
had proven the more general identity
$$
\sum_{l=1}^n S(k,l)=m_{(n)}^{(n)}(k),\quad \text{for }n\ge1,
\tag\AC$$
via a (rather impenetrable) recursively built bijection.
(See Appendix~A for a completely worked out special case of the identity.)
Here, the number $S(k,l)$ is a {\it Stirling number of the second kind},
giving the number of all partitions of a $k$-element set into $l$~blocks.
Clearly, the identity~(\AB) is the special case of~(\AC) where we
replace $k$ by~$2k$.

It is the article \cite{\BeriAA}, in which Berikkyzy et al\. study
fine properties of Halverson and Lewandowski's algorithm, which drew my
attention to~\cite{\HaLeAA}. The mentioned algorithm is an extremely
elegant deletion-insertion algorithm that proves~(\AA) bijectively.
Now, whenever I see (Robinson--Schensted) ``insertion", my immediate
reaction is: there must also be a presentation of the algorithm in
terms of Fomin's {\it growth diagrams} \cite{\FomiAZ, \FomiAB, \FomiAF}
(see \cite{\RobyAA, \RobyAD},
\cite{\SagaAQ, Sec.~5.2} and \cite{\StanBI, Sec.~7.13} for non-technical
expositions)! On my search in the literature whether this has already
been explained, I discovered that I had commented on that issue
in~\cite{\KratCE} (which I had completely forgotten\footnote{I believe
that the article \cite{\HaLeAA} had been brought to my attention by a referee.}):

\medskip
{\leftskip1cm\rightskip1cm\noindent\it
\dots\ we point out that \dots\ Robinson--Schensted
like algorithms between set partitions and oscillating sequences of
(integer) partitions have also been constructed by
Halverson and Lewandowski, with the completely
different motivation of explaining combinatorial identities 
arising from the representation theory of the partition algebra.
Halverson and Lewandowski provide both the insertion/deletion and the
growth diagram presentation of the algorithms. However, in their
considerations, Greene's theorem does not play any role.
\par\medskip
}

In retrospect: mostly correct, but not entirely. There are two
deletion-insertion algorithms in~\cite{\HaLeAA}: one that proves~(\AA),
and one that proves~(\AB). Halverson and Lewandowski do provide a
growth diagram description of the latter algorithm, but not of the
former. Second, the growth diagram description of the latter algorithm
is not quite the ``right one" in the sense that it does not extend to a bijective
proof of the more general identity~(\AC). Furthermore, by missing
these points, I also missed that, by looking into the growth diagram
perspective of the above identities, there are more identities to be
discovered in this context.

The purpose of the present article is to make up for these oversights.
After fixing notations and recalling the concept of growth diagrams in
the next section (largely copied from~\cite{\KratCE}) in order to be
self-contained, we present a
growth diagram bijection proving~(\AA) 
in Section~3. At this point, a disclaimer is in order: as opposed
to \cite{\KratCE}, which provides growth diagram versions of the
bijections in~\cite{\ChDDAB}, here I do not claim that our bijection
proving~(\AA) is a growth diagram version of Halverson and
Lewandowski's bijection. Rather, it seems that the two are not related
in any simple way. See the end of Section~3 for more specific comments.
On the other hand, they share a limiting property which is
investigated by Berikkyzy et al\. in~\cite{\BeriAA} for the algorithm
of Halverson and Lewandowski. See Theorem~\TC\ in the same section.

With this perspective in mind, one realises that~(\AA)
can be generalised to\footnote{This, corrected, version of the
identity is due to Catherine Yan; see the Acknowledgement.}
$$
f^\mu n^k=\sum_{\la\vdash n}f^\la m_{\mu}^\la(k),\quad \text{for }n\ge1,
\tag\AD$$
where $\mu$ is {\it any} fixed partition of~$n$. The corresponding bijective
proof using growth diagrams is presented in Section~4. (The
identity~(\AD) could also be proved by an appropriate adaptation of the
deletion-insertion algorithm of Halverson and Lewandowski.)
We move on in Section~5 to give a growth diagram bijection
proving~(\AC). That proof provides the inspiration for variations of that
identity. In Section~6, we present a growth diagram
bijection proving
$$
S(k,n)+S(k,n-1)=m_{(n)}^{(1^n)}(k),\quad \text{for }n\ge1,
\tag\AE$$
where $1^n$ is short for a sequence of $n$~1's.
(See Appendix~B for a completely worked out special case of the identity.)
Moreover, in Section~7, we argue that we have
$$
\sum_{l=1}^{n-2}S(k,l)+\sum_{l=1}^{n-1}l^2S(k,l)+(n-1)^2S(k,n)
=m_{(n-1,1)}^{(n-1,1)}(k),\quad \text{for }n\ge2.
\tag\AF$$
(See Appendix~C for a completely worked out special case of the identity.)
Obviously, more identities of this type could be derived in the same spirit.

We emphasise that, once one has understood the growth-diagram
machinery, all these proofs are {\it one-picture proofs}: the growth diagram
picture makes the corresponding identities immediately ``obvious".



\medskip
{\smc Acknowledgement.} I am indebted to Catherine Yan for pointing out
an error in Equation~(\AD) of the first version of this manuscript
(and, consequently, in the corresponding ``proof''), at the same time
providing the necessary correction.

\subhead 2. Definitions and notation \endsubhead
We start by fixing the standard partition notation (cf.\
e.g.\ \cite{\StanBI, Sec.~7.2}).
A {\it partition} is a weakly decreasing sequence
$\la=(\la_1,\la_2,\dots,\la_\ell)$ of positive integers.
This also includes the {\it empty partition} $()$, denoted by
$\emptyset$.
To each partition~$\la$, one
associates its {\it Ferrers diagram} (also called {\it Ferrers shape}), 
which is the left-justified
arrangement of squares with $\la_i$ squares in the $i$-th row,
$i=1,2,\dots$.
If $n=\la_1+\la_2+\dots+\la_\ell$, then we say that $\la$ is a
partition of (size)
~$n$, and we write $\la\vdash n$.
We define a {\it partial order} $\subseteq$ on partitions 
by containment of their Ferrers diagrams.
The {\it union} $\mu\cup\nu$ of two partitions $\mu$ and $\nu$
is the partition which arises by forming the union of the Ferrers
diagrams of $\mu$ and $\nu$. Thus, if
$\mu=(\mu_1,\mu_2,\dots)$ and
$\nu=(\nu_1,\nu_2,\dots)$, then $\mu\cup\nu$ is the partition 
$\la=(\la_1,\la_2,\dots)$, where $\la_i=\max\{\mu_i,\nu_i\}$ for 
$i=1,2,\dots$. The {\it intersection} 
$\mu\cap\nu$ of two partitions $\mu$ and $\nu$
is the partition which arises by forming the intersection of the Ferrers
diagrams of $\mu$ and $\nu$. Thus, if
$\mu=(\mu_1,\mu_2,\dots)$ and
$\nu=(\nu_1,\nu_2,\dots)$, then $\mu\cap\nu$ is the partition 
$\rh=(\rh_1,\rh_2,\dots)$, where $\rh_i=\min\{\mu_i,\nu_i\}$ for 
$i=1,2,\dots$.
The {\it conjugate of\/} a partition $\lambda$ is the partition 
$\la'=(\lambda^\prime_1, \dots, \lambda^\prime_{\lambda_1})$ where
$\lambda_j'$ is the length of the $j$-th column in the Ferrers
diagram of $\lambda$.

\medskip
Given a partition $\la=(\la_1,\la_2,\dots,\la_\ell)$, a {\it standard
{\rrm(}Young{\rm)} tableau of shape\/} $\la$ is a left-justified arrangement
of the integers $1,2,\dots,\la_1+\la_2+\dots+\la_\ell$
with $\la_i$ entries in row~$i$, $i=1,2,\dots$,
such that the entries along rows and columns are increasing.
By considering the sequence of partitions (Ferrers shapes)
$(\la^{i})_{i\ge0}$, where $\la^i$ is the shape formed by the
entries of $T$ which are at most~$i$,
$i=0,1,2,\dots$, one sees that standard tableaux of shape $\la$ are in
bijection with sequences 
$\emptyset=\la^0\subset\la^1\subset\dots\subset\la^n=\la$,
where $\la^{i-1}$ and $\la^i$ differ by exactly one square for all~$i$.


\medskip
{\it Growth diagrams} are certain labellings of arrangements of cells.
The arrangements of cells which we need here are arrangements which
are left-justified (that is, they have a straight vertical left
boundary),
bottom-justified (that is, they have a straight horizontal bottom
boundary), and rows and columns in the arrangement are ``without"
holes, that is, if we move along the top-right boundary of the
arrangement, we always move either to the right or to the bottom.
Figure~\FA.a shows an example of such a cell arrangement.

\midinsert
$$
\Einheit.4cm
\Pfad(0,18),666666666666666666111111111111111\endPfad
\Pfad(0,18),111666666666666666666\endPfad
\Pfad(0,15),111111666666666666666\endPfad
\Pfad(0,12),111111666111666666666\endPfad
\Pfad(0,9),111111111666666111666\endPfad
\Pfad(0,6),111111111666111111666\endPfad
\Pfad(0,3),111111111111111666\endPfad
\hbox{\hskip7.5cm}
\Pfad(0,18),666666666666666666111111111111111\endPfad
\Pfad(0,18),111666666666666666666\endPfad
\Pfad(0,15),111111666666666666666\endPfad
\Pfad(0,12),111111666111666666666\endPfad
\Pfad(0,9),111111111666666111666\endPfad
\Pfad(0,6),111111111666111111666\endPfad
\Pfad(0,3),111111111111111666\endPfad
\Label\ro{\text {\seventeenpoint $X$}}(4,4)
\Label\ro{\text {\seventeenpoint $X$}}(1,10)
\Label\ro{\text {\seventeenpoint $X$}}(13,1)
\hskip6cm
$$
\vskip10pt
\centerline{\eightpoint a. A cell arrangement
\hskip4cm
b. A filling of the cell arrangement}
\vskip6pt
\centerline{\smc Figure \FA}
\endinsert

We fill some cells of such an arrangement $C$ with crosses~$X$
such that every row and every column contains at most one~$X$.
See Figure~\FA.b for an example.

Next, the corners of the cells are labelled by partitions
such that the following two conditions are satisfied:

\roster
\item"(C1)" A partition is either equal to its right neighbour
or smaller by exactly one square, the same being true for a
partition and its top neighbour.
\item"(C2)" A partition and its right neighbour are equal if and only
if in the column of cells of $C$ below them there appears no $X$ and if
their bottom neighbours are also equal to each other.
Similarly, a partition and its top neighbour are equal if and only if
in the row of cells of $C$ to the left of them there appears no $X$ and if
their left neighbours are also equal to each other.
\endroster

See Figure~\FB\ for an example. (More examples can be found in
Figures~\FD--\FK.)
There, we use a short notation for partitions. For example,
$11$ is short for $(1,1)$.
Indeed, the
filling represented in Figure~\FB\ is the same as the one in
Figure~\FA.b.

Diagrams which obey the conditions (C1) and (C2) are called {\it
growth diagrams}.

\vskip10pt
\midinsert
$$
\Einheit.4cm
\Pfad(0,18),666666666666666666111111111111111\endPfad
\Pfad(0,18),111666666666666666666\endPfad
\Pfad(0,15),111111666666666666666\endPfad
\Pfad(0,12),111111666111666666666\endPfad
\Pfad(0,9),111111111666666111666\endPfad
\Pfad(0,6),111111111666111111666\endPfad
\Pfad(0,3),111111111111111666\endPfad
\Label\ro{\emptyset}(0,0)
\Label\ro{\emptyset}(0,3)
\Label\ro{\emptyset}(0,6)
\Label\ro{\emptyset}(0,9)
\Label\ro{\emptyset}(0,12)
\Label\ro{\emptyset}(0,15)
\Label\ro{\emptyset}(0,18)
\Label\ro{1}(3,18)
\Label\ro{1}(3,15)
\Label\ro{1}(3,12)
\Label\ro{\emptyset}(3,9)
\Label\ro{\emptyset}(3,6)
\Label\ro{\emptyset}(3,3)
\Label\ro{\emptyset}(3,0)
\Label\ro{\hphantom{1}11}(6,15)
\Label\ro{\hphantom{1}11}(6,12)
\Label\ro{1}(6,9)
\Label\ro{1}(6,6)
\Label\ro{\emptyset}(6,3)
\Label\ro{\emptyset}(6,0)
\Label\ro{1}(9,9)
\Label\ro{1}(9,6)
\Label\ro{\emptyset}(9,3)
\Label\ro{\emptyset}(9,0)
\Label\ro{\emptyset}(12,3)
\Label\ro{\emptyset}(12,0)
\Label\ro{1}(15,3)
\Label\ro{\emptyset}(15,0)
\Label\ro{\text {\seventeenpoint X}}(4,4)
\Label\ro{\text {\seventeenpoint X}}(1,10)
\Label\ro{\text {\seventeenpoint X}}(13,1)
\hskip6cm
$$
\vskip10pt
\centerline{\eightpoint A growth diagram}
\vskip6pt
\centerline{\smc Figure \FB}
\endinsert
\vskip10pt

We are interested in growth diagrams which
obey the following ({\it forward}) {\it local rules} 
(see Figure~\FC).

\vskip10pt
\midinsert
$$
\Pfad(0,0),1111222255556666\endPfad
\Label\lu{\rh}(0,0)
\Label\ru{\mu}(4,0)
\Label\lo{\nu}(0,4)
\Label\ro{\la}(4,4)
\hbox{\hskip6cm}
\Pfad(0,0),1111222255556666\endPfad
\Label\lu{\rh}(0,0)
\Label\ru{\mu}(4,0)
\Label\lo{\nu}(0,4)
\Label\ro{\la}(4,4)
\thicklines
\Pfad(1,1),33\endPfad
\raise.5pt\hbox to 0pt{\hbox{\hskip-.5pt}\Pfad(1,1),3\endPfad\hss}%
\raise.5pt\hbox to 0pt{\hbox{\hskip-.5pt}\Pfad(2,2),3\endPfad\hss}%
\raise1pt\hbox to 0pt{\hbox{\hskip-1pt}\Pfad(1,1),3\endPfad\hss}%
\raise1pt\hbox to 0pt{\hbox{\hskip-1pt}\Pfad(2,2),3\endPfad\hss}%
\raise-.5pt\hbox to 0pt{\hbox{\hskip.5pt}\Pfad(1,1),3\endPfad\hss}%
\raise-.5pt\hbox to 0pt{\hbox{\hskip.5pt}\Pfad(2,2),3\endPfad\hss}%
\raise-1pt\hbox to 0pt{\hbox{\hskip1pt}\Pfad(1,1),3\endPfad\hss}%
\raise-1pt\hbox to 0pt{\hbox{\hskip1pt}\Pfad(2,2),3\endPfad\hss}%
\Pfad(1,3),44\endPfad
\raise.5pt\hbox to 0pt{\hbox{\hskip.5pt}\Pfad(1,3),4\endPfad\hss}%
\raise.5pt\hbox to 0pt{\hbox{\hskip.5pt}\Pfad(2,2),4\endPfad\hss}%
\raise1pt\hbox to 0pt{\hbox{\hskip1pt}\Pfad(1,3),4\endPfad\hss}%
\raise1pt\hbox to 0pt{\hbox{\hskip1pt}\Pfad(2,2),4\endPfad\hss}%
\raise-.5pt\hbox to 0pt{\hbox{\hskip-.5pt}\Pfad(1,3),4\endPfad\hss}%
\raise-.5pt\hbox to 0pt{\hbox{\hskip-.5pt}\Pfad(2,2),4\endPfad\hss}%
\raise-1pt\hbox to 0pt{\hbox{\hskip-1pt}\Pfad(1,3),4\endPfad\hss}%
\raise-1pt\hbox to 0pt{\hbox{\hskip-1pt}\Pfad(2,2),4\endPfad\hss}%
\hskip2cm
$$
\vskip10pt
\centerline{\eightpoint a. A cell without cross\hskip3cm
b. A cell with cross}
\vskip6pt
\centerline{\smc Figure \FC}
\endinsert
\vskip10pt

\roster
\item"(F1)" If $\rh=\mu=\nu$, and if there is no cross in the
cell, then $\la=\rh$.
\item"(F2)" If $\rh=\mu\ne\nu$, then $\la=\nu$.
\item"(F3)" If $\rh=\nu\ne\mu$, then $\la=\mu$.
\item"(F4)" If $\rh,\mu,\nu$ are pairwise different, 
then $\la=\mu\cup\nu$.
\item"(F5)" If $\rh\ne \mu=\nu$, 
then $\la$ is formed by adding a square to the $(k+1)$-st row of
$\mu=\nu$, given that $\mu=\nu$ and $\rh$ differ in the $k$-th row.
\item"(F6)" If $\rh=\mu=\nu$, and if there is a cross in the
cell, then $\la$ is formed by adding a square to the first row of
$\rh=\mu=\nu$.
\endroster

Thus, if we label all the corners along the left and the bottom boundary
by empty partitions (which we shall always do in this paper), 
these rules allow one to determine all other labels
of corners uniquely.

It is not difficult to see that
the rules (F5) and (F6) are designed so that
one can also work one's way in the other direction, that is, given
$\la,\mu,\nu$, one can reconstruct $\rh$ {\it and\/} the filling of the cell.
The corresponding ({\it backward}) {\it local rules} are:

\roster
\item"(B1)" If $\la=\mu=\nu$, then $\rh=\la$.
\item"(B2)" If $\la=\mu\ne\nu$, then $\rh=\nu$.
\item"(B3)" If $\la=\nu\ne\mu$, then $\rh=\mu$.
\item"(B4)" If $\la,\mu,\nu$ are pairwise different, 
then $\rh=\mu\cap\nu$.
\item"(B5)" If $\la\ne \mu=\nu$, 
then $\rh$ is formed by deleting a square from the $(k-1)$-st row of
$\mu=\nu$, given that $\mu=\nu$ and $\la$ differ in the $k$-th row,
$k\ge2$.
\item"(B6)" If $\la\ne \mu=\nu$, and if $\la$ and $\mu=\nu$ differ in
the first row, then $\rh=\mu=\nu$.
\item""\hskip-36pt\vbox{\hsize13.76cm\noindent
In case (B6) the cell is filled with a cross.
In all other cases the cell is left empty.}
\endroster

Thus, given a labelling of the corners along the top-right boundary of a
cell arrangement, one can algorithmically reconstruct the labels of
the other corners of the cells {\it and\/} of the filling 
by working one's way to the left and to the bottom.
These observations lead to the following theorem.

\proclaim{Theorem \TA} 
Let $C$ be an arrangement of cells.
The fillings of $C$ with the property that every row and
every column contains at most one $X$ are in bijection with labellings
$(\emptyset=\la^0,\la^1,\dots,\la^{k}=\emptyset)$ of the corners of
cells appearing along the top-right boundary of~$C$, where $\la^{i-1}$
and $\la^i$ differ by at most one square, and
$\la^{i-1}\subseteq\la^i$ if $\la^{i-1}$ and $\la^i$ appear along
a horizontal edge, whereas
$\la^{i-1}\supseteq\la^i$ if $\la^{i-1}$ and $\la^i$ appear along
a vertical edge.
Moreover, $\la^{i-1}\subsetneqq\la^i$ if and only if there is an~$X$
in the column of cells of $C$ below the corners labelled by $\la^{i-1}$ and
$\la^i$, and $\la^{i-1}\supsetneqq\la^i$ if and only if there is an~$X$
in the row of cells of $C$ to the left of the corners labelled by
$\la^{i-1}$ and $\la^i$.
\endproclaim

\medskip
In addition to its local description,
the bijection of the above theorem has also a {\it global\/}
description. The latter is a consequence of a theorem of Greene
\cite{\GreCAA} (see also \cite{\BrFoAA, Theorems~2.1 and 3.2}). 
In order to formulate the
result, we need the following definitions: a {\it NE-chain} of a
filling is a sequence of $X$'s in the filling such that any $X$
in the sequence is above and to the right of the preceding $X$ in 
the sequence. Similarly, a {\it SE-chain} of a
filling is a set of $X$'s in the filling such that any $X$
in the sequence is below and to the right of the preceding $X$ in the sequence.

\proclaim{Theorem \TB} 
Given a growth diagram on a cell arrangement 
with empty partitions labelling all the
corners along the left boundary and the bottom boundary of the cell arrangement,
the partition $\la=(\la_1,\la_2,\dots,\la_\ell)$ labelling corner $c$ 
satisfies the following two properties:
\roster
\item"(G1)"For any $k$, the maximal cardinality of the union of $k$
NE-chains situated in the rectangular region to the left and
below of $c$ is equal to $\la_1+\la_2+\dots+\la_k$.
\item"(G2)"For any $k$, the maximal cardinality of the union of $k$
SE-chains situated in the rectangular region to the left and
below of $c$ is equal to $\la'_1+\la'_2+\dots+\la'_k$, where $\la'$
denotes the partition conjugate to $\la$.
\endroster
In particular, $\la_1$ is the length of the longest NE-chain
in the rectangular region to the left and below of $c$, and $\la'_1$
is the length of the longest SE-chain in the same
rectangular region.
\endproclaim

\subhead 3. Bijective proof of (\AA)\endsubhead
We apply Theorem~\TA\ to the cell arrangement which has row lengths
$n,n+1,\dots,n+k,n+k,\dots,n+k$ (from top to bottom), where $n+k$ occurs $n$~times.
Figure~\FD\ shows such a cell arrangement for $n=6$ and $k=3$.
(The reader should ignore the crosses and labellings at this point.)

\midinsert
\vskip10pt
$$
\Einheit.35cm
\Pfad(0,0),111111111111111111111111111\endPfad
\Pfad(0,3),111111111111111111111111111\endPfad
\Pfad(0,6),111111111111111111111111111\endPfad
\Pfad(0,9),111111111111111111111111111\endPfad
\Pfad(0,12),111111111111111111111111111\endPfad
\Pfad(0,15),111111111111111111111111111\endPfad
\Pfad(0,18),111111111111111111111111111\endPfad
\Pfad(0,21),111111111111111111111111\endPfad
\Pfad(0,24),111111111111111111111\endPfad
\Pfad(0,27),111111111111111111\endPfad
\Pfad(0,0),222222222222222222222222222\endPfad
\Pfad(3,0),222222222222222222222222222\endPfad
\Pfad(6,0),222222222222222222222222222\endPfad
\Pfad(9,0),222222222222222222222222222\endPfad
\Pfad(12,0),222222222222222222222222222\endPfad
\Pfad(15,0),222222222222222222222222222\endPfad
\Pfad(18,0),222222222222222222222222222\endPfad
\Pfad(21,0),222222222222222222222222\endPfad
\Pfad(24,0),222222222222222222222\endPfad
\Pfad(27,0),222222222222222222\endPfad
\PfadDicke{3pt}
\Pfad(0,18),111111111111111111111111111\endPfad
%
\Label\ro{\text {\seventeenpoint X}}(25,16)
\Label\ro{\text {\seventeenpoint X}}(22,13)
\Label\ro{\text {\seventeenpoint X}}(19,10)
\Label\ro{\text {\seventeenpoint X}}(16,7)
\Label\ro{\text {\seventeenpoint X}}(10,4)
\Label\ro{\text {\seventeenpoint X}}(1,1)
\Label\ro{\text {\seventeenpoint X}}(4,19)
\Label\ro{\text {\seventeenpoint X}}(7,25)
\Label\ro{\text {\seventeenpoint X}}(13,22)
\Label\u{\emptyset}(0,0)
\Label\u{\emptyset}(3,0)
\Label\u{\emptyset}(6,0)
\Label\u{\emptyset}(9,0)
\Label\u{\emptyset}(12,0)
\Label\u{\emptyset}(15,0)
\Label\u{\emptyset}(18,0)
\Label\u{\emptyset}(21,0)
\Label\u{\emptyset}(24,0)
\Label\u{\emptyset}(27,0)
\Label\ro{1}(27,3)
\Label\lo{1}(24,3)
\Label\lo{1}(21,3)
\Label\lo{1}(18,3)
\Label\lo{1}(15,3)
\Label\lo{1}(12,3)
\Label\lo{1}(9,3)
\Label\lo{1}(6,3)
\Label\lo{1}(3,3)
\Label\lo{\emptyset}(0,3)
\Label\ro{2}(27,6)
\Label\lo{2}(24,6)
\Label\lo{2}(21,6)
\Label\lo{2}(18,6)
\Label\lo{2}(15,6)
\Label\lo{2}(12,6)
\Label\lo{1}(9,6)
\Label\lo{1}(6,6)
\Label\lo{1}(3,6)
\Label\lo{\emptyset}(0,6)
\Label\ro{3}(27,9)
\Label\lo{3}(24,9)
\Label\lo{3}(21,9)
\Label\lo{3}(18,9)
\Label\lo{2}(15,9)
\Label\lo{2}(12,9)
\Label\lo{1}(9,9)
\Label\lo{1}(6,9)
\Label\lo{1}(3,9)
\Label\lo{\emptyset}(0,9)
\Label\ro{4}(27,12)
\Label\lo{4}(24,12)
\Label\lo{4}(21,12)
\Label\lo{3}(18,12)
\Label\lo{2}(15,12)
\Label\lo{2}(12,12)
\Label\lo{1}(9,12)
\Label\lo{1}(6,12)
\Label\lo{1}(3,12)
\Label\lo{\emptyset}(0,12)
\Label\ro{5}(27,15)
\Label\lo{5}(24,15)
\Label\lo{4}(21,15)
\Label\lo{3}(18,15)
\Label\lo{2}(15,15)
\Label\lo{2}(12,15)
\Label\lo{1}(9,15)
\Label\lo{1}(6,15)
\Label\lo{1}(3,15)
\Label\lo{\emptyset}(0,15)
\Label\ro{6}(27,18)
\Label\ro{5}(24,18)
\Label\lo{4}(21,18)
\Label\lo{3}(18,18)
\Label\lo{2}(15,18)
\Label\lo{2}(12,18)
\Label\lo{1}(9,18)
\Label\lo{1}(6,18)
\Label\lo{1}(3,18)
\Label\lo{\emptyset}(0,18)
\Label\ro{51\ }(24,21)
\Label\ro{41}(21,21)
\Label\lo{31\ }(18,21)
\Label\lo{21\ }(15,21)
\Label\lo{21\ }(12,21)
\Label\lo{2}(9,21)
\Label\lo{2}(6,21)
\Label\lo{1}(3,21)
\Label\lo{\emptyset}(0,21)
\Label\ro{42\ }(21,24)
\Label\ro{\ 32}(18,24)
\Label\lo{31\ }(15,24)
\Label\lo{21\ }(12,24)
\Label\lo{2}(9,24)
\Label\lo{2}(6,24)
\Label\lo{1}(3,24)
\Label\lo{\emptyset}(0,24)
\Label\ro{321\ }(18,27)
\Label\o{32\ }(15,27)
\Label\o{31\ }(12,27)
\Label\o{3}(9,27)
\Label\o{2}(6,27)
\Label\o{1}(3,27)
\Label\o{\emptyset}(0,27)
\hskip9.45cm
$$
\vskip10pt
\centerline{\eightpoint Growth diagram bijection for (\AA)}
\vskip6pt
\centerline{\smc Figure \FD}
\endinsert

We consider {\it exclusively} fillings of this cell arrangement
which have the following two properties:

\roster 
\item There is {\it exactly} one cross in each row and in each column.
\item In the last $n$ rows of the cell arrangement (this is the part
of the cell arrangement in which all row lengths are~$n+k$; in the
figure it is separated from the upper part by a thick line) the
crosses form a NE-chain.
\endroster

It should be observed that Properties~(1) and~(2) together imply
that the cross in the right-most column of the cell
arrangement must occur at the top of the column. See Figure~\FD\ for an
example of such a filling. 

By Theorems~\TA\ and~\TB, the forward growth diagram construction yields a bijection
between the above fillings and
sequences of partitions (read along the top-right boundary of the cell
arrangement) of the form
$$\gather
\emptyset=\la^0\subset \la^1\subset \dots\subset \la^n\kern8cm
\tag\BAa
\\
\supset \la^{n+1}\subset \la^{n+2}\supset \la^{n+3}\subset\dots
\supset \la^{n+2k-1}\subset \la^{n+2k}=(n)
\tag\BAb
\\
\kern7cm
\supset(n-1)\supset\dots
\supset(1)\supset\emptyset,
\tag\BAc
\endgather$$
where successive partitions in this sequence differ by exactly one square.
In other words, the images under the forward growth diagram construction
of the fillings satisfying Properties~(1)
and~(2) decompose into an increasing sequence from the empty partition
to $\la:=\la^n$ (the part in~(\BAa); it corresponds to a standard tableau of
shape~$\la=\la^n$), followed by a vacillating tableau of length~$k$ from $\la=\la^n$
to $(n)$ (the part in~(\BAb) together with $\la^n$ from the previous
line), followed by a completely determined decreasing sequence
from~$(n)$ to the empty partition (the part in~(\BAc)). In particular,
it is Property~(2) combined with Theorem~\TB\ which implies that the
last $n+1$ partitions are completely determined as indicated above.

Conversely, again using
Theorems~\TA\ and~\TB, by putting such a sequence along the top-right
boundary of the cell arrangement and applying the backward growth
diagram construction, one obtains a filling with Properties~(1) and~(2).

Hence, we see that the fillings satisfying Properties~(1) and~(2) are
in bijection with pairs $(T,V)$, where $T$ is a standard tableau
of some shape~$\la$ of size~$n$, and where $V$ is a vacillating
tableau of length~$k$ from~$\la$ to~$(n)$. Clearly, the number of these pairs
is exactly the sum on the right-hand side of~(\AA).

It remains to count the fillings of the above cell arrangements.
The key observation is that, once the crosses are filled in the first
$k$~rows of the cell arrangement (in Figure~\FD\ this is the part
above the thick line), then the filling is already completely
determined since, by Property~(1), there must be exactly one cross in
each row and in each column, and, by Property~(2), these crosses must
form a NE-chain. 

So, the remaining question is: how many ways are there to put crosses in
the first $k$~rows? There are $n$ possibilities to put a cross in the
first row (which has length~$n$). Once that cross has been placed,
since by Property~(1) there can be only one cross in a column, there
remain $n$~possibilities to place a cross in the second row (which has
length~$n+1$). Etc. Thus, in total we have $n^k$ possibilities to
place crosses in the first $k$~rows, explaining the left-hand side
of~(\AA).\quad \quad \qed

\medskip
The example in Figure~\FD\ illustrates this bijection. The filling
there gets mapped to the pair
$$
\left(\matrix 1&2&3\\4&5\\6\endmatrix,\
321\supset32\subset42\supset41\subset51\supset5\subset6\right).
$$

\medskip
As we pointed out in the introduction, our growth diagram bijection
does not seem to be related to the deletion-insertion bijection of
Halverson and Lewandowski in~\cite{\HaLeAA} in any simple way. Indeed,
Halverson and Lewandowski realise the left-hand side of~(\AA) combinatorially as
sequences $(i_1,i_2,\dots,i_k)\in\{1,2,\dots,n\}^k$, and
use these data for their deletion-insertion procedure. Obviously,
there are several ways to extract such a sequence from our
fillings. Since we agreed that the placement of crosses in the first
$k$~rows determines all other crosses uniquely by Properties~(1)
and~(2), the most straightforward way to read off such a sequence
would be to define $i_1$ to be the number of the column in which we find the
cross in the first row, to define $i_2$ to be the number of the column
in which we find the cross in the second row not counting column~$i_1$, etc. In the
example of Figure~\FD, we would read off the sequence $(3,4,2)$.
(The reader should observe that this is the reversal of the sequence
that is used as an example in Figure~3 of~\cite{\HaLeAA}; yet, the
images of the bijections here and in~\cite{\HaLeAA} seem completely
unrelated.)

The explanation for the unrelatedness of the two bijections lies
perhaps in the fact that Halverson and Lewandowski use jeu de taquin
to define their deletion. While there is a growth diagram description
of jeu de taquin (cf\. \cite{\StanBI, Appendix~1 of Chapter~7}), we do
not use it here, which seems to make the two bijections incomparable.

\medskip
Nevertheless, they share one property that is in the focus
of~\cite{\BeriAA}: a limiting property. In order to understand what is
meant by this, let us introduce the following notation: given a
partition $\la=(\la_1,\la_2,\dots,\la_\ell)$, define the truncation
$\la^*:=(\la_2,\dots,\la_\ell)$. Then we have the following result.

\proclaim{Theorem \TC}
Let $(i_1,i_2,\dots,i_k)$ be a given sequence in $\{1,2,\dots,n\}^k$.
Let 
$$
(n)=\la^n\supset \la^{n+1}\subset \la^{n+2}\supset \la^{n+3}\subset\dots
\supset \la^{n+2k-1}\subset \la^{n+2k}=(n)
$$
be the image of $(i_1,i_2,\dots,i_k)$ {\rm(}under the identification with
placements of crosses of cell arrangements that was explained just
above{\rm)} under the growth diagram bijection described in the proof
of\/~{\rm(\AA)} in this section. Then
$$
\emptyset=(\la^n)^*\supset (\la^{n+1})^*\subset (\la^{n+2})^*\supset
(\la^{n+3})^*\subset\dots
\supset (\la^{n+2k-1})^*\subset (\la^{n+2k})^*=\emptyset
\tag\BB$$
is always the same for $n\ge \max\{i_1,i_2+1\dots,i_k+k-1\}$.
\endproclaim
\remark{Remark}
The vacillating tableau (\BB) is called limiting vacillating tableau
in~\cite{\BeriAA}.  
\endremark
\demo{Proof of Theorem \TC}
This is completely obvious in view of Theorem~\TB. Once the NE-chain
in the lower part of the cell arrangement (the last $n$~rows) ``takes
over" (meaning that it is longer than~$k$), the first component of the
partitions read along the top-right boundary of the triangular part of the
cell arrangement will be equal to the length of that chain, while all
other components are smaller, and are determined --- again via
Theorem~\TB\ --- by the crosses in the first $k$~rows, which in their turn correspond
to the sequence $(i_1,i_2,\dots,i_k)$.\quad \quad \qed
\enddemo

\subhead 4. Bijective proof of (\AD)\endsubhead
The proof of (\AD) follows the one of~(\AA). There is one
notable difference, though. Since, now, along the staircase part of the cell
arrangement we want to read a vacillating tableau from~$\la$ to~$\mu$
(with $\mu$ not necessarily equal to~$(n)$), we must modify the
arrangement of crosses in the lower part of the cell arrangement.
To be precise, given a partition $\mu=(\mu_1,\mu_2,\dots,\mu_\ell)$ of
size~$n$, we consider {\it exclusively} fillings of the cell arrangement
as in the proof of~(\AA) which have the following two properties:

\roster 
\item"(1')" There is {\it exactly} one cross in each row and in each column.
\item"(2')" In rows $k+n,k+n-1,\dots,k+n-\mu_1+1$ (the bottom-most 
$\mu_1$~rows of\linebreak the cell arrangement), the crosses form a NE-chain;
in rows $k+n-\mu_1,\mathbreak k+n-\mu_1-1,\dots,k+n-\mu_1-\mu_2+1$,
the crosses form another NE-chain,
each of these crosses lying to the left of the right-most cross of the
previous\linebreak NE-chain; \dots;
in rows $k+n-\mu_1-\dots-\mu_{\ell-1},
k+n-\mu_1-\dots-\mu_{\ell-1}-1,\dots,\mathbreak
k+n-\mu_1-\dots-\mu_{\ell-1}-\mu_\ell+1=k+1$,
the crosses also form a NE-chain,
each of these crosses lying to the left of the right-most cross of the
previous NE-chain.
\endroster

See Figure~\FE\ for an
example of such a filling, where $n=6$, $k=3$, and $\mu=(3,2,1)$.

\midinsert
$$
\Einheit.35cm
\Pfad(0,0),111111111111111111111111111\endPfad
\Pfad(0,3),111111111111111111111111111\endPfad
\Pfad(0,6),111111111111111111111111111\endPfad
\Pfad(0,9),111111111111111111111111111\endPfad
\Pfad(0,12),111111111111111111111111111\endPfad
\Pfad(0,15),111111111111111111111111111\endPfad
\Pfad(0,18),111111111111111111111111111\endPfad
\Pfad(0,21),111111111
111111111111111\endPfad
\Pfad(0,24),111111111111111111111\endPfad
\Pfad(0,27),111111111111111111\endPfad
\Pfad(0,0),222222222222222222222222222\endPfad
\Pfad(3,0),222222222222222222222222222\endPfad
\Pfad(6,0),222222222222222222222222222\endPfad
\Pfad(9,0),222222222222222222222222222\endPfad
\Pfad(12,0),222222222222222222222222222\endPfad
\Pfad(15,0),222222222222222222222222222\endPfad
\Pfad(18,0),222222222222222222222222222\endPfad
\Pfad(21,0),222222222222222222222222\endPfad
\Pfad(24,0),222222222222222222222\endPfad
\Pfad(27,0),222222222222222222\endPfad
\PfadDicke{3pt}
\Pfad(0,18),111111111111111111111111111\endPfad
%
\Label\ro{\text {\seventeenpoint X}}(25,7)
\Label\ro{\text {\seventeenpoint X}}(22,4)
\Label\ro{\text {\seventeenpoint X}}(16,1)
\Label\ro{\text {\seventeenpoint X}}(19,13)
\Label\ro{\text {\seventeenpoint X}}(10,10)
\Label\ro{\text {\seventeenpoint X}}(1,16)
\Label\ro{\text {\seventeenpoint X}}(4,25)
\Label\ro{\text {\seventeenpoint X}}(7,22)
\Label\ro{\text {\seventeenpoint X}}(13,19)
\Label\u{\emptyset}(0,0)
\Label\u{\emptyset}(3,0)
\Label\u{\emptyset}(6,0)
\Label\u{\emptyset}(9,0)
\Label\u{\emptyset}(12,0)
\Label\u{\emptyset}(15,0)
\Label\u{\emptyset}(18,0)
\Label\u{\emptyset}(21,0)
\Label\u{\emptyset}(24,0)
\Label\u{\emptyset}(27,0)
\Label\ro{1}(27,3)
\Label\lo{1}(24,3)
\Label\lo{1}(21,3)
\Label\lo{1}(18,3)
\Label\lo{\emptyset}(15,3)
\Label\lo{\emptyset}(12,3)
\Label\lo{\emptyset}(9,3)
\Label\lo{\emptyset}(6,3)
\Label\lo{\emptyset}(3,3)
\Label\lo{\emptyset}(0,3)
\Label\ro{2}(27,6)
\Label\lo{2}(24,6)
\Label\lo{1}(21,6)
\Label\lo{1}(18,6)
\Label\lo{\emptyset}(15,6)
\Label\lo{\emptyset}(12,6)
\Label\lo{\emptyset}(9,6)
\Label\lo{\emptyset}(6,6)
\Label\lo{\emptyset}(3,6)
\Label\lo{\emptyset}(0,6)
\Label\ro{3}(27,9)
\Label\lo{2}(24,9)
\Label\lo{1}(21,9)
\Label\lo{1}(18,9)
\Label\lo{\emptyset}(15,9)
\Label\lo{\emptyset}(12,9)
\Label\lo{\emptyset}(9,9)
\Label\lo{\emptyset}(6,9)
\Label\lo{\emptyset}(3,9)
\Label\lo{\emptyset}(0,9)
\Label\ro{31}(27,12)
\Label\lo{21\ }(24,12)
\Label\lo{11\ }(21,12)
\Label\lo{11\ }(18,12)
\Label\lo{1}(15,12)
\Label\lo{1}(12,12)
\Label\lo{\emptyset}(9,12)
\Label\lo{\emptyset}(6,12)
\Label\lo{\emptyset}(3,12)
\Label\lo{\emptyset}(0,12)
\Label\ro{32}(27,15)
\Label\lo{22\ }(24,15)
\Label\lo{21\ }(21,15)
\Label\lo{11\ }(18,15)
\Label\lo{1}(15,15)
\Label\lo{1}(12,15)
\Label\lo{\emptyset}(9,15)
\Label\lo{\emptyset}(6,15)
\Label\lo{\emptyset}(3,15)
\Label\lo{\emptyset}(0,15)
\Label\ro{321}(27,18)
\Label\ro{\ \ 221}(24,18)
\Label\lo{211\ \ }(21,18)
\Label\lo{111\ \ }(18,18)
\Label\lo{11\ }(15,18)
\Label\lo{11\ }(12,18)
\Label\lo{1}(9,18)
\Label\lo{1}(6,18)
\Label\lo{1}(3,18)
\Label\lo{\emptyset}(0,18)
\Label\ro{222}(24,21)
\Label\ro{\ \ 221}(21,21)
\Label\lo{211\ \ }(18,21)
\Label\lo{21\ }(15,21)
\Label\lo{11\ }(12,21)
\Label\lo{1}(9,21)
\Label\lo{1}(6,21)
\Label\lo{1}(3,21)
\Label\lo{\emptyset}(0,21)
\Label\ro{2211}(21,24)
\Label\ro{\ \ 221}(18,24)
\Label\lo{22\ }(15,24)
\Label\lo{21\ }(12,24)
\Label\lo{2}(9,24)
\Label\lo{1}(6,24)
\Label\lo{1}(3,24)
\Label\lo{\emptyset}(0,24)
\Label\ro{2211}(18,27)
\Label\o{221}(15,27)
\Label\o{211}(12,27)
\Label\o{21}(9,27)
\Label\o{2}(6,27)
\Label\o{1}(3,27)
\Label\o{\emptyset}(0,27)
\hskip9.45cm
$$
\vskip10pt
\centerline{\eightpoint Growth diagram bijection for (\AD)}
\vskip6pt
\centerline{\smc Figure \FE}
\endinsert

Theorem~\TB\
implies that, if we apply the forward growth diagram construction,
then along the right boundary we read (from bottom to top) the partitions
$$\multline
\emptyset\subset(1)\subset\dots\subset(\mu_1)
\subset(\mu_1,1)\subset\dots\subset(\mu_1,\mu_2)\\
\subset\dots\subset
(\mu_1,\dots,\mu_{\ell-1},1),\dots,(\mu_1,\dots,\mu_{\ell-1},\mu_\ell)=\mu,
\endmultline\tag\CACA$$
that is, first the first component of the partitions grows up
to~$\mu_1$, then the second component grows up to~$\mu_2$, \dots, and
finally the last component grows up to~$\mu_\ell$, so that we end
in~$\mu$. See Figure~\FE, where $\mu=(3,2,1)$.

Consequently, one sees that, by the growth diagram construction,
fillings satisfying Properties~(1') and~(2') are in bijection with
sequences of partitions (read along the top-right boundary of the cell
arrangement) of the form
$$\gather
\emptyset=\la^0\subset \la^1\subset \dots\subset \la^n\kern8cm
\tag\CAa
\\
\supset \la^{n+1}\subset \la^{n+2}\supset \la^{n+3}\subset\dots
\supset \la^{n+2k-1}\subset \la^{n+2k}=\mu\kern2.5cm
\tag\CAb
\\
\supset(\mu_1,\dots,\mu_\ell-1)\supset\dots\supset(\mu_1,\dots,\mu_{\ell-1})\supset
\dots\supset(\mu_1)\supset\dots\supset(1)\supset\emptyset,
\tag\CAc
\endgather$$
where successive partitions in this sequence differ by exactly one square.
Here again, these sequences decompose into three parts,
namely an increasing sequence from the empty partition
to $\la:=\la^n$ (the part in~(\CAa), corresponding to a standard tableau of
shape~$\la=\la^n$), followed by a vacillating tableau of length~$k$ from $\la=\la^n$
to $\mu$ (the part in~(\CAb) together with $\la^n$ from the previous
line), followed by a completely determined decreasing sequence
from~$\mu$ to the empty partition (the part in~(\CAc)). 

Hence, we see that the fillings satisfying Properties~(1') and~(2') are
in bijection with pairs $(T,V)$, where $T$ is a standard tableau
of some shape~$\la$ of size~$n$, and where $V$ is a vacillating
tableau of length~$k$ from~$\la$ to~$\mu$. Clearly, the number of these pairs
is exactly the sum on the right-hand side of~(\AD).

In order to explain the left-hand side of~(\AD),
it remains to count the fillings satisfying Properties~(1') and~(2').
We follow the procedure of the proof of~(\AA) in the previous section.
That is, we first place crosses in the first
$k$~rows of the cell arrangement (in Figure~\FE\ this is the part
above the thick line). By the argument given in the previous section,
there are $n^k$~possibilities for placing these $k$~crosses.
However, in difference to the proof in the previous section, here
it is not true that, once the crosses in the first $k$~rows have been
placed, the remaining filling is uniquely determined.

We claim that there are $f^\mu$ possibilities to complete the filling
once the crosses in the first $k$~rows have been placed. Clearly, this
would explain the left-hand side of~(\AD), and thus complete the proof
of~(\AD). 

\midinsert
$$
\Einheit.35cm
\Pfad(0,0),111111111111111111\endPfad
\Pfad(0,3),111111111111111111\endPfad
\Pfad(0,6),111111111111111111\endPfad
\Pfad(0,9),111111111111111111\endPfad
\Pfad(0,12),111111111111111111\endPfad
\Pfad(0,15),111111111111111111\endPfad
\Pfad(0,18),111111111111111111\endPfad
\Pfad(0,0),222222222222222222\endPfad
\Pfad(3,0),222222222222222222\endPfad
\Pfad(6,0),222222222222222222\endPfad
\Pfad(9,0),222222222222222222\endPfad
\Pfad(12,0),222222222222222222\endPfad
\Pfad(15,0),222222222222222222\endPfad
\Pfad(18,0),222222222222222222\endPfad
\PfadDicke{3pt}
\Pfad(0,18),111111111111111111\endPfad
%
\Label\ro{\text {\seventeenpoint X}}(16,7)
\Label\ro{\text {\seventeenpoint X}}(13,4)
\Label\ro{\text {\seventeenpoint X}}(7,1)
\Label\ro{\text {\seventeenpoint X}}(10,13)
\Label\ro{\text {\seventeenpoint X}}(4,10)
\Label\ro{\text {\seventeenpoint X}}(1,16)
\Label\u{\emptyset}(0,0)
\Label\u{\emptyset}(3,0)
\Label\u{\emptyset}(6,0)
\Label\u{\emptyset}(9,0)
\Label\u{\emptyset}(12,0)
\Label\u{\emptyset}(15,0)
\Label\u{\emptyset}(18,0)
\Label\r{1\ }(18,3)
\Label\lo{1}(15,3)
\Label\lo{1}(12,3)
\Label\lo{1}(9,3)
\Label\lo{\emptyset}(6,3)
\Label\lo{\emptyset}(3,3)
\Label\lo{\emptyset}(0,3)
\Label\r{2\ }(18,6)
\Label\lo{2}(15,6)
\Label\lo{1}(12,6)
\Label\lo{1}(9,6)
\Label\lo{\emptyset}(6,6)
\Label\lo{\emptyset}(3,6)
\Label\lo{\emptyset}(0,6)
\Label\r{3\ }(18,9)
\Label\lo{2}(15,9)
\Label\lo{1}(12,9)
\Label\lo{1}(9,9)
\Label\lo{\emptyset}(6,9)
\Label\lo{\emptyset}(3,9)
\Label\lo{\emptyset}(0,9)
\Label\r{31}(18,12)
\Label\lo{21\ }(15,12)
\Label\lo{11\ }(12,12)
\Label\lo{11\ }(9,12)
\Label\lo{1}(6,12)
\Label\lo{\emptyset}(3,12)
\Label\lo{\emptyset}(0,12)
\Label\r{32}(18,15)
\Label\lo{22\ }(15,15)
\Label\lo{21\ }(12,15)
\Label\lo{11\ }(9,15)
\Label\lo{1}(6,15)
\Label\lo{\emptyset}(3,15)
\Label\lo{\emptyset}(0,15)
\Label\r{\ 321}(18,18)
\Label\o{221}(15,18)
\Label\o{211}(12,18)
\Label\o{111}(9,18)
\Label\o{11}(6,18)
\Label\o{1}(3,18)
\Label\lo{\emptyset}(0,18)
\hskip6.3cm
$$
\vskip10pt
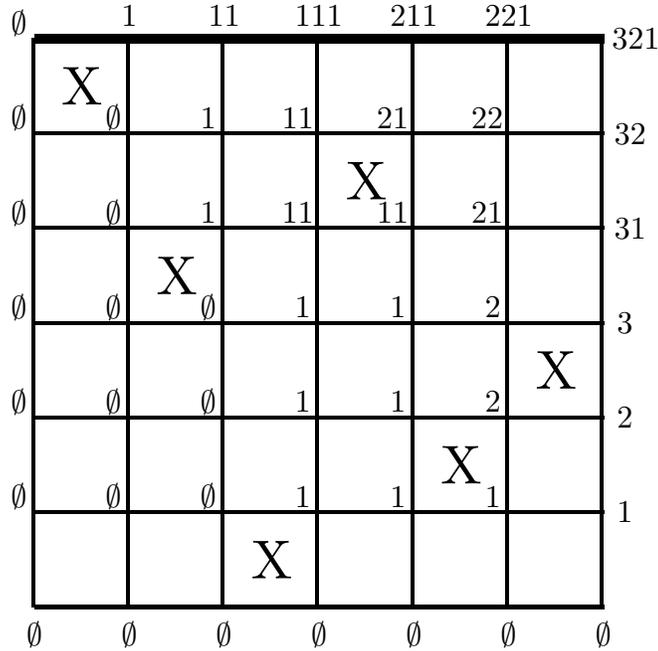
\centerline{\eightpoint Reduced growth diagram derived from
Figure~\FE}
\vskip6pt
\centerline{\smc Figure \FEa}
\endinsert

For establishing the claim, it is useful to consider the reduced growth
diagram where the first $k$~rows and the columns in which the crosses
of the first $k$~rows are located have been deleted. See Figure~\FEa\
for the reduced diagram that is obtained from the diagram in our
example in Figure~\FE. The crosses must be placed according to
Property~(2'). By Theorem~\TB, this is equivalent to the sequence of
partitions along the right boundary of the (square) growth diagram
being given 
by~(\CACA). On the other hand, along the top boundary we may have any
sequence of partitions $\emptyset\subset
\nu^1\subset\nu^2\subset\cdots\subset \nu^{n}=\mu$, and each such
sequence corresponds to a unique standard Young tableau of
shape~$\mu$. This proves 
the claim, and hence the theorem.\quad \quad \qed

\medskip
Figure~\FE, together with Figure~\FEa,
shows an example of the above bijection in which the pair
$$
\left((2,2,3)\ ,\ \matrix 1&4&6\\2&5\\3\endmatrix
\right)
$$
of a position sequence in $\{1,2,\dots,n\}^k=
\{1,2,\dots,6\}^3$ (corresponding to the vector of chosen
positions --- counted from the left --- in rows~1, 2, ‚\dots,~6,
in this order, of the growth diagram in Figure~\FE; it
should be kept in mind that, in the counting, we only take into
account  {\it available} positions)
and a standard Young tableau of shape $(3,2,1)$ (corresponding to the
sequence along the top boundary of the reduced growth diagram in Figure~\FEa)
gets mapped to the pair
$$
\left(\matrix 1&2&\\3&5\\4\\6\endmatrix,\
2211\supset221\subset2211\supset221\subset222\supset221\subset321\right).
$$

\subhead 5. Bijective proof of (\AC)\endsubhead
We choose again the cell arrangement from the proof of~(\AA), that is,
the arrangement with row lengths
$n,n+1,\dots,n+k,n+k,\dots,n+k$ (from top to bottom), where $n+k$ occurs $n$~times.
Figure~\FF\ shows such a cell arrangement for $n=3$ and $k=10$.

\midinsert
$$
\Einheit.3cm
\Pfad(0,0),111111111111111111111111111111111111111\endPfad
\Pfad(0,3),111111111111111111111111111111111111111\endPfad
\Pfad(0,6),111111111111111111111111111111111111111\endPfad
\Pfad(0,9),111111111111111111111111111111111111111\endPfad
\Pfad(0,12),111111111111111111111111111111111111\endPfad
\Pfad(0,15),111111111111111111111111111111111\endPfad
\Pfad(0,18),111111111111111111111111111111\endPfad
\Pfad(0,21),111111111111111111111111111\endPfad
\Pfad(0,24),111111111111111111111111\endPfad
\Pfad(0,27),111111111111111111111\endPfad
\Pfad(0,30),111111111111111111\endPfad
\Pfad(0,33),111111111111111\endPfad
\Pfad(0,36),111111111111\endPfad
\Pfad(0,39),111111111\endPfad
\Pfad(0,0),222222222222222222222222222222222222222\endPfad
\Pfad(3,0),222222222222222222222222222222222222222\endPfad
\Pfad(6,0),222222222222222222222222222222222222222\endPfad
\Pfad(9,0),222222222222222222222222222222222222222\endPfad
\Pfad(12,0),222222222222222222222222222222222222\endPfad
\Pfad(15,0),222222222222222222222222222222222\endPfad
\Pfad(18,0),222222222222222222222222222222\endPfad
\Pfad(21,0),222222222222222222222222222\endPfad
\Pfad(24,0),222222222222222222222222\endPfad
\Pfad(27,0),222222222222222222222\endPfad
\Pfad(30,0),222222222222222222\endPfad
\Pfad(33,0),222222222222222\endPfad
\Pfad(36,0),222222222222\endPfad
\Pfad(39,0),222222222\endPfad
\PfadDicke{3pt}
\Pfad(0,9),111111111111111111111111111111111111111\endPfad
\Pfad(9,0),222222222222222222222222222222222222222\endPfad
\Label\ro{\text {\seventeenpoint X}}(1,1)
\Label\ro{\text {\seventeenpoint X}}(4,22)
\Label\ro{\text {\seventeenpoint X}}(7,37)
\Label\ro{\text {\seventeenpoint X}}(10,34)
\Label\ro{\text {\seventeenpoint X}}(13,31)
\Label\ro{\text {\seventeenpoint X}}(16,28)
\Label\ro{\text {\seventeenpoint X}}(19,25)
\Label\ro{\text {\seventeenpoint X}}(22,19)
\Label\ro{\text {\seventeenpoint X}}(25,16)
\Label\ro{\text {\seventeenpoint X}}(28,4)
\Label\ro{\text {\seventeenpoint X}}(31,13)
\Label\ro{\text {\seventeenpoint X}}(34,10)
\Label\ro{\text {\seventeenpoint X}}(37,7)
\Label\u{\eightpoint\emptyset}(0,0)
\Label\u{\eightpoint\emptyset}(3,0)
\Label\u{\eightpoint\emptyset}(6,0)
\Label\u{\eightpoint\emptyset}(9,0)
\Label\u{\eightpoint\emptyset}(12,0)
\Label\u{\eightpoint\emptyset}(15,0)
\Label\u{\eightpoint\emptyset}(18,0)
\Label\u{\eightpoint\emptyset}(21,0)
\Label\u{\eightpoint\emptyset}(24,0)
\Label\u{\eightpoint\emptyset}(27,0)
\Label\u{\eightpoint\emptyset}(30,0)
\Label\u{\eightpoint\emptyset}(33,0)
\Label\u{\eightpoint\emptyset}(36,0)
\Label\u{\eightpoint\emptyset}(39,0)
\Label\ro{\eightpoint1}(39,3)
\Label\lo{\eightpoint1}(36,3)
\Label\lo{\eightpoint1}(33,3)
\Label\lo{\eightpoint1}(30,3)
\Label\lo{\eightpoint1}(27,3)
\Label\lo{\eightpoint1}(24,3)
\Label\lo{\eightpoint1}(21,3)
\Label\lo{\eightpoint1}(18,3)
\Label\lo{\eightpoint1}(15,3)
\Label\lo{\eightpoint1}(12,3)
\Label\lo{\eightpoint1}(9,3)
\Label\lo{\eightpoint1}(6,3)
\Label\lo{\eightpoint1}(3,3)
\Label\lo{\eightpoint\emptyset}(0,3)
\Label\ro{\eightpoint2}(39,6)
\Label\lo{\eightpoint2}(36,6)
\Label\lo{\eightpoint2}(33,6)
\Label\lo{\eightpoint2}(30,6)
\Label\lo{\eightpoint1}(27,6)
\Label\lo{\eightpoint1}(24,6)
\Label\lo{\eightpoint1}(21,6)
\Label\lo{\eightpoint1}(18,6)
\Label\lo{\eightpoint1}(15,6)
\Label\lo{\eightpoint1}(12,6)
\Label\lo{\eightpoint1}(9,6)
\Label\lo{\eightpoint1}(6,6)
\Label\lo{\eightpoint1}(3,6)
\Label\lo{\eightpoint\emptyset}(0,6)
\Label\ro{\eightpoint3}(39,9)
\Label\ro{\eightpoint2}(36,9)
\Label\lo{\eightpoint2}(33,9)
\Label\lo{\eightpoint2}(30,9)
\Label\lo{\eightpoint1}(27,9)
\Label\lo{\eightpoint1}(24,9)
\Label\lo{\eightpoint1}(21,9)
\Label\lo{\eightpoint1}(18,9)
\Label\lo{\eightpoint1}(15,9)
\Label\lo{\eightpoint1}(12,9)
\Label\lo{\eightpoint1}(9,9)
\Label\lo{\eightpoint1}(6,9)
\Label\lo{\eightpoint1}(3,9)
\Label\lo{\eightpoint\emptyset}(0,9)
\Label\ro{\eightpoint3}(36,12)
\Label\ro{\eightpoint2}(33,12)
\Label\lo{\eightpoint2}(30,12)
\Label\lo{\eightpoint1}(27,12)
\Label\lo{\eightpoint1}(24,12)
\Label\lo{\eightpoint1}(21,12)
\Label\lo{\eightpoint1}(18,12)
\Label\lo{\eightpoint1}(15,12)
\Label\lo{\eightpoint1}(12,12)
\Label\lo{\eightpoint1}(9,12)
\Label\lo{\eightpoint1}(6,12)
\Label\lo{\eightpoint1}(3,12)
\Label\lo{\eightpoint\emptyset}(0,12)
\Label\ro{\eightpoint3}(33,15)
\Label\ro{\eightpoint2}(30,15)
\Label\lo{\eightpoint1}(27,15)
\Label\lo{\eightpoint1}(24,15)
\Label\lo{\eightpoint1}(21,15)
\Label\lo{\eightpoint1}(18,15)
\Label\lo{\eightpoint1}(15,15)
\Label\lo{\eightpoint1}(12,15)
\Label\lo{\eightpoint1}(9,15)
\Label\lo{\eightpoint1}(6,15)
\Label\lo{\eightpoint1}(3,15)
\Label\lo{\eightpoint\emptyset}(0,15)
\Label\ro{\eightpoint21\ }(30,18)
\Label\ro{\eightpoint2}(27,18)
\Label\lo{\eightpoint1}(24,18)
\Label\lo{\eightpoint1}(21,18)
\Label\lo{\eightpoint1}(18,18)
\Label\lo{\eightpoint1}(15,18)
\Label\lo{\eightpoint1}(12,18)
\Label\lo{\eightpoint1}(9,18)
\Label\lo{\eightpoint1}(6,18)
\Label\lo{\eightpoint1}(3,18)
\Label\lo{\eightpoint\emptyset}(0,18)
\Label\ro{\eightpoint21\ }(27,21)
\Label\ro{\eightpoint2}(24,21)
\Label\lo{\eightpoint1}(21,21)
\Label\lo{\eightpoint1}(18,21)
\Label\lo{\eightpoint1}(15,21)
\Label\lo{\eightpoint1}(12,21)
\Label\lo{\eightpoint1}(9,21)
\Label\lo{\eightpoint1}(6,21)
\Label\lo{\eightpoint1}(3,21)
\Label\lo{\eightpoint\emptyset}(0,21)
\Label\ro{\eightpoint21\ }(24,24)
\Label\ro{\eightpoint2}(21,24)
\Label\lo{\eightpoint2}(18,24)
\Label\lo{\eightpoint2}(15,24)
\Label\lo{\eightpoint2}(12,24)
\Label\lo{\eightpoint2}(9,24)
\Label\lo{\eightpoint2}(6,24)
\Label\lo{\eightpoint1}(3,24)
\Label\lo{\eightpoint\emptyset}(0,24)
\Label\ro{\eightpoint3}(21,27)
\Label\ro{\eightpoint2}(18,27)
\Label\lo{\eightpoint2}(15,27)
\Label\lo{\eightpoint2}(12,27)
\Label\lo{\eightpoint2}(9,27)
\Label\lo{\eightpoint2}(6,27)
\Label\lo{\eightpoint1}(3,27)
\Label\lo{\eightpoint\emptyset}(0,27)
\Label\ro{\eightpoint3}(18,30)
\Label\ro{\eightpoint2}(15,30)
\Label\lo{\eightpoint2}(12,30)
\Label\lo{\eightpoint2}(9,30)
\Label\lo{\eightpoint2}(6,30)
\Label\lo{\eightpoint1}(3,30)
\Label\lo{\eightpoint\emptyset}(0,30)
\Label\ro{\eightpoint3}(15,33)
\Label\ro{\eightpoint2}(12,33)
\Label\lo{\eightpoint2}(9,33)
\Label\lo{\eightpoint2}(6,33)
\Label\lo{\eightpoint1}(3,33)
\Label\lo{\eightpoint\emptyset}(0,33)
\Label\ro{\eightpoint3}(12,36)
\Label\ro{\eightpoint2}(9,36)
\Label\lo{\eightpoint2}(6,36)
\Label\lo{\eightpoint1}(3,36)
\Label\lo{\eightpoint\emptyset}(0,36)
\Label\ro{\eightpoint3}(9,39)
\Label\o{\eightpoint2}(6,39)
\Label\o{\eightpoint1}(3,39)
\Label\o{\eightpoint\emptyset}(0,39)
\hskip12.3cm
$$
\vskip10pt
\centerline{\eightpoint Growth diagram bijection for (\AC)}
\vskip6pt
\centerline{\smc Figure \FF}
\endinsert

The fillings that we consider here are different, though. More
precisely, we consider {\it exclusively} fillings of this cell arrangement
which have the following three properties:

\roster 
\item There is {\it exactly} one cross in each row and in each column.
\item In the last $n$ rows of the cell arrangement (this is the part
of the cell arrangement in which all row lengths are~$n+k$; in the
figure it is separated from the upper part by a thick horizontal line) the
crosses form a NE-chain.
\item In the first $n$ columns of the cell arrangement (this is the part
of the cell arrangement in which all column lengths are~$n+k$; in the
figure it is separated from the right part by a thick vertical line) the
crosses form a NE-chain.
\endroster

It should be observed that Properties~(1) and~(2) together imply
that the cross in the right-most column of the cell
arrangement must occur at the top of the column, and that
Properties~(1) and~(3) imply that the cross in
the top row of the arrangement must be in the right-most cell of that
row. See Figure~\FF\ for an example of such a filling. 

By Theorems~\TA\ and~\TB, the forward growth diagram construction yields a bijection
between the above fillings and
sequences of partitions (read along the top-right boundary of the cell
arrangement) of the form
$$\gather
\emptyset\subset (1)\subset \dots\subset (n)=\la^n\kern8cm
\tag\DAa
\\
\supset \la^{n+1}\subset \la^{n+2}\supset \la^{n+3}\subset\dots
\supset \la^{n+2k-1}\subset \la^{n+2k}=(n)
\tag\DAb
\\
\kern7cm
\supset(n-1)\supset\dots
\supset(1)\supset\emptyset,
\tag\DAc
\endgather$$
where successive partitions in this sequence differ by exactly one square.
In other words, the images under the forward growth diagram construction
of the fillings satisfying Properties~(1)--(3) decompose into a
completely determined
increasing sequence from the empty partition 
to $(n)=\la^n$ (the part in~(\DAa)),
followed by a vacillating tableau of length~$k$ from $(n)$
to $(n)$ (the part in~(\DAb) together with $(n)$ from the previous
line), followed by a completely determined decreasing sequence
from~$(n)$ to the empty partition (the part in~(\DAc)). In particular,
it is Properties~(2) and~(3) combined with Theorem~\TB\ which imply that the
first $n+1$ and the last $n+1$ partitions are completely determined as
indicated above.


By definition, the above sequences are counted by~$m^{(n)}_{(n)}(k)$.
It remains to determine the number of fillings satisfying
Properties~(1)--(3).

The first observation is that such a filling will have a certain
number, say~$n-l$ for some integer~$l$ with $1\le l\le n$, of successive
crosses along the main diagonal of the cell arrangement. More
precisely, there will be crosses in rows and columns~$i$,
$i=1,2,\dots,n-l$ (counted from bottom-left), for some 
integer~$l$ with $1\le l\le n$. Further $l$~crosses
will have to be placed in columns~$n-l+1,\dots,n-1,n$ such that,
together with the aforementioned $n-l$~crosses along the main diagonal
they form a NE-chain, and further $l$~crosses
will have to be placed in rows~$n-l+1,\dots,n-1,n$ such that,
together with the aforementioned $n-l$~crosses along the main diagonal
they form a NE-chain.
In our example in Figure~\FF, there is just $1=3-2$ cross in the
bottom-left of the main diagonal, so that $l=2$.

\midinsert
$$
\Einheit.3cm
\Pfad(0,0),111111111111111111111111111111\endPfad
\Pfad(0,3),111111111111111111111111111\endPfad
\Pfad(0,6),111111111111111111111111\endPfad
\Pfad(0,9),111111111111111111111\endPfad
\Pfad(0,12),111111111111111111\endPfad
\Pfad(0,15),111111111111111\endPfad
\Pfad(0,18),111111111111\endPfad
\Pfad(0,21),111111111\endPfad
\Pfad(0,24),111111\endPfad
\Pfad(0,27),111\endPfad
\Pfad(0,0),222222222222222222222222222222\endPfad
\Pfad(3,0),222222222222222222222222222\endPfad
\Pfad(6,0),222222222222222222222222\endPfad
\Pfad(9,0),222222222222222222222\endPfad
\Pfad(12,0),222222222222222222\endPfad
\Pfad(15,0),222222222222222\endPfad
\Pfad(18,0),222222222222\endPfad
\Pfad(21,0),222222222\endPfad
\Pfad(24,0),222222\endPfad
\Pfad(27,0),222\endPfad
\PfadDicke{3pt}
\Pfad(0,0),111111111111111111111111111111\endPfad
\Pfad(0,0),222222222222222222222222222222\endPfad
%
\Label\ro{\text {\seventeenpoint X}}(1,25)
\Label\ro{\text {\seventeenpoint X}}(4,22)
\Label\ro{\text {\seventeenpoint X}}(7,19)
\Label\ro{\text {\seventeenpoint X}}(10,16)
\Label\ro{\text {\seventeenpoint X}}(13,10)
\Label\ro{\text {\seventeenpoint X}}(16,7)
\Label\ro{\text {\seventeenpoint X}}(22,4)
\Label\ro{\text {\seventeenpoint X}}(25,1)
\hskip9.8cm
$$
\vskip10pt
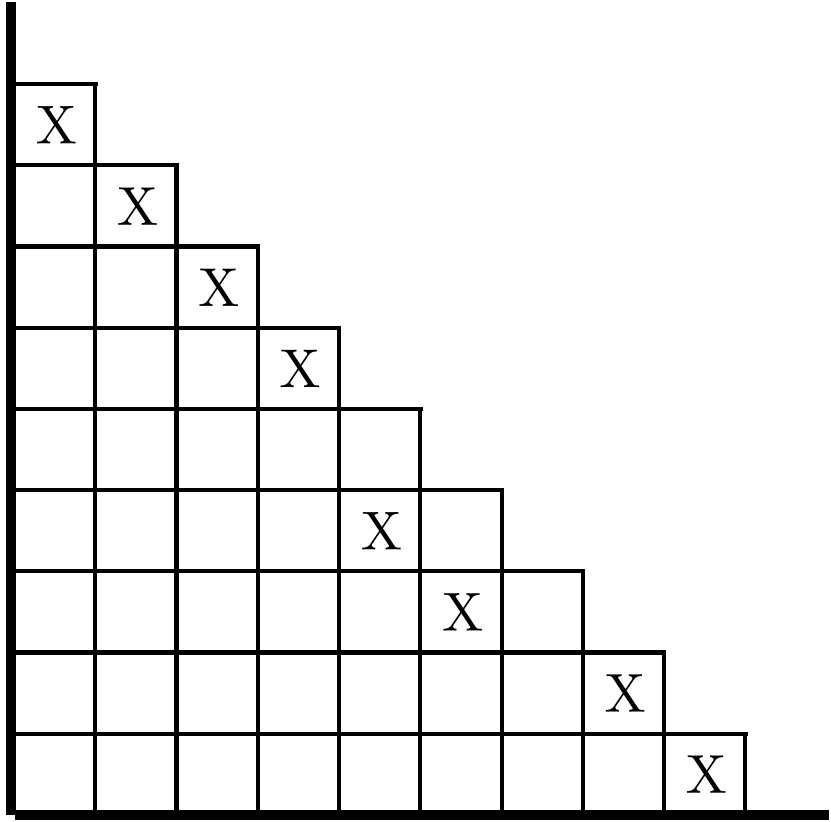
\centerline{\eightpoint The triangular region from Figure~\FF\ isolated}
\vskip6pt
\centerline{\smc Figure \FFa}
\endinsert

We now concentrate on the triangular region, $\De$ say, consisting of rows and
columns $n+1,n+2,\dots,n+k$ (again counted from bottom-left).
In Figure~\FF, this is the region to the right and above of the thick
lines. We have isolated that region in Figure~\FFa. The $l$~crosses in
columns $n-l+1,\dots,n-1,n$ discussed above occupy certain rows.
(One of the crosses must necessarily be placed at the end of the top
row of the cell arrangement.) These $l$~rows
must not be occupied by crosses in~$\De$. Similarly, the $l$~crosses in
rows $n-l+1,\dots,n-1,n$ discussed above occupy certain columns.
(One of the crosses must necessarily be placed at the top of the right-most
column of the cell arrangement.) These $l$~columns
must not be occupied by crosses in~$\De$. On the other hand, we must
place exactly one cross in each of the remaining $k-l$~rows and
$k-l$~columns. See Figures~\FF\ and~\FFa\ (where $n=3$, $k=10$, and~$l=2$).

The configuration of crosses in~$\De$ may be interpreted in a one-to-one fashion as a
set partition~$\pi$ of $\{1,2,\dots,k\}$. Indeed, if there is
a cross in column~$i$ of~$\De$ (counted from the left) and in row~$j$
of~$\De$ (counted from the top), then we declare $i$ and~$j$ to be
in the same block of~$\pi$. Thus, the configuration of crosses in
Figure~\FFa\ corresponds to the partition
$$
\{\{1,2,3,4,5,7\},\{6,8,9,10\}\}.
\tag\DB$$
This is indeed a one-to-one correspondence,
the inverse mapping being defined by ordering the
numbers in each block of the partition by size and placing a cross in
column~$i$ and row~$j$ whenever $i$ and~$j$ are successive elements in
a(n ordered) block (with $i<j$).

The proof of (\AC) can now be completed by observing that, under the
above described correspondence, configurations of $k-l$ crosses in~$\De$ are
in bijection with partitions of $\{1,2,\dots,k\}$ with $l$~blocks. This explains the
left-hand side of~(\AC) since the Stirling number $S(k,l)$ equals the
number of these partitions.\quad \quad \qed

\medskip
The bijection of the above proof  is illustrated in Figure~\FF,
mapping the partition in~(\DB) to the vacillating tableau
$$
3,2,3,2,3,2,3,2,3,2,21,2,21,2,21,2,3,2,3,2,3,
$$
and vice versa.

\subhead 6. Bijective proof of (\AE)\endsubhead
The proof of (\AE) is analogous to the one of~(\AC) in the previous
section. The only difference is that, here, Property~(3) of the
fillings gets replaced by

\roster 
\item "(3')" In the first $n$ columns of the cell arrangement (this is the part
of the cell arrangement in which all column lengths are~$n+k$) the
crosses form a SE-chain.
\endroster

As a consequence, there are two possibilities for the placement of the
crosses in the bottom-left square region
of the cell arrangement consisting of rows and columns
$1,2,\dots,n$ (counted from bottom-left): either there are no crosses in that region,
or there is exactly one cross, which is placed in row~1 and
column~$n$.

\midinsert
$$
\Einheit.3cm
\Pfad(0,0),111111111111111111111111111111111111111\endPfad
\Pfad(0,3),111111111111111111111111111111111111111\endPfad
\Pfad(0,6),111111111111111111111111111111111111111\endPfad
\Pfad(0,9),111111111111111111111111111111111111111\endPfad
\Pfad(0,12),111111111111111111111111111111111111\endPfad
\Pfad(0,15),111111111111111111111111111111111\endPfad
\Pfad(0,18),111111111111111111111111111111\endPfad
\Pfad(0,21),111111111111111111111111111\endPfad
\Pfad(0,24),111111111111111111111111\endPfad
\Pfad(0,27),111111111111111111111\endPfad
\Pfad(0,30),111111111111111111\endPfad
\Pfad(0,33),111111111111111\endPfad
\Pfad(0,36),111111111111\endPfad
\Pfad(0,39),111111111\endPfad
\Pfad(0,0),222222222222222222222222222222222222222\endPfad
\Pfad(3,0),222222222222222222222222222222222222222\endPfad
\Pfad(6,0),222222222222222222222222222222222222222\endPfad
\Pfad(9,0),222222222222222222222222222222222222222\endPfad
\Pfad(12,0),222222222222222222222222222222222222\endPfad
\Pfad(15,0),222222222222222222222222222222222\endPfad
\Pfad(18,0),222222222222222222222222222222\endPfad
\Pfad(21,0),222222222222222222222222222\endPfad
\Pfad(24,0),222222222222222222222222\endPfad
\Pfad(27,0),222222222222222222222\endPfad
\Pfad(30,0),222222222222222222\endPfad
\Pfad(33,0),222222222222222\endPfad
\Pfad(36,0),222222222222\endPfad
\Pfad(39,0),222222222\endPfad
\PfadDicke{3pt}
\Pfad(0,9),111111111111111111111111111111111111111\endPfad
\Pfad(9,0),222222222222222222222222222222222222222\endPfad
\Label\ro{\text {\seventeenpoint X}}(1,37)
\Label\ro{\text {\seventeenpoint X}}(4,31)
\Label\ro{\text {\seventeenpoint X}}(7,22)
\Label\ro{\text {\seventeenpoint X}}(10,34)
\Label\ro{\text {\seventeenpoint X}}(13,1)
\Label\ro{\text {\seventeenpoint X}}(16,28)
\Label\ro{\text {\seventeenpoint X}}(19,25)
\Label\ro{\text {\seventeenpoint X}}(22,19)
\Label\ro{\text {\seventeenpoint X}}(25,16)
\Label\ro{\text {\seventeenpoint X}}(28,4)
\Label\ro{\text {\seventeenpoint X}}(31,13)
\Label\ro{\text {\seventeenpoint X}}(34,10)
\Label\ro{\text {\seventeenpoint X}}(37,7)
\Label\u{\eightpoint\emptyset}(0,0)
\Label\u{\eightpoint\emptyset}(3,0)
\Label\u{\eightpoint\emptyset}(6,0)
\Label\u{\eightpoint\emptyset}(9,0)
\Label\u{\eightpoint\emptyset}(12,0)
\Label\u{\eightpoint\emptyset}(15,0)
\Label\u{\eightpoint\emptyset}(18,0)
\Label\u{\eightpoint\emptyset}(21,0)
\Label\u{\eightpoint\emptyset}(24,0)
\Label\u{\eightpoint\emptyset}(27,0)
\Label\u{\eightpoint\emptyset}(30,0)
\Label\u{\eightpoint\emptyset}(33,0)
\Label\u{\eightpoint\emptyset}(36,0)
\Label\u{\eightpoint\emptyset}(39,0)
\Label\ro{\eightpoint1}(39,3)
\Label\lo{\eightpoint1}(36,3)
\Label\lo{\eightpoint1}(33,3)
\Label\lo{\eightpoint1}(30,3)
\Label\lo{\eightpoint1}(27,3)
\Label\lo{\eightpoint1}(24,3)
\Label\lo{\eightpoint1}(21,3)
\Label\lo{\eightpoint1}(18,3)
\Label\lo{\eightpoint1}(15,3)
\Label\lo{\eightpoint\emptyset}(12,3)
\Label\lo{\eightpoint\emptyset}(9,3)
\Label\lo{\eightpoint\emptyset}(6,3)
\Label\lo{\eightpoint\emptyset}(3,3)
\Label\lo{\eightpoint\emptyset}(0,3)
\Label\ro{\eightpoint2}(39,6)
\Label\lo{\eightpoint2}(36,6)
\Label\lo{\eightpoint2}(33,6)
\Label\lo{\eightpoint2}(30,6)
\Label\lo{\eightpoint1}(27,6)
\Label\lo{\eightpoint1}(24,6)
\Label\lo{\eightpoint1}(21,6)
\Label\lo{\eightpoint1}(18,6)
\Label\lo{\eightpoint1}(15,6)
\Label\lo{\eightpoint\emptyset}(12,6)
\Label\lo{\eightpoint\emptyset}(9,6)
\Label\lo{\eightpoint\emptyset}(6,6)
\Label\lo{\eightpoint\emptyset}(3,6)
\Label\lo{\eightpoint\emptyset}(0,6)
\Label\ro{\eightpoint3}(39,9)
\Label\ro{\eightpoint2}(36,9)
\Label\lo{\eightpoint2}(33,9)
\Label\lo{\eightpoint2}(30,9)
\Label\lo{\eightpoint1}(27,9)
\Label\lo{\eightpoint1}(24,9)
\Label\lo{\eightpoint1}(21,9)
\Label\lo{\eightpoint1}(18,9)
\Label\lo{\eightpoint1}(15,9)
\Label\lo{\eightpoint\emptyset}(12,9)
\Label\lo{\eightpoint\emptyset}(9,9)
\Label\lo{\eightpoint\emptyset}(6,9)
\Label\lo{\eightpoint\emptyset}(3,9)
\Label\lo{\eightpoint\emptyset}(0,9)
\Label\ro{\eightpoint3}(36,12)
\Label\ro{\eightpoint2}(33,12)
\Label\lo{\eightpoint2}(30,12)
\Label\lo{\eightpoint1}(27,12)
\Label\lo{\eightpoint1}(24,12)
\Label\lo{\eightpoint1}(21,12)
\Label\lo{\eightpoint1}(18,12)
\Label\lo{\eightpoint1}(15,12)
\Label\lo{\eightpoint\emptyset}(12,12)
\Label\lo{\eightpoint\emptyset}(9,12)
\Label\lo{\eightpoint\emptyset}(6,12)
\Label\lo{\eightpoint\emptyset}(3,12)
\Label\lo{\eightpoint\emptyset}(0,12)
\Label\ro{\eightpoint3}(33,15)
\Label\ro{\eightpoint2}(30,15)
\Label\lo{\eightpoint1}(27,15)
\Label\lo{\eightpoint1}(24,15)
\Label\lo{\eightpoint1}(21,15)
\Label\lo{\eightpoint1}(18,15)
\Label\lo{\eightpoint1}(15,15)
\Label\lo{\eightpoint\emptyset}(12,15)
\Label\lo{\eightpoint\emptyset}(9,15)
\Label\lo{\eightpoint\emptyset}(6,15)
\Label\lo{\eightpoint\emptyset}(3,15)
\Label\lo{\eightpoint\emptyset}(0,15)
\Label\ro{\eightpoint21\ }(30,18)
\Label\ro{\eightpoint2}(27,18)
\Label\lo{\eightpoint1}(24,18)
\Label\lo{\eightpoint1}(21,18)
\Label\lo{\eightpoint1}(18,18)
\Label\lo{\eightpoint1}(15,18)
\Label\lo{\eightpoint\emptyset}(12,18)
\Label\lo{\eightpoint\emptyset}(9,18)
\Label\lo{\eightpoint\emptyset}(6,18)
\Label\lo{\eightpoint\emptyset}(3,18)
\Label\lo{\eightpoint\emptyset}(0,18)
\Label\ro{\eightpoint21\ }(27,21)
\Label\ro{\eightpoint2}(24,21)
\Label\lo{\eightpoint1}(21,21)
\Label\lo{\eightpoint1}(18,21)
\Label\lo{\eightpoint1}(15,21)
\Label\lo{\eightpoint\emptyset}(12,21)
\Label\lo{\eightpoint\emptyset}(9,21)
\Label\lo{\eightpoint\emptyset}(6,21)
\Label\lo{\eightpoint\emptyset}(3,21)
\Label\lo{\eightpoint\emptyset}(0,21)
\Label\ro{\eightpoint21\ }(24,24)
\Label\ro{\eightpoint11 }(21,24)
\Label\lo{\eightpoint11\ }(18,24)
\Label\lo{\eightpoint11\ }(15,24)
\Label\lo{\eightpoint1}(12,24)
\Label\lo{\eightpoint1}(9,24)
\Label\lo{\eightpoint\emptyset}(6,24)
\Label\lo{\eightpoint\emptyset}(3,24)
\Label\lo{\eightpoint\emptyset}(0,24)
\Label\ro{\eightpoint21\ }(21,27)
\Label\ro{\eightpoint11 }(18,27)
\Label\lo{\eightpoint11\ }(15,27)
\Label\lo{\eightpoint1}(12,27)
\Label\lo{\eightpoint1}(9,27)
\Label\lo{\eightpoint\emptyset}(6,27)
\Label\lo{\eightpoint\emptyset}(3,27)
\Label\lo{\eightpoint\emptyset}(0,27)
\Label\ro{\eightpoint21\ }(18,30)
\Label\ro{\eightpoint11}(15,30)
\Label\lo{\eightpoint1}(12,30)
\Label\lo{\eightpoint1}(9,30)
\Label\lo{\eightpoint\emptyset}(6,30)
\Label\lo{\eightpoint\emptyset}(3,30)
\Label\lo{\eightpoint\emptyset}(0,30)
\Label\ro{\eightpoint111\ }(15,33)
\Label\ro{\eightpoint11}(12,33)
\Label\lo{\eightpoint11\ }(9,33)
\Label\lo{\eightpoint1}(6,33)
\Label\lo{\eightpoint\emptyset}(3,33)
\Label\lo{\eightpoint\emptyset}(0,33)
\Label\ro{\eightpoint21\ }(12,36)
\Label\ro{\eightpoint11}(9,36)
\Label\lo{\eightpoint1}(6,36)
\Label\lo{\eightpoint\emptyset}(3,36)
\Label\lo{\eightpoint\emptyset}(0,36)
\Label\ro{\eightpoint111}(9,39)
\Label\o{\eightpoint11}(6,39)
\Label\o{\eightpoint1}(3,39)
\Label\o{\eightpoint\emptyset}(0,39)
\hskip12.3cm
$$
\vskip10pt
\centerline{\eightpoint Growth diagram bijection for (\AE), first case}
\vskip6pt
\centerline{\smc Figure \FG}
\endinsert

Figure~\FG\ provides an example for the former case, in which $n=3$
and $k=10$, and where the partition 
$$
\{\{1,2\},\{3,4,5,7\},\{6,8,9,10\}\}
$$
gets mapped to the vacillating tableau
$$\multline
111\supset11\subset21\supset11\subset111\supset11\subset21\supset11\subset
21\supset11\\
\subset21\supset2\subset21\supset2\subset21\supset2\subset
3\supset2\subset3\supset2\subset3,
\endmultline$$
while Figure~\FH\ shows an example for the latter case, in which again
$n=3$ and $k=10$, and where the partition 
$$
\{\{1,2,6,8,9,10\},\{3,4,5,7\}\}
$$
is mapped to
$$\multline
111\supset11\subset21\supset11\subset21\supset2\subset3\supset2
\subset3\supset2\subset21\\
\supset2\subset21\supset2\subset21\supset2\subset
3\supset2\subset3\supset2\subset3.
\endmultline$$
Clearly, here Theorem~\TB\ implies that the sequences that we read off
along the top-right boundary of the cell arrangement have the form
$$\gather
\emptyset\subset (1)\subset (1,1)\subset\dots\subset (1,1,\dots,1)=\la^n\kern5cm
\tag\EAa
\\
\supset \la^{n+1}\subset \la^{n+2}\supset \la^{n+3}\subset\dots
\supset \la^{n+2k-1}\subset \la^{n+2k}=(n)
\tag\EAb
\\
\kern7cm
\supset(n-1)\supset\dots
\supset(1)\supset\emptyset,
\tag\EAc
\endgather$$

\midinsert
$$
\Einheit.3cm
\Pfad(0,0),111111111111111111111111111111111111111\endPfad
\Pfad(0,3),111111111111111111111111111111111111111\endPfad
\Pfad(0,6),111111111111111111111111111111111111111\endPfad
\Pfad(0,9),111111111111111111111111111111111111111\endPfad
\Pfad(0,12),111111111111111111111111111111111111\endPfad
\Pfad(0,15),111111111111111111111111111111111\endPfad
\Pfad(0,18),111111111111111111111111111111\endPfad
\Pfad(0,21),111111111111111111111111111\endPfad
\Pfad(0,24),111111111111111111111111\endPfad
\Pfad(0,27),111111111111111111111\endPfad
\Pfad(0,30),111111111111111111\endPfad
\Pfad(0,33),111111111111111\endPfad
\Pfad(0,36),111111111111\endPfad
\Pfad(0,39),111111111\endPfad
\Pfad(0,0),222222222222222222222222222222222222222\endPfad
\Pfad(3,0),222222222222222222222222222222222222222\endPfad
\Pfad(6,0),222222222222222222222222222222222222222\endPfad
\Pfad(9,0),222222222222222222222222222222222222222\endPfad
\Pfad(12,0),222222222222222222222222222222222222\endPfad
\Pfad(15,0),222222222222222222222222222222222\endPfad
\Pfad(18,0),222222222222222222222222222222\endPfad
\Pfad(21,0),222222222222222222222222222\endPfad
\Pfad(24,0),222222222222222222222222\endPfad
\Pfad(27,0),222222222222222222222\endPfad
\Pfad(30,0),222222222222222222\endPfad
\Pfad(33,0),222222222222222\endPfad
\Pfad(36,0),222222222222\endPfad
\Pfad(39,0),222222222\endPfad
\PfadDicke{3pt}
\Pfad(0,9),111111111111111111111111111111111111111\endPfad
\Pfad(9,0),222222222222222222222222222222222222222\endPfad
\Label\ro{\text {\seventeenpoint X}}(1,37)
\Label\ro{\text {\seventeenpoint X}}(4,31)
\Label\ro{\text {\seventeenpoint X}}(7,1)
\Label\ro{\text {\seventeenpoint X}}(10,34)
\Label\ro{\text {\seventeenpoint X}}(13,22)
\Label\ro{\text {\seventeenpoint X}}(16,28)
\Label\ro{\text {\seventeenpoint X}}(19,25)
\Label\ro{\text {\seventeenpoint X}}(22,19)
\Label\ro{\text {\seventeenpoint X}}(25,16)
\Label\ro{\text {\seventeenpoint X}}(28,4)
\Label\ro{\text {\seventeenpoint X}}(31,13)
\Label\ro{\text {\seventeenpoint X}}(34,10)
\Label\ro{\text {\seventeenpoint X}}(37,7)
\Label\u{\eightpoint\emptyset}(0,0)
\Label\u{\eightpoint\emptyset}(3,0)
\Label\u{\eightpoint\emptyset}(6,0)
\Label\u{\eightpoint\emptyset}(9,0)
\Label\u{\eightpoint\emptyset}(12,0)
\Label\u{\eightpoint\emptyset}(15,0)
\Label\u{\eightpoint\emptyset}(18,0)
\Label\u{\eightpoint\emptyset}(21,0)
\Label\u{\eightpoint\emptyset}(24,0)
\Label\u{\eightpoint\emptyset}(27,0)
\Label\u{\eightpoint\emptyset}(30,0)
\Label\u{\eightpoint\emptyset}(33,0)
\Label\u{\eightpoint\emptyset}(36,0)
\Label\u{\eightpoint\emptyset}(39,0)
\Label\ro{\eightpoint1}(39,3)
\Label\lo{\eightpoint1}(36,3)
\Label\lo{\eightpoint1}(33,3)
\Label\lo{\eightpoint1}(30,3)
\Label\lo{\eightpoint1}(27,3)
\Label\lo{\eightpoint1}(24,3)
\Label\lo{\eightpoint1}(21,3)
\Label\lo{\eightpoint1}(18,3)
\Label\lo{\eightpoint1}(15,3)
\Label\lo{\eightpoint1}(12,3)
\Label\lo{\eightpoint1\ }(9,3)
\Label\lo{\eightpoint\emptyset}(6,3)
\Label\lo{\eightpoint\emptyset}(3,3)
\Label\lo{\eightpoint\emptyset}(0,3)
\Label\ro{\eightpoint2}(39,6)
\Label\lo{\eightpoint2}(36,6)
\Label\lo{\eightpoint2}(33,6)
\Label\lo{\eightpoint2}(30,6)
\Label\lo{\eightpoint1}(27,6)
\Label\lo{\eightpoint1}(24,6)
\Label\lo{\eightpoint1}(21,6)
\Label\lo{\eightpoint1}(18,6)
\Label\lo{\eightpoint1}(15,6)
\Label\lo{\eightpoint1}(12,6)
\Label\lo{\eightpoint1\ }(9,6)
\Label\lo{\eightpoint\emptyset}(6,6)
\Label\lo{\eightpoint\emptyset}(3,6)
\Label\lo{\eightpoint\emptyset}(0,6)
\Label\ro{\eightpoint3}(39,9)
\Label\ro{\eightpoint2}(36,9)
\Label\lo{\eightpoint2}(33,9)
\Label\lo{\eightpoint2}(30,9)
\Label\lo{\eightpoint1}(27,9)
\Label\lo{\eightpoint1}(24,9)
\Label\lo{\eightpoint1}(21,9)
\Label\lo{\eightpoint1}(18,9)
\Label\lo{\eightpoint1}(15,9)
\Label\lo{\eightpoint1}(12,9)
\Label\lo{\eightpoint1\ }(9,9)
\Label\lo{\eightpoint\emptyset}(6,9)
\Label\lo{\eightpoint\emptyset}(3,9)
\Label\lo{\eightpoint\emptyset}(0,9)
\Label\ro{\eightpoint3}(36,12)
\Label\ro{\eightpoint2}(33,12)
\Label\lo{\eightpoint2}(30,12)
\Label\lo{\eightpoint1}(27,12)
\Label\lo{\eightpoint1}(24,12)
\Label\lo{\eightpoint1}(21,12)
\Label\lo{\eightpoint1}(18,12)
\Label\lo{\eightpoint1}(15,12)
\Label\lo{\eightpoint1}(12,12)
\Label\lo{\eightpoint1\ }(9,12)
\Label\lo{\eightpoint\emptyset}(6,12)
\Label\lo{\eightpoint\emptyset}(3,12)
\Label\lo{\eightpoint\emptyset}(0,12)
\Label\ro{\eightpoint3}(33,15)
\Label\ro{\eightpoint2}(30,15)
\Label\lo{\eightpoint1}(27,15)
\Label\lo{\eightpoint1}(24,15)
\Label\lo{\eightpoint1}(21,15)
\Label\lo{\eightpoint1}(18,15)
\Label\lo{\eightpoint1}(15,15)
\Label\lo{\eightpoint1}(12,15)
\Label\lo{\eightpoint1\ }(9,15)
\Label\lo{\eightpoint\emptyset}(6,15)
\Label\lo{\eightpoint\emptyset}(3,15)
\Label\lo{\eightpoint\emptyset}(0,15)
\Label\ro{\eightpoint21\ }(30,18)
\Label\ro{\eightpoint2}(27,18)
\Label\lo{\eightpoint1}(24,18)
\Label\lo{\eightpoint1}(21,18)
\Label\lo{\eightpoint1}(18,18)
\Label\lo{\eightpoint1}(15,18)
\Label\lo{\eightpoint1}(12,18)
\Label\lo{\eightpoint1\ }(9,18)
\Label\lo{\eightpoint\emptyset}(6,18)
\Label\lo{\eightpoint\emptyset}(3,18)
\Label\lo{\eightpoint\emptyset}(0,18)
\Label\ro{\eightpoint21\ }(27,21)
\Label\ro{\eightpoint2}(24,21)
\Label\lo{\eightpoint1}(21,21)
\Label\lo{\eightpoint1}(18,21)
\Label\lo{\eightpoint1}(15,21)
\Label\lo{\eightpoint1}(12,21)
\Label\lo{\eightpoint1\ }(9,21)
\Label\lo{\eightpoint\emptyset}(6,21)
\Label\lo{\eightpoint\emptyset}(3,21)
\Label\lo{\eightpoint\emptyset}(0,21)
\Label\ro{\eightpoint21\ }(24,24)
\Label\ro{\eightpoint2}(21,24)
\Label\lo{\eightpoint2}(18,24)
\Label\lo{\eightpoint2}(15,24)
\Label\lo{\eightpoint1}(12,24)
\Label\lo{\eightpoint1\ }(9,24)
\Label\lo{\eightpoint\emptyset}(6,24)
\Label\lo{\eightpoint\emptyset}(3,24)
\Label\lo{\eightpoint\emptyset}(0,24)
\Label\ro{\eightpoint3}(21,27)
\Label\ro{\eightpoint2}(18,27)
\Label\lo{\eightpoint2}(15,27)
\Label\lo{\eightpoint1}(12,27)
\Label\lo{\eightpoint1\ }(9,27)
\Label\lo{\eightpoint\emptyset}(6,27)
\Label\lo{\eightpoint\emptyset}(3,27)
\Label\lo{\eightpoint\emptyset}(0,27)
\Label\ro{\eightpoint3}(18,30)
\Label\ro{\eightpoint2}(15,30)
\Label\lo{\eightpoint1}(12,30)
\Label\lo{\eightpoint1\ }(9,30)
\Label\lo{\eightpoint\emptyset}(6,30)
\Label\lo{\eightpoint\emptyset}(3,30)
\Label\lo{\eightpoint\emptyset}(0,30)
\Label\ro{\eightpoint21\ }(15,33)
\Label\ro{\eightpoint11}(12,33)
\Label\lo{\eightpoint11\ \ }(9,33)
\Label\lo{\eightpoint1}(6,33)
\Label\lo{\eightpoint\emptyset}(3,33)
\Label\lo{\eightpoint\emptyset}(0,33)
\Label\ro{\eightpoint21\ }(12,36)
\Label\ro{\eightpoint11}(9,36)
\Label\lo{\eightpoint1}(6,36)
\Label\lo{\eightpoint\emptyset}(3,36)
\Label\lo{\eightpoint\emptyset}(0,36)
\Label\ro{\eightpoint111}(9,39)
\Label\o{\eightpoint11}(6,39)
\Label\o{\eightpoint1}(3,39)
\Label\o{\eightpoint\emptyset}(0,39)
\hskip12.3cm
$$
\vskip10pt
\centerline{\eightpoint Growth diagram bijection for (\AE), second case}
\vskip6pt
\centerline{\smc Figure \FH}
\endinsert

We leave the details to the reader.\quad \quad \qed

\subhead 7. Sketch of bijective proof of (\AF)\endsubhead
We proceed as in the previous proofs of (\AC) and (\AE).
Here, we consider fillings of the cell arrangement with row lengths 
$n,n+1,\dots,n+k,n+k,\dots,n+k$ (from top to bottom) which have the
following properties:

\roster 
\item There is {\it exactly} one cross in each row and in each column.
\item The crosses in rows $k+1,k+2,\dots,k+n-1$ of the cell arrangement 
form a NE-chain, but the cross in row~$k+n$ (the last row) does not
extend this NE-chain.
\item The crosses in columns $2,3,\dots,n$ of the cell arrangement
form a NE-chain, but the cross in the first column does not extend
this NE-chain.
\endroster

\midinsert
$$
\Einheit.3cm
\Pfad(0,0),111111111111111111111111111\endPfad
\Pfad(0,3),111111111111111111111111111\endPfad
\Pfad(0,6),111111111111111111111111111\endPfad
\Pfad(0,9),111111111111111111111111111\endPfad
\Pfad(0,12),111111111111111111111111111\endPfad
\Pfad(0,15),111111111111111111111111111\endPfad
\Pfad(0,18),111111111111111111111111\endPfad
\Pfad(0,21),111111111111111111111\endPfad
\Pfad(0,24),111111111111111111\endPfad
\Pfad(0,27),111111111111111\endPfad
\Pfad(0,0),222222222222222222222222222\endPfad
\Pfad(3,0),222222222222222222222222222\endPfad
\Pfad(6,0),222222222222222222222222222\endPfad
\Pfad(9,0),222222222222222222222222222\endPfad
\Pfad(12,0),222222222222222222222222222\endPfad
\Pfad(15,0),222222222222222222222222222\endPfad
\Pfad(18,0),222222222222222222222222\endPfad
\Pfad(21,0),222222222222222222222\endPfad
\Pfad(24,0),222222222222222222\endPfad
\Pfad(27,0),222222222222222\endPfad
\PfadDicke{3pt}
\Pfad(0,15),111111111111111111111111111\endPfad
\Pfad(15,0),222222222222222222222222222\endPfad
\Label\ro{\text {\seventeenpoint X}}(1,4)
\Label\ro{\text {\seventeenpoint X}}(4,1)
\Label\ro{\text {\seventeenpoint X}}(7,7)
\Label\ro{\text {\seventeenpoint X}}(22,10)
\Label\ro{\text {\seventeenpoint X}}(25,13)
\Label\ro{\text {\seventeenpoint X}}(10,16)
\Label\ro{\text {\seventeenpoint X}}(13,25)
\Label\ro{\text {\seventeenpoint X}}(16,22)
\Label\ro{\text {\seventeenpoint X}}(19,19)
\Label\u{\eightpoint\emptyset}(0,0)
\Label\u{\eightpoint\emptyset}(3,0)
\Label\u{\eightpoint\emptyset}(6,0)
\Label\u{\eightpoint\emptyset}(9,0)
\Label\u{\eightpoint\emptyset}(12,0)
\Label\u{\eightpoint\emptyset}(15,0)
\Label\u{\eightpoint\emptyset}(18,0)
\Label\u{\eightpoint\emptyset}(21,0)
\Label\u{\eightpoint\emptyset}(24,0)
\Label\ru{\eightpoint\emptyset}(27,0)
\Label\ro{\eightpoint1}(27,3)
\Label\lo{\eightpoint1}(24,3)
\Label\lo{\eightpoint1}(21,3)
\Label\lo{\eightpoint1}(18,3)
\Label\lo{\eightpoint1\ }(15,3)
\Label\lo{\eightpoint1}(12,3)
\Label\lo{\eightpoint1}(9,3)
\Label\lo{\eightpoint1}(6,3)
\Label\lo{\eightpoint\emptyset}(3,3)
\Label\lo{\eightpoint\emptyset}(0,3)
\Label\ro{\eightpoint\ 11}(27,6)
\Label\lo{\eightpoint11\ }(24,6)
\Label\lo{\eightpoint11\ }(21,6)
\Label\lo{\eightpoint11\ }(18,6)
\Label\lo{\eightpoint11\ \ }(15,6)
\Label\lo{\eightpoint11\ }(12,6)
\Label\lo{\eightpoint11\ }(9,6)
\Label\lo{\eightpoint11\ }(6,6)
\Label\lo{\eightpoint1}(3,6)
\Label\lo{\eightpoint\emptyset}(0,6)
\Label\ro{\eightpoint\ 21}(27,9)
\Label\lo{\eightpoint21\ }(24,9)
\Label\lo{\eightpoint21\ }(21,9)
\Label\lo{\eightpoint21\ }(18,9)
\Label\lo{\eightpoint21\ \ }(15,9)
\Label\lo{\eightpoint21\ }(12,9)
\Label\lo{\eightpoint21\ }(9,9)
\Label\lo{\eightpoint11\ }(6,9)
\Label\lo{\eightpoint1}(3,9)
\Label\lo{\eightpoint\emptyset}(0,9)
\Label\ro{\eightpoint\ 31}(27,12)
\Label\lo{\eightpoint31\ }(24,12)
\Label\lo{\eightpoint21\ }(21,12)
\Label\lo{\eightpoint21\ }(18,12)
\Label\lo{\eightpoint21\ \ }(15,12)
\Label\lo{\eightpoint21\ }(12,12)
\Label\lo{\eightpoint21\ }(9,12)
\Label\lo{\eightpoint11\ }(6,12)
\Label\lo{\eightpoint1}(3,12)
\Label\lo{\eightpoint\emptyset}(0,12)
\Label\ro{\eightpoint\ 41}(27,15)
\Label\ro{\eightpoint31}(24,15)
\Label\lo{\eightpoint21\ }(21,15)
\Label\lo{\eightpoint21\ }(18,15)
\Label\lo{\eightpoint21\ \ }(15,15)
\Label\lo{\eightpoint21\ }(12,15)
\Label\lo{\eightpoint21\ }(9,15)
\Label\lo{\eightpoint11\ }(6,15)
\Label\lo{\eightpoint1}(3,15)
\Label\lo{\eightpoint\emptyset}(0,15)
\Label\ro{\eightpoint32}(24,18)
\Label\ro{\eightpoint31}(21,18)
\Label\lo{\eightpoint31\ }(18,18)
\Label\lo{\eightpoint31\ \ }(15,18)
\Label\lo{\eightpoint31\ }(12,18)
\Label\lo{\eightpoint21\ }(9,18)
\Label\lo{\eightpoint11\ }(6,18)
\Label\lo{\eightpoint1}(3,18)
\Label\lo{\eightpoint\emptyset}(0,18)
\Label\ro{\eightpoint41}(21,21)
\Label\ro{\eightpoint31}(18,21)
\Label\lo{\eightpoint31\ \ }(15,21)
\Label\lo{\eightpoint31\ }(12,21)
\Label\lo{\eightpoint21\ }(9,21)
\Label\lo{\eightpoint11\ }(6,21)
\Label\lo{\eightpoint1}(3,21)
\Label\lo{\eightpoint\emptyset}(0,21)
\Label\ro{\eightpoint41}(18,24)
\Label\ro{\eightpoint31}(15,24)
\Label\lo{\eightpoint31\ }(12,24)
\Label\lo{\eightpoint21\ }(9,24)
\Label\lo{\eightpoint11\ }(6,24)
\Label\lo{\eightpoint1}(3,24)
\Label\lo{\eightpoint\emptyset}(0,24)
\Label\ro{\eightpoint41}(15,27)
\Label\lo{\eightpoint31\ }(12,27)
\Label\lo{\eightpoint21\ }(9,27)
\Label\lo{\eightpoint11\ }(6,27)
\Label\lo{\eightpoint1}(3,27)
\Label\lo{\eightpoint\emptyset}(0,27)
\hskip8.1cm
$$
\vskip10pt
\centerline{\eightpoint Growth diagram bijection for (\AF), first case}
\vskip6pt
\centerline{\smc Figure \FI}
\endinsert

\midinsert
$$
\Einheit.3cm
\Pfad(0,0),111111111111111111111111111\endPfad
\Pfad(0,3),111111111111111111111111111\endPfad
\Pfad(0,6),111111111111111111111111111\endPfad
\Pfad(0,9),111111111111111111111111111\endPfad
\Pfad(0,12),111111111111111111111111111\endPfad
\Pfad(0,15),111111111111111111111111111\endPfad
\Pfad(0,18),111111111111111111111111\endPfad
\Pfad(0,21),111111111111111111111\endPfad
\Pfad(0,24),111111111111111111\endPfad
\Pfad(0,27),111111111111111\endPfad
\Pfad(0,0),222222222222222222222222222\endPfad
\Pfad(3,0),222222222222222222222222222\endPfad
\Pfad(6,0),222222222222222222222222222\endPfad
\Pfad(9,0),222222222222222222222222222\endPfad
\Pfad(12,0),222222222222222222222222222\endPfad
\Pfad(15,0),222222222222222222222222222\endPfad
\Pfad(18,0),222222222222222222222222\endPfad
\Pfad(21,0),222222222222222222222\endPfad
\Pfad(24,0),222222222222222222\endPfad
\Pfad(27,0),222222222222222\endPfad
\PfadDicke{3pt}
\Pfad(0,15),111111111111111111111111111\endPfad
\Pfad(15,0),222222222222222222222222222\endPfad
\Label\ro{\text {\seventeenpoint X}}(1,22)
\Label\ro{\text {\seventeenpoint X}}(25,1)
\Label\ro{\text {\seventeenpoint X}}(4,4)
\Label\ro{\text {\seventeenpoint X}}(19,10)
\Label\ro{\text {\seventeenpoint X}}(22,13)
\Label\ro{\text {\seventeenpoint X}}(7,7)
\Label\ro{\text {\seventeenpoint X}}(10,16)
\Label\ro{\text {\seventeenpoint X}}(13,25)
\Label\ro{\text {\seventeenpoint X}}(16,19)
\Label\u{\eightpoint\emptyset}(0,0)
\Label\u{\eightpoint\emptyset}(3,0)
\Label\u{\eightpoint\emptyset}(6,0)
\Label\u{\eightpoint\emptyset}(9,0)
\Label\u{\eightpoint\emptyset}(12,0)
\Label\u{\eightpoint\emptyset}(15,0)
\Label\u{\eightpoint\emptyset}(18,0)
\Label\u{\eightpoint\emptyset}(21,0)
\Label\u{\eightpoint\emptyset}(24,0)
\Label\ru{\eightpoint\emptyset}(27,0)
\Label\ro{\eightpoint1}(27,3)
\Label\lo{\eightpoint\emptyset}(24,3)
\Label\lo{\eightpoint\emptyset}(21,3)
\Label\lo{\eightpoint\emptyset}(18,3)
\Label\lo{\eightpoint\emptyset}(15,3)
\Label\lo{\eightpoint\emptyset}(12,3)
\Label\lo{\eightpoint\emptyset}(9,3)
\Label\lo{\eightpoint\emptyset}(6,3)
\Label\lo{\eightpoint\emptyset}(3,3)
\Label\lo{\eightpoint\emptyset}(0,3)
\Label\ro{\eightpoint\ 11}(27,6)
\Label\lo{\eightpoint1}(24,6)
\Label\lo{\eightpoint1}(21,6)
\Label\lo{\eightpoint1}(18,6)
\Label\lo{\eightpoint1\ }(15,6)
\Label\lo{\eightpoint1}(12,6)
\Label\lo{\eightpoint1}(9,6)
\Label\lo{\eightpoint1}(6,6)
\Label\lo{\eightpoint\emptyset}(3,6)
\Label\lo{\eightpoint\emptyset}(0,6)
\Label\ro{\eightpoint\ 21}(27,9)
\Label\lo{\eightpoint2}(24,9)
\Label\lo{\eightpoint2}(21,9)
\Label\lo{\eightpoint2}(18,9)
\Label\lo{\eightpoint2\ }(15,9)
\Label\lo{\eightpoint2}(12,9)
\Label\lo{\eightpoint2}(9,9)
\Label\lo{\eightpoint1}(6,9)
\Label\lo{\eightpoint\emptyset}(3,9)
\Label\lo{\eightpoint\emptyset}(0,9)
\Label\ro{\eightpoint\ 31}(27,12)
\Label\lo{\eightpoint3}(24,12)
\Label\lo{\eightpoint3}(21,12)
\Label\lo{\eightpoint2}(18,12)
\Label\lo{\eightpoint2\ }(15,12)
\Label\lo{\eightpoint2}(12,12)
\Label\lo{\eightpoint2}(9,12)
\Label\lo{\eightpoint1}(6,12)
\Label\lo{\eightpoint\emptyset}(3,12)
\Label\lo{\eightpoint\emptyset}(0,12)
\Label\ro{\eightpoint\ 41}(27,15)
\Label\ro{\eightpoint4}(24,15)
\Label\lo{\eightpoint3}(21,15)
\Label\lo{\eightpoint2}(18,15)
\Label\lo{\eightpoint2\ }(15,15)
\Label\lo{\eightpoint2}(12,15)
\Label\lo{\eightpoint2}(9,15)
\Label\lo{\eightpoint1}(6,15)
\Label\lo{\eightpoint\emptyset}(3,15)
\Label\lo{\eightpoint\emptyset}(0,15)
\Label\ro{\eightpoint41}(24,18)
\Label\ro{\eightpoint31}(21,18)
\Label\lo{\eightpoint3}(18,18)
\Label\lo{\eightpoint3\ }(15,18)
\Label\lo{\eightpoint3}(12,18)
\Label\lo{\eightpoint2}(9,18)
\Label\lo{\eightpoint1}(6,18)
\Label\lo{\eightpoint\emptyset}(3,18)
\Label\lo{\eightpoint\emptyset}(0,18)
\Label\ro{\eightpoint41}(21,21)
\Label\ro{\eightpoint4}(18,21)
\Label\lo{\eightpoint3\ }(15,21)
\Label\lo{\eightpoint3}(12,21)
\Label\lo{\eightpoint2}(9,21)
\Label\lo{\eightpoint1}(6,21)
\Label\lo{\eightpoint\emptyset}(3,21)
\Label\lo{\eightpoint\emptyset}(0,21)
\Label\ro{\eightpoint41}(18,24)
\Label\ro{\eightpoint31}(15,24)
\Label\lo{\eightpoint31\ }(12,24)
\Label\lo{\eightpoint21\ }(9,24)
\Label\lo{\eightpoint11\ }(6,24)
\Label\lo{\eightpoint1}(3,24)
\Label\lo{\eightpoint\emptyset}(0,24)
\Label\ro{\eightpoint41}(15,27)
\Label\lo{\eightpoint31\ }(12,27)
\Label\lo{\eightpoint21\ }(9,27)
\Label\lo{\eightpoint11\ }(6,27)
\Label\lo{\eightpoint1}(3,27)
\Label\lo{\eightpoint\emptyset}(0,27)
\hskip8.1cm
$$
\vskip10pt
\centerline{\eightpoint Growth diagram bijection for (\AF), second case}
\vskip6pt
\centerline{\smc Figure \FJ}
\endinsert

\midinsert
$$
\Einheit.3cm
\Pfad(0,0),111111111111111111111111111\endPfad
\Pfad(0,3),111111111111111111111111111\endPfad
\Pfad(0,6),111111111111111111111111111\endPfad
\Pfad(0,9),111111111111111111111111111\endPfad
\Pfad(0,12),111111111111111111111111\endPfad
\Pfad(0,15),111111111111111111111\endPfad
\Pfad(0,18),111111111111111111\endPfad
\Pfad(0,21),111111111111111\endPfad
\Pfad(0,24),111111111111\endPfad
\Pfad(0,27),111111111\endPfad
\Pfad(0,0),222222222222222222222222222\endPfad
\Pfad(3,0),222222222222222222222222222\endPfad
\Pfad(6,0),222222222222222222222222222\endPfad
\Pfad(9,0),222222222222222222222222222\endPfad
\Pfad(12,0),222222222222222222222222\endPfad
\Pfad(15,0),222222222222222222222\endPfad
\Pfad(18,0),222222222222222222\endPfad
\Pfad(21,0),222222222222222\endPfad
\Pfad(24,0),222222222222\endPfad
\Pfad(27,0),222222222\endPfad
\PfadDicke{3pt}
\Pfad(0,9),111111111111111111111111111\endPfad
\Pfad(9,0),222222222222222222222222222\endPfad
\Label\ro{\text {\seventeenpoint X}}(1,22)
\Label\ro{\text {\seventeenpoint X}}(25,1)
\Label\ro{\text {\seventeenpoint X}}(13,4)
\Label\ro{\text {\seventeenpoint X}}(19,10)
\Label\ro{\text {\seventeenpoint X}}(16,13)
\Label\ro{\text {\seventeenpoint X}}(22,7)
\Label\ro{\text {\seventeenpoint X}}(4,16)
\Label\ro{\text {\seventeenpoint X}}(7,25)
\Label\ro{\text {\seventeenpoint X}}(10,19)
\Label\u{\eightpoint\emptyset}(0,0)
\Label\u{\eightpoint\emptyset}(3,0)
\Label\u{\eightpoint\emptyset}(6,0)
\Label\u{\eightpoint\emptyset}(9,0)
\Label\u{\eightpoint\emptyset}(12,0)
\Label\u{\eightpoint\emptyset}(15,0)
\Label\u{\eightpoint\emptyset}(18,0)
\Label\u{\eightpoint\emptyset}(21,0)
\Label\u{\eightpoint\emptyset}(24,0)
\Label\ru{\eightpoint\emptyset}(27,0)
\Label\ro{\eightpoint1}(27,3)
\Label\lo{\eightpoint\emptyset}(24,3)
\Label\lo{\eightpoint\emptyset}(21,3)
\Label\lo{\eightpoint\emptyset}(18,3)
\Label\lo{\eightpoint\emptyset}(15,3)
\Label\lo{\eightpoint\emptyset}(12,3)
\Label\lo{\eightpoint\emptyset\ }(9,3)
\Label\lo{\eightpoint\emptyset}(6,3)
\Label\lo{\eightpoint\emptyset}(3,3)
\Label\lo{\eightpoint\emptyset}(0,3)
\Label\ro{\eightpoint\ 11}(27,6)
\Label\lo{\eightpoint1}(24,6)
\Label\lo{\eightpoint1}(21,6)
\Label\lo{\eightpoint1}(18,6)
\Label\lo{\eightpoint1\ }(15,6)
\Label\lo{\eightpoint\emptyset}(12,6)
\Label\lo{\eightpoint\emptyset\ }(9,6)
\Label\lo{\eightpoint\emptyset}(6,6)
\Label\lo{\eightpoint\emptyset}(3,6)
\Label\lo{\eightpoint\emptyset}(0,6)
\Label\ro{\eightpoint\ 21}(27,9)
\Label\ro{\eightpoint2}(24,9)
\Label\lo{\eightpoint1}(21,9)
\Label\lo{\eightpoint1}(18,9)
\Label\lo{\eightpoint1}(15,9)
\Label\lo{\eightpoint\emptyset}(12,9)
\Label\lo{\eightpoint\emptyset\ }(9,9)
\Label\lo{\eightpoint\emptyset}(6,9)
\Label\lo{\eightpoint\emptyset}(3,9)
\Label\lo{\eightpoint\emptyset}(0,9)
\Label\ro{\eightpoint21}(24,12)
\Label\ro{\eightpoint2}(21,12)
\Label\lo{\eightpoint1}(18,12)
\Label\lo{\eightpoint1}(15,12)
\Label\lo{\eightpoint\emptyset}(12,12)
\Label\lo{\eightpoint\emptyset\ }(9,12)
\Label\lo{\eightpoint\emptyset}(6,12)
\Label\lo{\eightpoint\emptyset}(3,12)
\Label\lo{\eightpoint\emptyset}(0,12)
\Label\ro{\eightpoint21}(21,15)
\Label\ro{\eightpoint2}(18,15)
\Label\lo{\eightpoint1}(15,15)
\Label\lo{\eightpoint\emptyset}(12,15)
\Label\lo{\eightpoint\emptyset\ }(9,15)
\Label\lo{\eightpoint\emptyset}(6,15)
\Label\lo{\eightpoint\emptyset}(3,15)
\Label\lo{\eightpoint\emptyset}(0,15)
\Label\ro{\eightpoint21}(18,18)
\Label\ro{\eightpoint11}(15,18)
\Label\lo{\eightpoint1\ }(12,18)
\Label\lo{\eightpoint1\ }(9,18)
\Label\lo{\eightpoint1}(6,18)
\Label\lo{\eightpoint\emptyset}(3,18)
\Label\lo{\eightpoint\emptyset}(0,18)
\Label\ro{\eightpoint21}(15,21)
\Label\ro{\eightpoint2}(12,21)
\Label\lo{\eightpoint1\ }(9,21)
\Label\lo{\eightpoint1}(6,21)
\Label\lo{\eightpoint\emptyset}(3,21)
\Label\lo{\eightpoint\emptyset}(0,21)
\Label\ro{\eightpoint21\ }(12,24)
\Label\ro{\eightpoint11}(9,24)
\Label\lo{\eightpoint11\ }(6,24)
\Label\lo{\eightpoint1}(3,24)
\Label\lo{\eightpoint\emptyset}(0,24)
\Label\ro{\eightpoint21\ }(9,27)
\Label\lo{\eightpoint11\ }(6,27)
\Label\lo{\eightpoint1}(3,27)
\Label\lo{\eightpoint\emptyset}(0,27)
\hskip8.1cm
$$
\vskip10pt
\centerline{\eightpoint Growth diagram bijection for (\AF), third case}
\vskip6pt
\centerline{\smc Figure \FK}
\endinsert

By Theorems~\TA\ and~\TB, the forward growth diagram construction yields a bijection
between the above fillings and
sequences of partitions (read along the top-right boundary of the cell
arrangement) of the form
$$\gather
\emptyset\subset
(1)\subset(1,1)\subset(2,1)\subset \dots\subset (n-1,1)=\la^n\kern4cm
\tag\FAa
\\
\supset \la^{n+1}\subset \la^{n+2}\supset \la^{n+3}\subset\dots
\supset \la^{n+2k-1}\subset \la^{n+2k}=(n-1,1)
\tag\FAb
\\
\kern4cm
\supset(n-2,1)\supset\dots\supset(2,1)\supset(1,1)
\supset(1)\supset\emptyset,
\tag\FAc
\endgather$$
where successive partitions in this sequence differ by exactly one square.
In other words, the images under the forward growth diagram construction
of the fillings satisfying Properties~(1)--(3) decompose into a
completely determined
increasing sequence from the empty partition 
to $(n-1,1)=\la^n$ (the part in~(\FAa)),
followed by a vacillating tableau of length~$k$ from $(n-1,1)$
to $(n-1,1)$ (the part in~(\FAb) together with $(n-1,1)$ from the previous
line), followed by a completely determined decreasing sequence
from~$(n-1,1)$ to the empty partition (the part in~(\FAc)). In particular,
it is Properties~(2) and~(3) combined with Theorem~\TB\ which imply that the
first $n+1$ and the last $n+1$ partitions are completely determined as
indicated above.


By definition, the above sequences are counted by~$m^{(n-1,1)}_{(n-1,1)}(k)$.
It remains to determine the number of fillings satisfying
Properties~(1)--(3).

The set of these fillings decomposes
into three pairwise disjoint subsets according to three structurally
different possibilities to place the crosses. These three
possibilities are exemplified in Figures~\FI--\FK, respectively. The
meaning of the thick lines in the figures is 
the same as in Figures~\FF--\FH.

The first possibility (see Figure~\FI, where $n=5$ and
$k=4$) is to have a cross in the
last row and the second column of the arrangement, a cross in
the next-to-last row and the first column, and crosses
along the main diagonal in rows and columns~$i$, for $i=3,4,\dots,n-l$
(counted from bottom-left), for some 
integer~$l$ with $1\le l\le n-2$, while further $l$~crosses
are placed in columns~$n-l+1,\dots,n-1,n$ such that,
together with the aforementioned $n-l-1$~crosses in columns $2,3,\dots,n-l$
they form a NE-chain, and further $l$~crosses
are placed in rows~$n-l+1,\dots,n-1,n$ such that,
together with the aforementioned $n-l-1$~crosses in rows
$2,3,\dots,n-l$ they form a NE-chain.
(In Figure~\FI, we have $l=2$.)
The number of these fillings is given by the first term on the
left-hand side of~(\AF).

The second possibility (see Figure~\FJ, where $n=5$ and $k=4$) is to have crosses
along the main diagonal in rows and columns~$i$, for $i=2,3,\dots,n-l+1$
(counted from bottom-left), for some 
integer~$l$ with $1\le l\le n-1$, while further $l$~crosses
are placed in columns~$n-l+2,\dots,n-1,n$ such that,
together with the aforementioned $n-l$~crosses in columns $2,3,\dots,n-l+1$
they form a NE-chain, further $l$~crosses
are placed in rows~$n-l+2,\dots,n-1,n$ such that,
together with the aforementioned $n-l$~crosses in rows
$2,3,\dots,n-l+1$ they form a NE-chain, and finally a cross is placed
in the first column within the first $k$~columns, and a cross is
placed in the last row within the last $k$~columns.
(In Figure~\FJ, we have $l=3$.)
The number of
these fillings is given by the second term on the left-hand side
of~(\AF). The multiplicative factor~$l^2$ has its explanation in the freedom to 
place the ``special crosses" in the first column and the last row in
relation to the two sets of ``further $l$~crosses" mentioned above.

The third (and last) possibility (see Figure~\FK, where $n=3$ and
$k=6$) is to place no
crosses into the bottom-left square region consisting of the first
$n$~columns and last $n$~rows, to place a NE-chain of crosses in
columns $2,3,\dots,n$ avoiding that square region, together with a cross in
the first column that does not extend this NE-chain, and to place a NE-chain of
crosses in rows $2,3,\dots,n$ (counted from bottom) avoiding that
square region, together with a cross in the last row that does not extend
this NE-chain. The number of
these fillings is given by the third term on the left-hand side
of~(\AF). Here, the multiplicative factor is $(n-1)^2$ (and not~$n^2$)
since there is one option less for the relative arrangement of the
crosses in the first $n$~columns and last $n$~rows.

We leave the details to the reader.\quad \quad \qed

\bigskip
\centerline{\bf Appendix A}
\medskip

Here we work out the special case of (\AC) where $n=3$ and $k=5$.
We have $S(5,1)=1$, $S(5,2)=15$, and $S(5,3)=25$. Hence, on the
left-hand side of~(\AC) we obtain $S(5,1)+S(5,2)+S(5,3)=41$. The 41
vacillating tableaux of length~5 from~$(3)$ to~$(3)$ that must exist according to
the right-hand side of~(\AC) are
$$\gather
3\supset 2\subset 3\supset 2\subset 3\supset 2\subset 3\supset 2\subset 3\supset  2\subset 3\\
3\supset 2\subset 3\supset 2\subset 3\supset 2\subset 3\supset 2\subset 21\supset 2\subset 3\\
3\supset 2\subset 3\supset 2\subset 3\supset 2\subset 21\supset 2\subset 3\supset 2\subset 3\\
3\supset 2\subset 3\supset 2\subset 3\supset 2\subset 21\supset 2\subset 21\supset 2\subset 3\\
3\supset 2\subset 3\supset 2\subset 3\supset 2\subset 21\supset 11\subset 21\supset 2\subset 3\\
3\supset 2\subset 3\supset 2\subset 21\supset 2\subset 3\supset 2\subset 3\supset  2\subset 3\\
3\supset 2\subset 3\supset 2\subset 21\supset 2\subset 3\supset 2\subset 21\supset 2\subset 3\\
3\supset 2\subset 3\supset 2\subset 21\supset 2\subset 21\supset 2\subset 3\supset 2\subset 3\\
3\supset 2\subset 3\supset 2\subset 21\supset 2\subset 21\supset 2\subset 21\supset 2\subset 3\\
3\supset 2\subset 3\supset 2\subset 21\supset 2\subset 21\supset 11\subset 21\supset 2\subset 3\\
3\supset 2\subset 3\supset 2\subset 21\supset 11\subset 21\supset 2\subset 3\supset 2\subset 3\\
3\supset 2\subset 3\supset 2\subset 21\supset 11\subset 21\supset 2\subset 21\supset 2\subset 3\\
3\supset 2\subset 3\supset 2\subset 21\supset 11\subset 21\supset 11\subset 21\supset 2\subset 3\\
3\supset 2\subset 21\supset 2\subset 3\supset 2\subset 3\supset 2\subset 3\supset  2\subset 3
\endgather$$
$$\gather
3\supset 2\subset 21\supset 2\subset 3\supset 2\subset 3\supset 2\subset 21\supset 2\subset 3\\
3\supset 2\subset 21\supset 2\subset 3\supset 2\subset 21\supset 2\subset 3\supset 2\subset 3\\
3\supset 2\subset 21\supset 2\subset 3\supset 2\subset 21\supset 2\subset 21\supset 2\subset 3\\
3\supset 2\subset 21\supset 2\subset 3\supset 2\subset 21\supset 11\subset 21\supset 2\subset 3\\
3\supset 2\subset 21\supset 2\subset 21\supset 2\subset 3\supset 2\subset 3\supset  2\subset 3\\
3\supset 2\subset 21\supset 2\subset 21\supset 2\subset 3\supset 2\subset 21\supset 2\subset 3\\
3\supset 2\subset 21\supset 2\subset 21\supset 2\subset 21\supset 2\subset 3\supset 2\subset 3\\
3\supset 2\subset 21\supset 2\subset 21\supset 2\subset 21\supset 2\subset 21\supset 2\subset 3\\
3\supset 2\subset 21\supset 2\subset 21\supset 2\subset 21\supset 11\subset 21\supset 2\subset 3\\
3\supset 2\subset 21\supset 2\subset 21\supset 11\subset 21\supset 2\subset 3\supset 2\subset 3\\
3\supset 2\subset 21\supset 2\subset 21\supset 11\subset 21\supset 2\subset 21\supset 2\subset 3\\
3\supset 2\subset 21\supset 2\subset 21\supset 11\subset 21\supset 11\subset 21\supset 2\subset 3\\
3\supset 2\subset 21\supset 11\subset 21\supset 2\subset 3\supset 2\subset 3\supset  2\subset 3\\
3\supset 2\subset 21\supset 11\subset 21\supset 2\subset 3\supset 2\subset 21\supset 2\subset 3\\
3\supset 2\subset 21\supset 11\subset 21\supset 2\subset 21\supset 2\subset 3\supset 2\subset 3\\
3\supset 2\subset 21\supset 11\subset 21\supset 2\subset 21\supset 2\subset 21\supset 2\subset 3\\
3\supset 2\subset 21\supset 11\subset 21\supset 2\subset 21\supset 11\subset 21\supset 2\subset 3\\
3\supset 2\subset 21\supset 11\subset 21\supset 11\subset 21\supset 2\subset 3\supset 2\subset 3\\
3\supset 2\subset 21\supset 11\subset 21\supset 11\subset 21\supset 2\subset 21\supset 2\subset 3\\
3\supset 2\subset 21\supset 11\subset 21\supset 11\subset 21\supset 11\subset 21\supset 2\subset 3\\
3\supset 2\subset 3\supset 2\subset 21\supset 11\subset 111\supset 11\subset 21\supset 2\subset 3\\
3\supset 2\subset 21\supset 2\subset 21\supset 11\subset 111\supset 11\subset 21\supset 2\subset 3\\
3\supset 2\subset 21\supset 11\subset 111\supset 11\subset 21\supset 2\subset 3\supset 2\subset 3\\
3\supset 2\subset 21\supset 11\subset 111\supset 11\subset 21\supset 2\subset 21\supset 2\subset 3\\
3\supset 2\subset 21\supset 11\subset 111\supset 11\subset 21\supset 11\subset 21\supset 2\subset 3\\
3\supset 2\subset 21\supset 11\subset 21\supset 11\subset 111\supset 11\subset 21\supset 2\subset 3\\
3\supset 2\subset 21\supset 11\subset 111\supset 11\subset 111\supset 11\subset 21\supset 2\subset 3
\endgather$$

\bigskip
\centerline{\bf Appendix B}
\medskip

Here we work out the special case of (\AE) where $n=3$ and $k=5$.
We have $S(5,2)=15$ and $S(5,3)=25$. Hence, on the
left-hand side of~(\AE) we obtain $S(5,2)+S(5,3)=40$. The 40
vacillating tableaux of length~5 from~$(3)$ to~$(1,1,1)$ that must exist according to
the right-hand side of~(\AE) are
$$\gather
3\supset 2\subset 3\supset 2\subset 3\supset 2\subset 3\supset 2\subset 21\supset 11\subset 111\\
3\supset 2\subset 3\supset 2\subset 3\supset 2\subset 21\supset 2\subset 21\supset 11\subset 111\\
3\supset 2\subset 3\supset 2\subset 3\supset 2\subset 21\supset 11\subset 21\supset 11\subset 111\\
\endgather$$
$$\gather
3\supset 2\subset 3\supset 2\subset 21\supset 2\subset 3\supset 2\subset 21\supset 11\subset 111\\
3\supset 2\subset 3\supset 2\subset 21\supset 2\subset 21\supset 2\subset 21\supset 11\subset 111\\
3\supset 2\subset 3\supset 2\subset 21\supset 2\subset 21\supset 11\subset 21\supset 11\subset 111\\
3\supset 2\subset 3\supset 2\subset 21\supset 11\subset 21\supset 2\subset 21\supset 11\subset 111\\
3\supset 2\subset 3\supset 2\subset 21\supset 11\subset 21\supset 11\subset 21\supset 11\subset 111\\
3\supset 2\subset 21\supset 2\subset 3\supset 2\subset 3\supset 2\subset 21\supset 11\subset 111\\
3\supset 2\subset 21\supset 2\subset 3\supset 2\subset 21\supset 2\subset 21\supset 11\subset 111\\
3\supset 2\subset 21\supset 2\subset 3\supset 2\subset 21\supset 11\subset 21\supset 11\subset 111\\
3\supset 2\subset 21\supset 2\subset 21\supset 2\subset 3\supset 2\subset 21\supset 11\subset 111\\
3\supset 2\subset 21\supset 2\subset 21\supset 2\subset 21\supset 2\subset 21\supset 11\subset 111\\
3\supset 2\subset 21\supset 2\subset 21\supset 2\subset 21\supset 11\subset 21\supset 11\subset 111\\
3\supset 2\subset 21\supset 2\subset 21\supset 11\subset 21\supset 2\subset 21\supset 11\subset 111\\
3\supset 2\subset 21\supset 2\subset 21\supset 11\subset 21\supset 11\subset 21\supset 11\subset 111\\
3\supset 2\subset 21\supset 11\subset 21\supset 2\subset 3\supset 2\subset 21\supset 11\subset 111\\
3\supset 2\subset 21\supset 11\subset 21\supset 2\subset 21\supset 2\subset 21\supset 11\subset 111\\
3\supset 2\subset 21\supset 11\subset 21\supset 2\subset 21\supset 11\subset 21\supset 11\subset 111\\
3\supset 2\subset 21\supset 11\subset 21\supset 11\subset 21\supset 2\subset 21\supset 11\subset 111\\
3\supset 2\subset 21\supset 11\subset 21\supset 11\subset 21\supset 11\subset 21\supset 11\subset 111\\
3\supset 2\subset 3\supset 2\subset 21\supset 11\subset 111\supset 11\subset 21\supset 11\subset 111\\
3\supset 2\subset 21\supset 2\subset 21\supset 11\subset 111\supset 11\subset 21\supset 11\subset 111\\
3\supset 2\subset 21\supset 11\subset 111\supset 11\subset 21\supset 2\subset 21\supset 11\subset 111\\
3\supset 2\subset 21\supset 11\subset 111\supset 11\subset 21\supset 11\subset 21\supset 11\subset 111\\
3\supset 2\subset 21\supset 11\subset 21\supset 11\subset 111\supset 11\subset 21\supset 11\subset 111\\
3\supset 2\subset 21\supset 11\subset 111\supset 11\subset 111\supset 11\subset 21\supset 11\subset 111\\
3\supset 2\subset 3\supset 2\subset 3\supset 2\subset 21\supset 11\subset 111\supset 11\subset 111\\
3\supset 2\subset 3\supset 2\subset 21\supset 2\subset 21\supset 11\subset 111\supset 11\subset 111\\
3\supset 2\subset 3\supset 2\subset 21\supset 11\subset 21\supset 11\subset 111\supset 11\subset 111\\
3\supset 2\subset 21\supset 2\subset 3\supset 2\subset 21\supset 11\subset 111\supset 11\subset 111\\
3\supset 2\subset 21\supset 2\subset 21\supset 2\subset 21\supset 11\subset 111\supset 11\subset 111\\
3\supset 2\subset 21\supset 2\subset 21\supset 11\subset 21\supset 11\subset 111\supset 11\subset 111\\
3\supset 2\subset 21\supset 11\subset 21\supset 2\subset 21\supset 11\subset 111\supset 11\subset 111\\
3\supset 2\subset 21\supset 11\subset 21\supset 11\subset 21\supset 11\subset 111\supset 11\subset 111\\
3\supset 2\subset 3\supset 2\subset 21\supset 11\subset 111\supset 11\subset 111\supset 11\subset 111\\
3\supset 2\subset 21\supset 2\subset 21\supset 11\subset 111\supset 11\subset 111\supset 11\subset 111\\
3\supset 2\subset 21\supset 11\subset 111\supset 11\subset 21\supset 11\subset 111\supset 11\subset 111\\
\endgather$$
$$\gather
3\supset 2\subset 21\supset 11\subset 21\supset 11\subset 111\supset 11\subset 111\supset 11\subset 111\\
3\supset 2\subset 21\supset 11\subset 111\supset 11\subset 111\supset 11\subset 111\supset 11\subset 111
\endgather$$

\bigskip
\centerline{\bf Appendix C}
\medskip

Here we work out the special case of (\AF) where $n=3$ and $k=3$.
We have $S(3,1)=1$, $S(3,2)=3$, and $S(3,3)=1$. Hence, on the
left-hand side of~(\AF) we obtain $S(3,1)+1^2\cdot S(3,1)+2^2\cdot S(3,2)
+(3-1)^2S(3,3)=18$. The 18
vacillating tableaux of length~3 from~$(2,1)$ to~$(2,1)$ that must exist according to
the right-hand side of~(\AF) are
$$\gather
21\supset 2\subset 21\supset 2\subset 21\supset 2\subset 21\\
21\supset 2\subset 21\supset 2\subset 21\supset 11\subset 21\\
21\supset 2\subset 21\supset 11\subset 21\supset 2\subset 21\\
21\supset 11\subset 21\supset 2\subset 21\supset 2\subset 21\\
21\supset 2\subset 21\supset 11\subset 21\supset 11\subset 21\\
21\supset 11\subset 21\supset 11\subset 21\supset 2\subset 21\\
21\supset 11\subset 21\supset 2\subset 21\supset 11\subset 21\\
21\supset 11\subset 21\supset 11\subset 21\supset 11\subset 21\\
X
21\supset 2\subset 3\supset 2\subset 21\supset 2\subset 21\\
21\supset 2\subset 21\supset 2\subset 3\supset 2\subset 21\\
21\supset 2\subset 3\supset 2\subset 3\supset 2\subset 21\\
Y
21\supset 2\subset 3\supset 2\subset 21\supset 11\subset 21\\
Z
21\supset 11\subset 21\supset 2\subset 3\supset 2\subset 21\\
U
21\supset 2\subset 21\supset 11\subset 111\supset 11\subset 21\\
21\supset 11\subset 111\supset 11\subset 21\supset 2\subset 21\\
21\supset 11\subset 111\supset 11\subset 21\supset 11\subset 21\\
21\supset 11\subset 21\supset 11\subset 111\supset 11\subset 21\\
21\supset 11\subset 111\supset 11\subset 111\supset 11\subset 21\\
\endgather$$


\Refs

\ref\no \BeriAA\by Z. Berikkyzy, P. E. Harris, A. Pun, C. Yan, and C. Zhao
\paper On the limiting vacillating tableaux for integer sequences
\jour preprint, {\tt ar$\chi$iv:2208.13091}\endref

\ref\no \BrFoAA\by T. Britz and S. Fomin\paper Finite posets and Ferrers shapes
\jour Adv\. Math\. \vol 158\yr 2001\pages 86--127\endref

\ref\no \ChDDAB\by W.Y.C. Chen, E. Y. P. Deng, R. R. X. Du,
R. P. Stanley and C. H. Yan \paper Crossings and nestings of
matchings and partitions\jour
Trans\. Amer\. Math\. Soc\. \vol359 \yr 2007\pages1555--1575\endref

\ref\no \FomiAZ\by S. V. Fomin\paper Generalized Robinson--Schensted--Knuth
correspondence\paperinfo (Russian)\jour
Zap.\ Nauchn.\ Sem.\ Leningrad.\ Otdel.\ Mat.\ Inst.\ 
Steklov.\ (LOMI) \vol 155\yr
1986\finalinfo translation in  J. Soviet Math.\ {\bf 41} (1988),
979--991\endref

\ref\no \FomiAB\by S.    Fomin\paper 
Schensted algorithms for graded graphs\jour
J. Alg.\ Combin\. \vol 4\yr 1995\pages 5--45\endref

\ref\no \FomiAF\by S.    Fomin\paper 
Schur operators and Knuth correspondences\jour
J.~Combin.\ Theory Ser.~A \vol 72\yr 1995\pages 277--292\endref

\ref\no \GreCAA\by C. Greene\paper An extension of Schensted's theorem
\jour Adv.\ Math\.\vol 14\yr 1974\pages 254--265\endref

\ref\no \HaLeAA\by T. Halverson and T. Lewandowski\paper
RSK insertion for set partitions and diagram algebras\jour
Electron\. J. Combin\. \vol 11 \yr 2004/06\pages  Research Paper~24, 24~pp
\endref

\ref\no \KratCE\by C.    Krattenthaler \yr 2006 \paper Growth
diagrams, and increasing and decreasing chains in fillings of Ferrers
shapes\jour Adv\. Appl\. Math\.\vol 37\pages 404-431\endref  

\ref\no \MaRoAA\by P. P. Martin and G. Rollet\paper The Potts model
representation and a Robinson--Schensted correspondence for the partition algebra
\jour Compositio Math\. \vol 112\yr 1998\pages 237--254\endref

\ref\no \RobyAA\by T. W. Roby\book Applications and extensions of 
Fomin's generalization 
of the Robinson--Schensted correspondence to differential posets\publ
Ph.D. thesis, M.I.T.\publaddr Cambridge, Massachusetts\yr 1991\endref

\ref\no \RobyAD\by T. W. Roby\paper 
The connection between the Robinson--Schensted
correspondence for skew oscillating tableaux and graded graphs\jour
Discrete Math\.\vol 139\yr 1995\pages 481--485\endref

\ref\no \SagaAQ\by B. E. Sagan\book The symmetric group
\bookinfo 2nd edition\publ Springer--Verlag\publaddr New York
\yr 2001\endref

\ref\no \StanBI\by R. P. Stanley \yr 1999 \book Enumerative
Combinatorics\bookinfo vol.~2\publ Cambridge University Press\publaddr
Cambridge\endref 

\endRefs
\enddocument